\documentclass[a4paper,12pt]{article}

\usepackage{amsmath,
            mleftright,
            comment,
            amssymb,
            amsthm,
            nicefrac,
            xcolor,
            enumerate,
            authblk,
            color,
            etoolbox,
            mathrsfs,
            mathtools,
            bbm,
            xparse,
            stmaryrd,
            enumitem,
            mathtools,
            marginnote,
            tikz,
            xurl,
            graphicx,
            scalerel,
            relsize,
            stackengine
            }

\usepackage[colorlinks=true]{hyperref}

\usepackage[sort,capitalize]{cleveref}

\crefformat{equation}{(#2#1#3)}
\crefname{subsection}{Subsection}{Subsections}
\crefname{enumi}{item}{items}
\crefname{equation}{}{}

\usepackage[left=20mm,right=18mm,top=1.5cm,bottom=1cm,includeheadfoot]{geometry}

\usepackage[T1]{fontenc}
\usepackage[utf8]{inputenc}
\usepackage[english]{babel}
\usepackage[sort,capitalize]{cleveref}

\usepackage{csquotes}

\newtheorem{theo}{Theorem}[section]

\newtheorem{lemma}[theo]{Lemma}
\newtheorem{propo}[theo]{Proposition}
\newtheorem{corollary}[theo]{Corollary}
\newtheorem{setting}[theo]{Setting}

\theoremstyle{definition}

\newtheorem{definition}[theo]{Definition}

\numberwithin{equation}{section}

\newcommand*{\ol}{\overline}

\newcommand*{\dm}{d}
\newcommand*{\fdm}{\mathfrak{d}}

\newcommand*{\w}{\mathfrak{w}}  
\newcommand*{\W}{\mathfrak{W}}

\newcommand*{\Borel}{\mathcal{B}}

\newenvironment{mproof}[1]{\noindent \textit{Proof of {#1}.}}{\hfill \qed}

\title{Deep neural networks can provably solve \\ Bellman equations for Markov decision processes \\ without the curse of dimensionality}

\author{Arnulf Jentzen$^{1,2}$,
Konrad Kleinberg$^{3}$, and Thomas Kruse$^{4}$\bigskip\\
\small{$^1$ School of Data Science and Shenzhen Research Institute of Big Data,}
\vspace{-0.1cm}\\
\small{The Chinese University of Hong Kong, Shenzhen (CUHK-Shenzhen),}\vspace{-0.1cm}\\
\small{China; e-mail: \texttt{ajentzen}\textcircled{\texttt{a}}\texttt{cuhk.edu.cn}}\smallskip\\
\small{$^2$ Applied Mathematics: Institute for Analysis and Numerics, }
\vspace{-0.1cm}\\ 
\small{University of M\"unster, Germany; e-mail: \texttt{ajentzen}\textcircled{\texttt{a}}\texttt{uni-muenster.de}}\smallskip\\
\small{$^3$ Department of Mathematics \& Informatics, University of Wuppertal,}
\vspace{-0.1cm}\\
\small{Germany; e-mail: \texttt{kleinberg}\textcircled{\texttt{a}}\texttt{uni-wuppertal.de}}\smallskip\\
\small{$^4$ Department of Mathematics \& Informatics, University of Wuppertal,}
\vspace{-0.1cm}\\
\small{Germany; e-mail: \texttt{tkruse}\textcircled{\texttt{a}}\texttt{uni-wuppertal.de}}
}

\DeclareFontEncoding{LS1}{}{}
\DeclareFontSubstitution{LS1}{stix}{m}{n}
\DeclareMathAlphabet{\mathscr}{LS1}{stixscr}{m}{n}

\newcommand*{\N}{\mathbb{N}}
\newcommand*{\R}{\mathbb{R}}
\newcommand*{\Z}{\mathbb{Z}}

\newcommand*{\E}{\mathbb{E}}
\newcommand*{\F}{\mathbb{F}}
\newcommand*{\X}{\mathbb{X}}
\newcommand*{\Y}{\mathbb{Y}}
\newcommand*{\I}{\mathbb{I}}

\newcommand*{\bF}{\mathbf{F}}
\newcommand*{\bG}{\mathbf{G}}

\newcommand*{\bI}{\mathbf{I}}

\newcommand*{\bM}{\mathbf{M}}
\newcommand*{\bN}{\mathbf{N} \cfadd{def:ANN}}

\newcommand*{\bQ}{\mathbf{Q}}

\newcommand*{\bX}{\mathbf{X}}

\newcommand{\cB}{\mathcal{B}}

\newcommand{\cD}{\mathcal{D} \cfadd{def:ANN}}

\newcommand{\cF}{\mathcal{F}}

\newcommand{\cH}{\mathcal{H}}
\newcommand{\cI}{\mathcal{I} \cfadd{def:ANN}}

\newcommand{\cL}{\mathcal{L} \cfadd{def:ANN}}

\newcommand{\cO}{\mathcal{O} \cfadd{def:ANN}}

\newcommand{\cR}{\mathcal{R} \cfadd{def:realization}}

\newcommand{\fJ}{{\bf\mathfrak{J}}}
\newcommand{\fK}{\mathfrak{K}}
\newcommand{\fL}{\mathfrak{L}}

\newcommand{\fS}{\mathfrak{S} \cfadd{def:sum_transpose_anns}}
\newcommand{\fT}{\mathfrak{T} \cfadd{def:sum_transpose_anns}}

\newcommand{\fW}{\mathfrak{W}}

\newcommand{\fc}{\mathfrak{c}}
\newcommand{\fd}{\mathfrak{d}}

\newcommand{\fm}{\mathfrak{m}}

\newcommand{\ft}{\mathfrak{t}}

\newcommand{\bfx}{\mathbf{x}}
\newcommand{\bfy}{\mathbf{y}}

\newcommand{\scrn}{\mathscr{n}}

\newcommand{\err}{(0,1]}

\newcommand{\induct}{\dashrightarrow}
\newcommand{\with}{\curvearrowleft}

\newcommand{\eps}{\varepsilon}

\newcommand{\dxx}{{\rm d}}

\NewDocumentCommand{\fnorm}{sO{}m}{%
  {\IfBooleanTF{#1}
    {\fnormhelp{\left|}{\right|}{#3}}
    {\fnormhelp{#2|}{#2|}{#3}}}
}
\makeatletter
\newcommand{\fnormhelp}[3]{\mathpalette\fnormhelp@i{{#1}{#2}{#3}}}
\newcommand{\fnormhelp@i}[2]{\fnormhelp@ii#1#2}
\newcommand{\fnormhelp@ii}[4]{%
  \sbox\z@{$\m@th#1#2#4#3$}%
  \sbox\tw@{$\m@th\|$}%
  \mathopen{\hbox to\wd\tw@{\hss\vrule height \ht\z@ depth \dp\z@ width .3\wd\tw@\hss}}%
  #4
  \mathclose{\hbox to\wd\tw@{\hss\vrule height \ht\z@ depth \dp\z@ width .3\wd\tw@\hss}}%
}

\newcommand{\tnorm}[1]{{ \left\vert\kern-0.25ex \left\vert\kern-0.25ex \left\vert #1 \right\vert\kern-0.25ex \right\vert\kern-0.25ex \right\vert } \cfadd{def:norms}}

\newcommand{\tnorms}[2]{ \mathopen{#1|\kern-0.25ex#1|\kern-0.25ex#1|}
#2
\mathclose{#1|\kern-0.25ex#1|\kern-0.25ex#1|}\cfadd{def:norms}}

\newcommand{\normmm}[1]{{\left\vert\kern-0.25ex\left\vert\kern-0.25ex\left\vert #1 
    \right\vert\kern-0.25ex\right\vert\kern-0.25ex\right\vert}
    \cfadd{DNN_norm}}

\newcommand{\tker}{\kappa \cfadd{def:stochastic_kernel}}
\newcommand{\tkerANN}{\rho}

\newcommand{\meas}{\mathcal{M} \cfadd{def:wasserstein_space}}
\newcommand{\coup}{\Pi \cfadd{def:couplings}}
\newcommand{\wdist}{\mathbb{W} \cfadd{def:wasserstein}}

\newcommand{\mlfpconst}{\gamma}
\newcommand{\mlfppower}{\alpha}

\DeclarePairedDelimiter{\pr}{(}{)}
\DeclarePairedDelimiter{\cu}{\{}{\}}
\DeclarePairedDelimiter{\br}{[}{]}

\DeclarePairedDelimiter{\abs}{\lvert}{\rvert}
\DeclarePairedDelimiter{\norm}{\lVert}{\rVert}
\DeclarePairedDelimiter{\ceil}{\lceil}{\rceil}

\newcommand{\lipf}{L}
\newcommand{\lipFANN}{\ell}

\newcommand{\mlfp}{\mathcal{V}}
\newcommand{\varmlfp}{\mathcal{U}}

\newcommand{\qfunc}{\mathscr{q}}

\newcommand{\action}{\mathscr{a}}

\newcommand{\genConst}{\mathfrak{c}}
\newcommand{\genConstVar}{\mathscr{c}}

\newcommand*{\id}{\operatorname{I} \cfadd{def:identity_matrix}}

\newcommand{\varint}{\textstyle\int}
\newcommand{\smallsum}{\textstyle\sum}
\newcommand{\SmallSum}[2]{ {\textstyle\sum\limits_{#1}^{#2}}}

\newcommand{\activation}{\phi}

\newcommand{\activationDim}[1]{\mathfrak{M}_{\activation,#1} \cfadd{def:multi}}

\newcommand{\paramANN}{\mathcal{P} \cfadd{def:ANN}}

\newcommand{\ANNcomp}{\bullet \cfadd{def:compositions_of_anns}}
\newcommand{\compANN}[2]{{#1 \bullet \allowbreak #2} \cfadd{def:compositions_of_anns}}

\newcommand{\paraANN}[1]{\mathbf{P}_{\!#1} \cfadd{def:parallelizations_of_anns_same_depth}}
\newcommand{\paraLANN}[2]{\mathrm{P}_{\!#1,#2} \cfadd{def:parallelizations_of_anns_with_different_layer_structure}}

\newcommand{\longerANN}[1]{\mathcal{E}_{#1} \cfadd{def:extensions_of_anns}}

\DeclareMathOperator{\idfunc}{id}

\newcommand{\dimANNlevel}{\mathbb{D} \cfadd{def:ANN}}
\newcommand{\power}[2]{#1^{\bullet #2} \cfadd{def:powers_of_anns}}
\newcommand{\scalarANN}{\circledast \cfadd{def:scalar_multiplications_of_anns}}
\newcommand{\scalar}[2]{ #1 \circledast #2 \cfadd{def:scalar_multiplications_of_anns}}

\newcommand{\oSum}{\oplus \cfadd{def:sum_of_anns_same_length}}
\newcommand{\OSum}[2]{{\mathop\oplus\limits_{#1}^{#2} } \cfadd{def:sum_of_anns_same_length} }
\newcommand{\bSum}{{\mathop\boxplus} \cfadd{def:sum_of_anns_diff_length}}
\newcommand{\BSum}[3]{{\mathop\boxplus\limits_{#1,#2}^{#3}} \cfadd{def:sum_of_anns_diff_length}}

\allowdisplaybreaks

\ExplSyntaxOn

\NewDocumentCommand{\enum}{ O{;} m o }
 {
  \my_enum:nnn { #1 } { #2 } { #3 }
 }

\seq_new:N \l__my_enum_seq
\tl_new:N \l__my_enum_item_tl
\prop_new:N \l__verbs
\prop_put:Nnn \l__verbs {show} {shows}
\prop_put:Nnn \l__verbs {imply} {implies}
\prop_put:Nnn \l__verbs {demonstrate} {demonstrates}
\prop_put:Nnn \l__verbs {prove} {proves}
\prop_put:Nnn \l__verbs {establish} {establishes}
\prop_put:Nnn \l__verbs {ensure} {ensures}
\prop_put:Nnn \l__verbs {assure} {assures}
\prop_put:Nnn \l__verbs {yield} {yields}

\int_new:N \l__number_of_args

\cs_new_protected:Nn \my_enum:nnn
 {
  \seq_set_split:Nnn \l__my_enum_seq { #1 } { #2 }
  \seq_remove_all:Nn \l__my_enum_seq {}
  \int_set_eq:NN \l__number_of_args { \seq_count:N \l__my_enum_seq }
  \seq_use:Nnnn \l__my_enum_seq { ~and~ } { ,~ } { ,~and~ }
  \IfNoValueTF{#3}{}{
    \space
    \int_compare:nNnTF{ \l__number_of_args } < {2}{ \prop_item:Nn \l__verbs {#3} }{ #3 }
  }
 }

\ExplSyntaxOff

\ExplSyntaxOn

\seq_new:N \g_cflist_loaded

\seq_new:N \g_cflist_pending

\NewDocumentCommand{\cfadd}{ m }

{

  \seq_if_in:NnF \g_cflist_loaded { #1 } {

    \seq_if_in:NnF \g_cflist_pending { #1 } {

      \seq_gput_right:Nn \g_cflist_pending { #1 }

    }

  }

}

\NewDocumentCommand{\cfload}{ o }

{

  \seq_if_empty:NTF \g_cflist_pending {\unskip} {

    (cf.\ \cref{\seq_use:Nn \g_cflist_pending {,}})\IfValueTF{#1}{#1~}{\unskip}

    \seq_gconcat:NNN \g_cflist_loaded \g_cflist_loaded \g_cflist_pending

    \seq_gclear:N \g_cflist_pending

  }

}

\NewDocumentCommand{\cfclear} {} {

  \seq_gclear:N \g_cflist_loaded

  \seq_gclear:N \g_cflist_pending

}

\NewDocumentCommand{\cfout}{ o }

{

  \seq_if_empty:NTF \g_cflist_pending {\unskip} {

    (cf.\ \cref{\seq_use:Nn \g_cflist_pending {,}})\IfValueTF{#1}{#1~}{\unskip}

    \seq_gclear:N \g_cflist_pending

  }

}

\NewDocumentCommand{\ifnocf} { m } {

  \seq_if_empty:NT \g_cflist_pending { #1 }

}

\NewDocumentCommand{\cfconsiderloaded}{ m }{

  \seq_gput_right:Nn \g_cflist_loaded {#1}

}

\ExplSyntaxOff

\NewDocumentEnvironment {athm} {m m} {%
\begin{#1}\label{#2}\global\def\loc{#2}%
}{%
\end{#1}%
}

\NewDocumentEnvironment{aproof} {} {%
\begin{proof}[Proof~of~\cref{\loc}]%
}{%
\finishproofthus
\end{proof}%
}

\ExplSyntaxOff

\newcommand{\finishproofthus}{The proof of \cref{\loc} is thus complete.}

\ExplSyntaxOn

\ExplSyntaxOff

\ExplSyntaxOn

\bool_new:N \g_noteobserve

\NewDocumentCommand{\setnote}{}{
  \bool_gset_true:N \g_noteobserve
}

\NewDocumentCommand{\setobserve}{}{
  \bool_gset_false:N \g_noteobserve
}

\NewDocumentCommand{\nobs}{ o }{
  \IfValueT{#1}{
    \str_if_eq:noTF {note} {#1} {
      \bool_gset_true:N \g_noteobserve
    } {
      \str_if_eq:noTF {Note} {#1} {
        \bool_gset_true:N \g_noteobserve
      } {
        \bool_gset_false:N \g_noteobserve
      }
    }
  }
  \bool_if:nTF { \g_noteobserve } {
    \bool_gset_false:N \g_noteobserve
    note
  } {
    \bool_gset_true:N \g_noteobserve
    observe
  }
  \IfValueF{#1}{~}
}

\NewDocumentCommand{\Nobs}{ o }{
  \IfValueT{#1}{
    \str_if_eq:noTF {note} {#1} {
      \bool_gset_true:N \g_noteobserve
    } {
      \str_if_eq:noTF {Note} {#1} {
        \bool_gset_true:N \g_noteobserve
      } {
        \bool_gset_false:N \g_noteobserve
      }
    }
  }
  \bool_if:nTF { \g_noteobserve } {
    \bool_gset_false:N \g_noteobserve
    Note
  } {
    \bool_gset_true:N \g_noteobserve
    Observe
  }
  \IfValueF{#1}{~}
}

\bool_new:N \g_hencetherefore

\NewDocumentCommand{\hence}{ o }{
  \IfValueT{#1}{
    \str_if_eq:noTF {hence} {#1} {
      \bool_gset_true:N \g_hencetherefore
    } {
      \str_if_eq:noTF {Hence} {#1} {
        \bool_gset_true:N \g_hencetherefore
      } {
        \bool_gset_false:N \g_hencetherefore
      }
    }
  }
  \bool_if:nTF { \g_hencetherefore } {
    \bool_gset_false:N \g_hencetherefore
    hence
  } {
    \bool_gset_true:N \g_hencetherefore
    therefore
  }
  \IfValueF{#1}{~}
}

\NewDocumentCommand{\Hence}{ o }{
  \IfValueT{#1}{
    \str_if_eq:noTF {hence} {#1} {
      \bool_gset_true:N \g_hencetherefore
    } {
      \str_if_eq:noTF {Hence} {#1} {
        \bool_gset_true:N \g_hencetherefore
      } {
        \bool_gset_false:N \g_hencetherefore
      }
    }
  }
  \bool_if:nTF { \g_hencetherefore } {
    \bool_gset_false:N \g_hencetherefore
    Hence
  } {
    \bool_gset_true:N \g_hencetherefore
    Therefore
  }
  \IfValueF{#1}{~}
}

\int_new:N \g_furthermore

\NewDocumentCommand{\Moreover}{ o o }{
  \IfValueT{#1}{
    \str_case:nn {#1} {
	  {Next} {\int_gset:Nn {\g_furthermore} {0}}      
      {Furthermore} {\int_gset:Nn {\g_furthermore} {1}}
      {Moreover} {\int_gset:Nn {\g_furthermore} {2}}
      {In~addition} {\int_gset:Nn {\g_furthermore} {3}}
      {note} {\bool_gset_true:N \g_noteobserve}
      {observe} {\bool_gset_false:N \g_noteobserve}
    }
    \IfValueT{#2}{
      \str_case:nn {#2} {
	    {Next} {\int_gset:Nn {\g_furthermore} {0}}        
        {Furthermore} {\int_gset:Nn {\g_furthermore} {1}}
        {Moreover} {\int_gset:Nn {\g_furthermore} {2}}
        {In~addition} {\int_gset:Nn {\g_furthermore} {3}}
        {note} {\bool_gset_true:N \g_noteobserve}
        {observe} {\bool_gset_false:N \g_noteobserve}
      }
    }
  }
  \int_case:nn { \int_mod:nn {\g_furthermore} {4} } {
	{ 0 } { Next,~\nobs that}    
    { 1 } { Furthermore,~\nobs that}
    { 2 } { Moreover,~\nobs that}
    { 3 } { In~addition,~\nobs that}
  }
  \int_gincr:N \g_furthermore
  \IfValueF{#1}{~}
}

\ExplSyntaxOff

\ExplSyntaxOn

\int_new:N \g_enumit

\NewDocumentCommand{\Enum}{ m o }{
  \IfValueT{#2}{
    \str_case:nn {#2} {
      {ensure} {\int_gset:Nn {\g_enumit} {0}}
      {ensures} {\int_gset:Nn {\g_enumit} {0}}       
      {assure} {\int_gset:Nn {\g_enumit} {1}}
      {assures} {\int_gset:Nn {\g_enumit} {1}}
      {prove} {\int_gset:Nn {\g_enumit} {2}}
      {proves} {\int_gset:Nn {\g_enumit} {2}}
      {show} {\int_gset:Nn {\g_enumit} {3}}
      {shows} {\int_gset:Nn {\g_enumit} {3}}
      {demonstrate} {\int_gset:Nn {\g_enumit} {4}}
      {demonstrates} {\int_gset:Nn {\g_enumit} {4}}
      {yield} {\int_gset:Nn {\g_enumit} {5}}
      {yields} {\int_gset:Nn {\g_enumit} {5}}
      {establish} {\int_gset:Nn {\g_enumit} {6}}
      {establishes} {\int_gset:Nn {\g_enumit} {6}}
      {imply} {\int_gset:Nn {\g_enumit} {7}}
      {implies} {\int_gset:Nn {\g_enumit} {7}}
    }
  }
  \int_case:nn { \int_mod:nn {\g_enumit} {8} } {
    { 0 } { \enum{ #1 }[ensure]{ } }
    { 1 } { \enum{ #1 }[assure] }  
    { 2 } { \enum{ #1 }[prove] }  
    { 3 } { \enum{ #1 }[show] }
    { 4 } { \enum{ #1 }[demonstrate] }
    { 5 } { \enum{ #1 }[yield] }
    { 6 } { \enum{ #1 }[establish] }
    { 7 } { \enum{ #1 }[imply] }
  }
  \int_gincr:N \g_enumit
} 

\ExplSyntaxOff

\begin{document}
\maketitle

\begin{abstract}
Discrete time stochastic optimal control problems and Markov decision processes (MDPs) are fundamental models for sequential decision-making under uncertainty and as such provide the mathematical framework underlying reinforcement learning theory. A central tool for solving MDPs is the Bellman equation and its solution, the so-called $Q$-function. In this article, we construct deep neural network (DNN) approximations for $Q$-functions associated to MDPs with infinite time horizon and finite control set $A$. More specifically, we show that if the the payoff function and the random transition dynamics of the MDP can be suitably approximated by DNNs with leaky rectified linear unit (ReLU) activation, then the solutions $Q_d\colon \R^d\to \R^{|A|}$, $d\in \mathbb{N}$, of the associated Bellman equations can also be approximated in the $L^2$-sense by DNNs with leaky ReLU activation whose numbers of parameters grow at most polynomially in both the dimension $d\in \mathbb{N}$ of the state space and the reciprocal $1/\varepsilon$ of the prescribed error $\varepsilon\in (0,1)$. Our proof relies on the recently introduced full-history recursive multilevel fixed-point (MLFP) approximation scheme.
\end{abstract}

\newpage

\tableofcontents

\section{Introduction}

Reinforcement Learning (RL) constitutes a powerful framework for modeling sequential decision-making with widespread applications in robotics, engineering, operations research, or economics (see, e.g., Sutton \& Barto \cite{sutton2018reinforcement}, Bertsekas \cite{bertsekas2019reinforcement}, and Bertsekas \& Tistsiklis \cite{bertsekas1996neuro} for comprehensive overviews). The mathematical underpinnings of RL are provided by the theory of stochastic optimal control and, in particular, Markov decision processes (MDPs; see, e.g., Bertsekas \& Shreve \cite{bertsekas1996stochastic}, Powell \cite{powell2007approximate}, and Puterman \cite{puterman2014markov}). At the core of many RL algorithms lies the so-called $Q$-function, which represents the expected reward of taking a given action in a given state and continuing optimally afterwards. The accurate estimation of the $Q$-function is pivotal for learning near-optimal policies and predicting maximal expected rewards.

However, approximating the $Q$-function in high-dimensional state spaces poses significant computational challenges. Many algorithms exhibit an exponential growth in the required computational resources to achieve a given accuracy as the dimension increases -- a phenomenon referred to as curse of dimensionality (cf., e.g., Bellman \cite{Bellman57}, Novak \& Wo\'zniakowski \cite[Chapter 1]{NovakWozniakowski2008I}, and Novak \& Ritter \cite{MR1485004}).

Empirically, the use of deep neural networks (DNNs) in reinforcement learning has demonstrated remarkable success across a wide range of high-dimensional tasks, including game playing, robotics, and control (see, e.g., \cite{mnih2015human}, \cite{silver2016mastering}, \cite{lillicrap2015continuous}, and \cite{kaufmann2023champion}).
Apart from this widespread empirical success there are several theoretical analyses that aim to explain the effectiveness of deep reinforcement learning (see, e.g.,  \cite{fan2020theoretical}, \cite{cai2019neural}, and \cite{xu2020finite}).
Yet, a comprehensive theoretical understanding of deep Q-learning is still incomplete.

This work contributes to closing this gap by establishing that, under suitable structural assumptions on the payoff function, the $Q$-function can be approximated by DNNs with a complexity that grows polynomially in the reciprocal $1/\eps$ of the prescribed approximation accuracy $\eps \in (0,1]$ and the dimension $d\in \N=\{1,2,3,\ldots\}$ of the state space. Our analysis builds on the recently developed full-history recursive multilevel Picard (MLP) approximations which have been proven to overcome the curse of dimensionality in the numerical approximation of certain semilinear partial differential equations (PDEs) (see, e.g., \cite{EHutzenthalerJentzenKruse2021, HJKNW2018, becker2020numerical, EHutzenthalerJentzenKruse2017, hutzenthaler2019overcoming, giles20019generalised, HJKN20, HutzenthalerKruse2017, hjk2019overcoming, beck2019overcoming, beck2024overcoming,beck2020nonlinear, hutzenthaler2020zovercoming, neufeld2025multilevel,neufeld2025bmultilevel}).
By combining these MLP schemes with compositional and approximation properties of DNNs it has been demonstrated that DNNs can overcome the curse of dimensionality in the numerical approximation of broad classes of PDEs (see, e.g., \cite{hutzenthaler2019proof, cioica2022deep, neufeld2023deep, ackermann2023deep,ackermann2024deep, neufeld2024rectified, neufeld2024multilevel}). The recent work \cite{beck2023nonlinctrl} introduced a variant of MLP approximations tailored to Bellman equations. In the present work, we use the results of \cite{beck2023nonlinctrl} to construct DNN approximations of the $Q$-function that avoid the curse of dimensionality.

To briefly sketch the contribution of this article within this introductory section, we now present in the following result, \cref{thm:intro} below, a special case of \cref{thm:main_simple}, the main result of this article. Below  \cref{thm:intro} we explain in words the statement of  \cref{thm:intro} as well as the mathematical objects appearing in \cref{thm:intro}.

\begin{samepage}
\begin{theo} \label{thm:intro}
    Let $\beta \in [0, \infty) \backslash \cu{1}$,
    for every $\dm \in \N$, $x \in \pr{x_1, \dots, x_\dm} \in \R^\dm$ let $\fm\pr{x} \in \R^\dm$ satisfy $\fm \pr{ x } =  \pr{ \max \cu{ x_1, \beta x_1 }, \dots, \max \cu{x_\dm, \beta x_\dm} } $,
    for every $\dm, \fdm \in \N$ let
    \begin{equation} \label{eq:anns_intro}
        \textstyle
        \bN_{\dm, \fdm} =
        \textstyle
        \cup_{H \in \N} \cup_{\pr{l_0, l_1, \dots, l_{H + 1}} \in \cu{\dm} \times \N^H \times \cu{\fdm} } \pr{\times_{k = 1}^{H + 1}   \pr{\R^{l_k \times l_{k - 1}} \times \R^{l_k}}},
    \end{equation}
    for every $L \in \N$, $l_0, l_1, \dots, l_L \in \N$, $\Phi = \pr{ \pr{W_1, B_1}, \dots, \pr{W_L, B_L} } \in \pr{ \times_{k = 1}^L \pr{ \R^{ l_{k} \times l_{k - 1} } \times \R^{l_k} } }$ let $\cD \pr{\Phi} \in \N^{L + 1}$, $\paramANN\pr{\Phi} \in \N$, and $\cR\pr{\Phi} \colon \R^{l_0} \to \R^{l_L}$ satisfy for all $p_0 \in \R^{l_0}$, $p_1 \in \R^{l_1}$, $\dots$, $p_{L} \in \R^{l_{L}}$ with $\forall \, k \in \cu{1,2,\dots, L} \colon p_k = \fm \pr{ W_k p_{k - 1} + B_k }$ that 
    \begin{equation}
        \cD \pr{\Phi} = \pr{l_0, l_1, \dots, l_L}, \quad \paramANN\pr{\Phi} = \textstyle \sum_{k = 1}^L l_k \pr{l_{k - 1} + 1}, \quad \text{and} \quad\pr{ \cR \pr{\Phi} } \pr{p_0} = W_L p_{L - 1} + B_L,
    \end{equation}
    let $A$ be a finite set,
    let $\action \in A$,
    $\delta \in \pr{0,1}$,
    $\eta \in (0,1/\delta)$,
    $\genConst \in [ 1 , \infty)$,
    let $(\Omega, \mathcal{F}, \mathbb{P})$ be a probability space,
    for every $\dm \in \N$ let $g_\dm \colon \R^\dm \to \R^{\abs{A}}$ and $\xi_\dm \colon \Omega \to \R^\dm$ be measurable,
    for every $\dm \in \N$, $a \in A$ let $t^{\pr{a}}_\dm \colon \R^{2\dm} \to \R^\dm$ be a measurable function which\footnote{Note that for all $\dm \in \N$, $x = \pr{x_1, \dots, x_\dm} \in \R^\dm$ it holds that $\norm{x} = \br[\big]{\sum_{i=1}^{\dm} \pr{x_i}^{2}}^{\nicefrac{1}{2}}$ (standard norm).} satisfies $\sup_{x \in \R^\dm} \E \br[\big]{ \norm{t_\dm^{\pr{a}} \pr{x, \xi_\dm}}  } < \infty$,
    for every $\dm \in \N$, $\eps \in \err$, $a \in A$ let $\ft_{\dm, \eps}^{\pr{a}} \in \bN_{2\dm, \dm}$ satisfy $\cD\pr{ \ft_{\dm, \eps}^{\pr{a}}} = \cD \pr{ \ft_{\dm, \eps}^{\pr{\action}}}$,
    and for every $\dm \in \N$, $\eps \in \err$ let $\bG_{\dm, \eps} \in \bN_{\dm, \abs{A}}$ satisfy for all $x, y, p \in \R^\dm$, $a \in A$ that
    \begin{equation}
        \begin{split} \label{eq:intro_approximation_assumption}
            \varepsilon \norm{ g_\dm\pr{x}} + \norm[]{ \pr[]{ \cR (\bG_{\dm, \eps}) }(x) -  g_\dm (x) } + \norm[]{ \pr[]{\cR \pr{\ft_{\dm, \eps}^{\pr{a}}}} \pr{x, p} - t_\dm^{\pr{a}}\pr{x, p} } \le \eps \genConst \dm^{\genConst},
        \end{split} 
    \end{equation}
    \begin{equation}
        \begin{split} \label{eq:intro_lipschitz}
            \norm{ g_\dm \pr{x} - g_\dm\pr{y}} \le \genConst \norm{ x - y}, \qquad
            \norm[]{ t_\dm^{\pr{a}}\pr{ x, p } - t_\dm^{\pr{a}}\pr{ y, p }  }  \le \eta \norm{x - y}, 
            %
        \end{split}
    \end{equation}
    and $\paramANN \pr{\bG_{\dm, \eps}} + \paramANN \pr{ \ft_{\dm, \eps }^{\pr{\action}} }  \le \genConst \dm^\genConst \eps^{- \genConst}$.
    Then there exists $c \in \R$ such that for every $\dm \in \N$, $\eps \in \err$ and every probability measure $\mu \colon \Borel(\R^\dm) \to \br{0,1}$ it holds that
    \begin{enumerate}[label=(\roman*)]
        \item \label{item:intro_ex_uniq_sol} there exists a unique bounded measurable $\qfunc \colon \R^\dm \to \R^{\abs{A}}$ which satisfies for all $x \in \R^\dm$, $a \in A$ that
        \begin{equation} \label{eq:intro_bellman_eq}
            \pr{\qfunc \pr{x}} \pr{a} = \pr{g_\dm\pr{x}} \pr{a} + \delta \, \E  \br*{ \max_{b \in A} \pr{ \qfunc \pr{ t_{\dm}^{\pr{a}}\pr{x, \xi_\dm} } } \pr{b} } 
        \end{equation}
        and
        \item \label{item:intro_ex_ann} there exists $\bQ \in \bN_{\dm, \abs{A}}$ which satisfies
        \begin{equation} \label{eq:intro_approximating_ann}
            \textstyle \pr*{ \int_{\R^\dm} \norm{\qfunc \pr{x} - \pr*{ \cR(\bQ) }(x)}^2 \, \mu(\dxx x) }^{\nicefrac{1}{2}} \le \eps \qquad \text{and} \qquad \paramANN \pr{\bQ} \le c \dm^c \eps^{-c}.
        \end{equation}
    \end{enumerate}
\end{theo}
\end{samepage}


\cref{thm:intro} is an immediate consequence of \cref{cor:bellman} in \cref{subsec:bellman}. \cref{cor:bellman} follows from \cref{cor:activation_transition} in \cref{subsec:activation_transition}. \cref{cor:activation_transition}, in turn, is a consequence of \cref{thm:main_simple} in \cref{subsec:general_fp_equation}, the main result of this article.
In the following we add some explanatory sentences on the mathematical objects appearing in \cref{thm:intro} above.

The function $\pr{ \cup_{\dm \in \N} \R^\dm } \ni x \mapsto \fm\pr{x} \in \pr{ \cup_{\dm, \in \N} \R^\dm}$ represents a suitable multidimensional version of the one-dimensional activation function $\R \ni x \mapsto \max \cu{x, \beta x} \in \R$.
The real number $\beta \in [0, \infty) \backslash \cu{1}$ is a parameter that determines the choice of activation function. In particular, in the case $\beta = 0$ we have that $\fm$ is the multidimensional version of the ReLU activation and in the case $\beta \in \pr{0,1}$ we have that $\fm$ is the multidimensional version of the leaky ReLU activation with leaky factor $\beta$. 

The sets $\bN_{\dm, \fdm}$, $\dm, \fdm \in \N$, in \eqref{eq:anns_intro} in \cref{thm:intro} denote the sets of all artificial neural networks (ANNs) we consider in the approximation of the solutions of the Bellman equations for the $Q$-function discussed. Note that for all $\dm, \fdm \in \bN$ and for every ANN $\Phi \in \bN_{\dm, \fdm}$ the tuple of natural numbers $\cD \pr{\Phi} \in \cup_{H \in \N} \pr{\cu{\dm} \times \N^H \times \cu{\fdm}}$ denotes the architecture of the ANN $\Phi$. Moreover, we note that for all $\dm, \fdm \in \N$ and every ANN $\Phi \in \bN_{\dm, \fdm}$ we have that the function
\begin{equation}
    \cR \pr{\Phi} \colon \R^\dm \to \R^\fdm
\end{equation}
corresponds to the realization of the ANN $\Phi$ with the activation function $\R \ni x \mapsto \max \cu{x, \beta x} \in \R$.
Furthermore, observe that for all $\dm, \fdm \in \N$ and every ANN $\Phi \in \bN_{\dm, \fdm}$ the natural number $\paramANN\pr{\Phi} \in \N$ is the number of real numbers used to describe the ANN $\Phi$. In particular, we note that for all $\dm, \fdm \in \N$ and every $\Phi \in \bN_{\dm, \fdm}$ we have that $\paramANN \pr{\Phi}$ is corresponding in a rough sense to the amount of computer memory that is needed in order to store the ANN $\Phi$.

The triple $\pr{\Omega, \mathcal{F}, \mathbb{P}}$ in \cref{thm:intro} is the probability space on which all appearing random variables are defined. In \cref{thm:intro} we consider for every $\dm \in \N$ an MDP with state space $\pr{\R^\dm, \Borel \pr{\R^\dm}}$. The finite set $A$ denotes the common control set for all elements of the sequence of MDPs. For every $\dm \in \N$, $x \in \R^\dm$, $a \in A$ the one-step transition of the controlled Markov chain of the MDP with state space $\pr{\R^\dm, \Borel \pr{\R^\dm}}$ by choosing control $a$ in state $x$ evolves as the composition of the function $t^{\pr{a}}_\dm \pr{x, \cdot} \colon \R^{\dm} \to \R^\dm$ and the random variable $\xi_\dm \colon \Omega \to \R^\dm$, namely by
\begin{equation} \label{eq:intro_one_step}
    t_\dm^{\pr{a}} \pr{x, \xi_\dm} \colon \Omega \to \R^\dm.
\end{equation}
For every $\dm \in \N$ the function $g_\dm \colon \R^\dm \to \R^{\abs{A}}$ in \cref{thm:intro} denotes the reward function of the MDP with state space $\pr{\R^\dm, \Borel \pr{\R^\dm}}$. We assume that for all $\dm \in \N$ the function $g_\dm$ is uniformly bounded and that the sequence of these bounds grows at most polynomially in the state space dimension $\dm \in \N$. The real number $\delta \in \pr{0,1}$ is the common discount factor for all elements of the sequence of MDPs.

The real numbers $\eta \in \pr{0, \nicefrac{1}{\delta}}$ and $\genConst \in [1, \infty)$ are constants that we utilize to formulate our regularity and approximation assumptions in \cref{thm:intro}.
We assume there exist ANNs $\bG_{\dm, \eps} \in \bN_{\dm, \abs{A}}$, $\dm \in \N$, $\eps \in \err$, such that for every $\dm \in \N$, $\eps \in \err$ the realization function $\cR \pr{\bG_{\dm, \eps}} \colon \R^\dm \to \R^{\abs{A}}$ of the ANN $\bG_{\dm, \eps}$ is a suitable approximation of the function $g_\dm$ and that the number of parameters $\paramANN \pr{\bG_{\dm, \eps}}$ of the approximating ANN $\bG_{\dm, \eps}$ is bounded by $\genConst \dm^\genConst \eps^{-\genConst}$.
Furthermore, we assume there exist $\ft_{\dm, \eps}^{\pr{a}} \in \bN_{2\dm, \dm}$, $\dm \in \N$, $\eps \in \err$, $a \in A$, such that for every $\dm \in \N$, $\eps \in \err$, $a \in A$ the realization function $\cR \pr{ \ft_{\dm, \eps}^{\pr{a}} } \colon \R^{2\dm} \to \R^\dm$ of the ANN $\ft_{\dm, \eps}^{\pr{a}}$ is a suitable approximation of the function $t_\dm^{\pr{a}}$ and that the number of parameters $\paramANN \pr{\ft_{\dm, \eps}^{\pr{a}}}$ of the approximating ANN $\ft_{\dm \eps}^{\pr{a}}$ is bounded by $\genConst \dm^\genConst \eps^{-\genConst}$.

\cref{thm:intro} establishes in \cref{item:intro_ex_uniq_sol} the essentially well-known result that under the above assumptions it holds that for every $\dm \in \N$ the Bellman equation \eqref{eq:intro_bellman_eq} for the $Q$-function associated to the MDP with state space $(\R^\dm, \Borel\pr{\R^\dm})$ has a unique solution $\qfunc \colon \R^\dm \to \R^{\abs{A}}$.

Moreover, \cref{thm:intro} concludes in \cref{item:intro_ex_ann} that there exists a constant $c \in \R$ such that for every state space dimension $\dm \in \N$, every prescribed approximation accuracy $\eps \in \err$, and every probability measure $\mu \colon \Borel \pr{\R^\dm} \to \br{0,1}$ we have that there exists an ANN $\bQ \in \bN_{\dm, \abs{A}}$ such that the realization
\begin{equation}
    \cR \pr{\bQ} \colon \R^\dm \to \R^{\abs{A}}
\end{equation}
of the ANN $\bQ$ approximates the solution $\qfunc \colon \R^\dm \to \R^{\abs{A}}$ of the Bellman equation for the $Q$-function in the $L^2$-sense with respect to the probability measure $\mu$ on $\R^\dm$ with prescribed approximation accuracy $\eps$ and such that the number of parameters $\paramANN \pr{\bQ}$ of the approximating ANN $\bQ$ is bounded by $c \dm^c \eps^{-c}$.

The remainder of this article is structured as followed. In \cref{sec:stability} below we study stability properties for solutions of certain stochastic fixed point equations. In \cref{sec:ann_calculus} below we recall the necessary calculus for ANNs. In \cref{sec:ann_representations} we establish ANN representations for MLFP approximations for certain stochastic fixed-point equations. In \cref{sec:ANN_approximations_fixed_point_equations} we combine the stability result in \cref{sec:stability} and the ANN representations of the MFLP approximation scheme in \cref{sec:ann_representations} to prove \cref{thm:intro} and its generalizations in \cref{subsec:general_fp_equation}, \cref{subsec:activation_transition}, and \cref{subsec:bellman}.

\section{A stability result for solutions of certain stochastic fixed point equations}
\label{sec:stability}

\subsection{Setting}
\label{subsec:stability_setting}

\begin{definition}[Couplings] \label{def:couplings}
    Let $\dm \in \N$,
    let $\mu, \nu \colon \Borel(\R^\dm) \to \br{0,1}$ be probability measures.
    Then we denote by $\coup_\dm\pr{\mu, \nu}$ the set of all couplings of $\mu$ and $\nu$ given by
    \begin{equation}
        \coup_\dm\pr{\mu, \nu} = \cu*{\gamma : \substack{\gamma \text{ is a probability measure on } (\R^\dm \times \R^\dm, \Borel(\R^\dm) \otimes \Borel (\R^\dm)), \\ \text{for all } A \in \Borel(\R^\dm) \text{ it holds that } \mu(A) = \gamma (A \times \R^\dm) \text{ and } \nu(A) = \gamma(\R^\dm \times A) }}.
    \end{equation}
\end{definition}

\cfclear
\begin{definition}[Wasserstein-1 space] \label{def:wasserstein_space}
    \cfconsiderloaded{def:wasserstein_space}
    Let $\dm \in \N$.
    Then we denote by $\meas_\dm$ the set given by
    \begin{equation}
        \meas_\dm = \cu*{ \pr{\mu \colon \Borel(\R^\dm) \to \br{0,1}} : \mu \text{ is a probability measure with } \varint_{\R^\dm}\norm{x}\mu(\dxx x) < \infty}.
    \end{equation}
\end{definition}

\cfclear
\begin{definition}[Wasserstein-1 distance] \label{def:wasserstein}
    \cfconsiderloaded{def:wasserstein}
    Let $\dm \in \N$. Then we denote by $\wdist_\dm \colon \meas_\dm \times \meas_\dm \to [0, \infty)$ the function which satisfies for all $\mu, \nu \in \meas_\dm$ that
    \begin{equation}
        \wdist_\dm \pr{\mu, \nu} = \inf_{\pi \in \coup_\dm\pr{\mu, \nu}} \int_{\R^\dm \times \R^\dm} \norm{x - y} \pi \pr{ \dxx \pr{x, y}}
    \end{equation}
    \cfout.
\end{definition}

\begin{setting} \label{set:stability_w}
    Let $\dm \in \N$,
    let $\w \colon \R^\dm \to (0, \infty)$ be measurable,
    let $A$ be a nonempty finite set,
    and let
    \begin{equation}
        \fW = \cu[\bigg]{(u \colon \R^\dm \rightarrow \R^{\abs{A}}) : u \text{ is measurable, }  \sup_{\substack{x \in \R^\dm, a \in A}} \tfrac{| (u(x))(a) |}{|\w(x)|} < \infty}.
    \end{equation}
\end{setting}

\subsection{A stability result for solutions of certain stochastic fixed point equations}
\label{subsec:stability_results}

\cfclear
\begin{definition}[Stochastic kernel] \label{def:stochastic_kernel} 
    Let $\pr{\X, \mathcal{X}}$ and $\pr{\Y, \mathcal{Y}}$ be nonempty measurable spaces
    and let $\tker \colon \X \times \mathcal{Y} \to \br{0,1}$ satisfy for all $M \in \mathcal{Y}$ that $\X \ni x \mapsto \tker\pr{x, M} \in \br{0,1}$ is measurable and for all $x \in \X$ that $\mathcal{Y} \ni M \mapsto \tker\pr{x, M} \in \br{0,1}$ is a probability measure on $\pr{\Y, \mathcal{Y}}$. Then we say that $\tker$ is a stochastic kernel from $\pr{\X, \mathcal{X}}$ to $\pr{\Y, \mathcal{Y}}$.
\end{definition}

\cfclear
\begin{lemma} \label{lem:stability_kernel_w_nonlinearity}
    \cfconsiderloaded{lem:stability_kernel_w_nonlinearity}
    Assume \cref{set:stability_w},
    let $c, L_i \in [0, \infty)$, $i \in \{1,2\}$, satisfy for all $i \in \{1,2\}$ that $cL_i < 1$,
    let $f_i \colon \R^\dm \times \R^{\abs{A}} \to \R$, $i \in \{1,2\}$, be measurable,
    assume for all $i\in \{1,2\}$, $x \in \R^\dm$, $r,s\in \R^{\abs{A}}$ that $\abs{ f_i(x,r) - f_i(x,s) } \le L_i \max_{a \in A} \abs{ r(a) - s(a)}$,
    let $\tker^{(a)} \colon \R^\dm \times \Borel(\R^\dm) \to \br{0,1}$, $a \in A$, be stochastic kernels,
    assume for all $x \in \R^\dm$, $a \in A$ that $\int_{\R^\dm} \abs{\w(y)} \tker^{\pr{a}}(x, \dxx y) \le c \w(x)$,
    assume for all $i \in \cu{1,2}$ that
    \begin{equation}
        \begin{split}
            \sup_{(x, a)\in \R^\dm \times A} \abs{\w(x)}^{-1}\int_{\R^\dm} \abs{ f_i(y, 0)} \tker^{(a)}(x, \dxx y) < \infty.
        \end{split}
    \end{equation}
    Then
    \begin{enumerate}[label=(\roman *)]
        \item \label{item:stability_kernel_w_nonlinearity_ex_uniq} there exist unique $u_i \in \W$, $i \in \cu{1,2}$, such that for every $i \in \cu{1,2}$, $x \in \R^\dm$, $a \in A$ it holds that
        \begin{equation} \label{eq:stability_kernel_w_nonlinearity_fp_eq}
            \int_{\R^\dm} \abs*{f_i\pr{y, u_i(y)}} \tker^{\pr{a}}\pr{x, \dxx y} y < \infty \qquad \text{and} \qquad (u_i(x))(a) = \int_{\R^\dm} f_i(y, u_i(y)) \tker^{(a)}(x, \dxx y)
        \end{equation}
        and
        \item \label{item:stability_kernel_w_nonlinearity_result} it holds that
        \begin{equation} \label{eq:stability_kernel_w_nonlinearity}
            \br[\bigg]{\sup_{(x,a) \in \R^\dm \times A}\!\!\!\! \tfrac{\abs[]{ (u_1(x))(a) - (u_2(x))(a) }}{ \abs{\w(x)} }} \le \tfrac{c}{1-c (\min\{L_1, L_2\})} \br[\bigg]{\sup_{(y,r) \in \R^\dm \times \R^{\abs{A}}}\!\!\!\! \tfrac{ \abs{ f_1(y,r) - f_2(y,r) } }{ \abs{\w(y)} }}
        \end{equation}
    \end{enumerate}    
    \cfout.
\end{lemma}

\begin{mproof}{\cref{lem:stability_kernel_w_nonlinearity}}
    First, note that the assumption that $\w$ is measurable,
    the assumption that for all $i \in \cu{1,2}$ it holds that $cL_i < 1$,
    the assumption that for all $i \in \cu{1,2}$, $x \in \R^\dm$, $r,s \in \R^{\abs{A}}$ it holds that $\abs{f_i \pr{x,r} - f_i \pr{x,s}} \le L_i \max_{a \in A} \abs{r\pr{a} - s\pr{a}}$,
    the assumption that for all $x \in \R^\dm$, $a \in A$ it holds that $\int_{\R^\dm} \abs{\w(y)} \tker^{\pr{a}}(x, \dxx y) \le c \w(x)$,
    the assumption that for all $i \in \cu{1,2}$ it holds that $\sup_{(x, a)\in \R^\dm \times A} \abs{\w(x)}^{-1}\int_{\R^\dm} \abs{ f_i(y, 0)} \tker^{(a)}(x, \dxx y) < \infty$,
    and \cite[Lemma~2.2]{beck2023nonlinctrl} (applied for every $i\in \cu{1,2}$ with $c\with c$, $L \with L_i$, $\pr{\X, \mathcal{X}} \with \pr{\R^\dm, \cB\pr{\R^\dm}}$, $A \with A$, $\pr{\tker_a}_{a \in A} \with \pr{ \tker^{\pr{a}} }_{a \in A}$, $f \with f_i$, $\w \with \br{ \R^\dm \ni x \mapsto \br{A \ni a \mapsto \w(x) \in \pr{0, \infty}} \in \pr{0,\infty}^{\abs{A}} }$ in the notation of \cite[Lemma~2.2]{beck2023nonlinctrl})
    establish \cref{item:stability_kernel_w_nonlinearity_ex_uniq}.
    Combining this,
    the assumption that for all $x \in \R^\dm$, $r,s \in \R^{\abs{A}}$ it holds that $\abs{f_1(x,r) - f_1(x,s)} \le L_1 \max_{a \in A} \abs{r(a) - s(a)}$,
    the assumption that for all $x \in \R^\dm$, $a \in A$ it holds that $\int_{\R^\dm} \abs{\w(y)} \tker^{\pr{a}}(x, \dxx y) \le c \w(x)$,
    and the triangle inequality yields for all $x \in \R^\dm$, $a \in A$ that
    \begin{align}
        \begin{split}\label{eq:stability_kernel_w_nonlinearity_for_x_a}
            \tfrac{\abs{ (u_1(x))(a) - (u_2(x))(a) }}{ \abs{\w(x)} } &= \tfrac{1}{\abs{\w(x)}} \abs*{ \int_{\R^\dm} f_1\pr{y, u_1(y)} \tker^{\pr{a}}(x, \dxx y) - \int_{\R^\dm} f_2(y, u_2(y)) \tker^{\pr{a}}(x, \dxx y) } \\
            &= \tfrac{1}{\abs{\w(x)}} \abs*{ \int_{\R^\dm} f_1(y, u_1(y)) - f_2(y, u_2(y)) \tker^{(a)}(x, \dxx y)}\\
            &\le \tfrac{1}{\abs{\w(x)}} \pr[\bigg]{ \int_{\R^\dm} \abs{f_1 \pr{y, u_1(y)} - f_1\pr{y, u_2(y)} } \tker^{\pr{a}}(x, \dxx y) \\
            &\quad+ \int_{\R^\dm} \abs{ f_1\pr{ y, u_2(y) } - f_2\pr{y, u_2(y)} } \tker^{\pr{a}}(x, \dxx y) } \\
            &\le \tfrac{1}{\abs{\w(x)}} \pr[\bigg]{ L_1 \int_{\R^\dm} \max_{b \in A} \abs{(u_1(y))(b) - (u_2(y))(b)} \tker^{\pr{a}}(x, \dxx y) \\
            &\quad+ \int_{\R^\dm} \abs{ f_1\pr{y, u_2(y) } - f_2\pr{y, u_2(y)} } \tker^{\pr{a}}(x, \dxx y) } \\
            &\le \tfrac{L_1}{\abs{\w(x)}} \br[\bigg]{ \sup_{(z, c) \in \R^\dm \times A} \!\!\!\! \tfrac{\abs{ (u_1(z))(c) - (u_2(z))(c) }}{ \abs{\w(z)}} } \int_{\R^\dm} \abs{\w(y)} \tker^{\pr{a}}(x, \dxx y) \\
            &\quad+ \tfrac{1}{\abs{\w(x)}} \br[\bigg]{ \sup_{(z, r) \in \R^\dm \times \R^{\abs{A}}} \!\!\!\! \tfrac{ \abs{ f_1(z,r) - f_2(z, r) } }{ \abs{\w(z)} } } \int_{\R^\dm} \abs{\w(y)} \tker^{\pr{a}}(x, \dxx y) \\
            &\le cL_1 \br[\bigg]{ \sup_{(z, c) \in \R^\dm \times A} \!\!\!\! \tfrac{\abs{ (u_1(z))(c) - (u_2(z))(c) }}{ \abs{\w(z)}} } + c \br[\bigg]{ \sup_{(z, r) \in \R^\dm \times \R^{\abs{A}}} \!\!\!\! \tfrac{ \abs{ f_1(z,r) - f_2(z, r) } }{ \abs{\w(z)} } }.
        \end{split}
    \end{align}
    Moreover, note that \cref{item:stability_kernel_w_nonlinearity_ex_uniq} proves that 
    \begin{equation}
        \br[\bigg]{ \sup_{(y, b) \in \R^\dm \times A} \!\!\! \tfrac{ \abs{ (u_1(y))(b) - (u_2(y))(b) } }{ \abs{\w(y)} } } < \infty.
    \end{equation}
    This and \eqref{eq:stability_kernel_w_nonlinearity_for_x_a} imply that
    \begin{equation}
        (1 - cL_1) \br[\bigg]{ \sup_{(x, a) \in \R^\dm \times A} \!\!\!\! \tfrac{\abs{ (u_1(x))(a) - (u_2(x))(a) }}{ \abs{\w(x)}} } \le c \br[\bigg]{ \sup_{(y, r) \in \R^\dm \times \R^{\abs{A}}} \!\!\!\! \tfrac{ \abs{ f_1(y,r) - f_2(y, r) } }{ \abs{\w(y)}}}.
    \end{equation}
    This and the assumption that $cL_1 < 1$ demonstrate that 
    \begin{equation} \label{eq:stability_kernel_w_nonlinearity_1}
        \br[\bigg]{ \sup_{(x, a) \in  \R^\dm \times A} \!\!\!\! \tfrac{\abs{ (u_1(x))(a) - (u_2(x))(a) }}{ \abs{\w(x)}} } \le \tfrac{c}{1 - cL_1} \br[\bigg]{ \sup_{(y, r) \in \R^\dm \times \R^{\abs{A}}} \!\!\!\! \tfrac{ \abs{ f_1(y,r) - f_2(y, r) } }{ \abs{\w(y)}}}.
    \end{equation}
    Symmetry yields
    \begin{equation}
        \br[\bigg]{ \sup_{(x, a) \in \R^\dm \times A} \!\!\!\! \tfrac{\abs{ (u_1(x))(a) - (u_2(x))(a) }}{ \abs{\w(x)}} } \le \tfrac{c}{1 - cL_2} \br[\bigg]{ \sup_{(y, r) \in \R^\dm \times \R^{\abs{A}}} \!\!\!\! \tfrac{ \abs{ f_1(y,r) - f_2(y, r) } }{ \abs{\w(y)}}}.
    \end{equation}
    This and \eqref{eq:stability_kernel_w_nonlinearity_1} establish that
    \begin{equation}
        \br[\bigg]{ \sup_{(x, a) \in \R^\dm \times A} \!\!\!\! \tfrac{\abs{ (u_1(x))(a) - (u_2(x))(a) }}{ \abs{\w(x)}} } \le \tfrac{c}{1 - c(\min\{ L_1, L_2 \})} \br[\bigg]{ \sup_{(y, r) \in \R^\dm \times \R^{\abs{A}}} \!\!\!\! \tfrac{ \abs{ f_1(y,r) - f_2(y, r) } }{ \abs{\w(y)}}}.
    \end{equation}
    This proves \cref{item:stability_kernel_w_nonlinearity_result}. The proof of \cref{lem:stability_kernel_w_nonlinearity} is thus complete.
\end{mproof}

\cfclear
\begin{lemma} \label{lem:lipschitz_continuity_of_solutions_kernel_w}
    \cfconsiderloaded{lem:lipschitz_continuity_of_solutions_kernel_w}
    Assume \cref{set:stability_w},
    let $\eta, c, L, K \in [0, \infty)$,
    assume $cL < 1$ and $\eta L < 1$,
    let $f \colon \R^\dm \times \R^{\abs{A}} \to \R$ be measurable,
    assume for all $x, y \in \R^\dm$, $r, s \in \R^{\abs{A}}$ that $\abs{ f(x,r) - f(x, s) } \le L \max_{a \in A} \abs{r(a) - s(a)}$ and $\abs{f(x,r) - f(y, r)} \le K \norm{x - y}$,
    let $\tker^{\pr{a}} \colon \R^\dm \times \Borel(\R^\dm) \to \br{0,1}$, $a \in A$, be stochastic kernels,
    assume for all $x \in \R^\dm$, $a \in A$ that $\int_{\R^\dm} \abs{\w(y)} \tker^{\pr{a}}(x, \dxx y) \le c \w(x)$ and $\int_{\R^\dm} \norm{y} \tker^{\pr{a}}(x, \dxx y) < \infty$,
    assume
    \begin{equation}
        \begin{split}
            \sup_{(x, a) \in \R^\dm \times A} \abs{\w(x)}^{-1} \int_{\R^\dm} \abs{f\pr{ y, 0 }} \tker^{\pr{a}}(x, \dxx y) < \infty,
        \end{split}
    \end{equation}
    assume for all $x, y \in \R^\dm$, $a \in A$ that $\wdist_\dm \pr{ \tker^{\pr{a}}(x, \cdot), \tker^{\pr{a}}(y, \cdot) } \le \eta \norm{x - y}$.
    Then
    \begin{enumerate}[label=(\roman *)]
        \item \label{item:lipschitz_continuity_of_solutions_kernel_w_ex_uniq} there exists a unique $u \in \W$ such that for every $x \in \R^\dm$, $a \in A$ it holds that
        \begin{equation}
            \int_{\R^\dm} \abs{ f\pr{y, u(y)} } \tker^{\pr{a}}\pr{x, \dxx y} < \infty \qquad \text{and} \qquad \pr{u\pr{x}}\pr{a} = \int_{\R^\dm} f\pr{y, u(y)} \tker^{\pr{a}}\pr{x, \dxx y}
        \end{equation}
        and
        \item \label{item:lipschitz_continuity_of_solutions_kernel_w_result} it holds for all $x, y \in \R^\dm$, $a \in A$ that 
        \begin{equation} \label{eq:lipschitz_continuity_of_solutions_kernel_w}
            \abs{ (u(x))(a) - (u(y))(a)} \le \tfrac{\eta K}{1 - \eta L} \norm{x - y}
        \end{equation}
    \end{enumerate}
    \cfout.    
\end{lemma}

\begin{mproof}{\cref{lem:lipschitz_continuity_of_solutions_kernel_w}}
    First, observe that the assumption that $\w$ is measurable,
    the assumption that $cL < 1$,
    the assumption that for all $x \in \R^\dm$, $r,s \in \R^{\abs{A}}$ it holds that $\abs{ f\pr{x,r} - f\pr{x,s} } \le L \max_{a \in A} \abs{r\pr{a} - s\pr{a}}$,
    the assumption that for all $x \in \R^\dm$, $a \in A$ it holds that $\int_{\R^\dm} \abs{\w(y)} \tker^{\pr{a}}\pr{x, \dxx y} \le c \w(x)$,
    the assumption that $\sup_{(x, a) \in \R^\dm \times A} \abs{\w(x)}^{-1} \int_{\R^\dm} \abs{f\pr{ y, 0 }} \tker^{\pr{a}}(x, \dxx y) < \infty$,
    and \cite[Lemma~2.2]{beck2023nonlinctrl} (applied with $c \with c$, $L \with L$, $\pr{\X, \mathcal{X}} \with \pr{\R^\dm, \cB(\R^\dm)}$, $A \with A$, $\pr{\tker_a}_{a \in A} \with \pr{\tker^{\pr{a}}}_{a \in A}$, $f \with f$, $\w \with \br{ \R^\dm \ni x \mapsto \br{ A \ni a \mapsto \w(x) \in \pr{0,\infty} } \in \pr{0, \infty}^{\abs{A}} }$ in the notation of \cite[Lemma~2.2]{beck2023nonlinctrl}) demonstrates \cref{item:lipschitz_continuity_of_solutions_kernel_w_ex_uniq}.
    Throughout the remainder of this proof let $u_n \colon \R^\dm \rightarrow \R^{\abs{A}}$, $ n \in \N_0$, satisfy for all $n \in \N$, $x \in \R^\dm$, $a \in A$ that $u_0 = 0$ and 
    \begin{equation} \label{eq:lipschitz_continuity_of_solutions_kernel_w_def_iterate}
        \begin{split}
            (u_n(x))(a) = \int_{\R^\dm} f\pr{y, u_{n - 1}(y)} \tker^{\pr{a}}(x, \dxx y).
        \end{split}
    \end{equation}
    Next, note that the assumption that for all $a \in A$ it holds that $\tker^{\pr{a}}$ is a stochastic kernel,
    the assumption that for all $x \in \R^\dm$, $r, s \in \R^{\abs{A}}$ it holds that $|f(x, r) - f(x, s)| \le L \max_{a \in A} |r(a)-s(a)|$,
    the assumption that $cL \in [0,1)$,
    the assumption that for all $x \in \R^\dm$, $a \in A$ it holds that $\int_{\R^\dm} \abs{\w(y)} \tker^{\pr{a}}(x, \dxx y) \le c \w(x)$,
    the assumption that $\sup_{(x, a) \in \R^\dm \times A} \abs{\w(x)}^{-1} \int_{\R^\dm} \abs{f(y, 0)} \tker^{\pr{a}}(x, \dxx y) < \infty$,
    and induction ensure for all $n \in \N_0$ that $u_n \in \fW$.
    Combining the assumption that for all $x \in \R^\dm$, $r, s \in \R^A$ it holds that $ |f(x,r) - f(x, s)| \le L \max_{a \in A} | r(a) - s(a) |$
    and the assumption that $cL \in [0,1)$
    with \cref{item:lipschitz_continuity_of_solutions_kernel_w_ex_uniq}
    yields that $ \lim_{n \to \infty} \big( \sup_{(x,a) \in \R^\dm \times A}  \abs{\w(x)}^{-1} \abs{(u_n(x))(a) - (u(x))(a)} \big) = 0$.
    We prove by induction that for all $n \in \N_0$, $x, y \in \R^\dm$, $a \in A$ it holds that 
    \begin{equation} \label{eq:lipschitz_continuity_of_fixed_point_iterates_kernel_w}
        \abs{ (u_n(x))(a) - (u_n(y))(a) } \le  \br*{\eta K \SmallSum{i = 0}{n - 1} ( \eta L )^i} \norm{x - y}.
    \end{equation}
    Note that for all $x, y \in \R^\dm$, $a \in A$ it holds that $ | (u_0(x))(a) - (u_0(y))(a) | = 0 = [\eta K \sum_{i = 0}^{-1} (\eta L)^i] \|x - y\| $.
    This proves \eqref{eq:lipschitz_continuity_of_fixed_point_iterates_kernel_w} for the base case $n = 0$.
    For the induction step $\N_0 \ni (n - 1) \induct n \in \N$ fix $n \in \N$ and assume for all $l \in \N_0 \cap [0, n - 1] $, $x, y \in \R^\dm$, $a \in A$ that
    \begin{equation}
        \abs{ (u_l(x))(a) - (u_l(y))(a) } \le \br*{\eta K \SmallSum{i = 0}{l - 1} ( \eta L )^i} \norm{x - y}.
    \end{equation}
    Moreover, note that for all $x, y \in \R^\dm$, $a \in A$, $\gamma \in \coup_\dm( \tker^{\pr{a}}(x, \cdot ), \tker^{\pr{a}}(y, \cdot))$, $B \in \Borel(\R^\dm)$ it holds that $\gamma(B \times \R^\dm) = \tker^{\pr{a}}(x, B)$ and $\gamma(\R^\dm \times B) = \tker^{\pr{a}}(y, B)$.
    Combining this, \eqref{eq:lipschitz_continuity_of_solutions_kernel_w_def_iterate},
    the triangle inequality,
    the assumption that for all $x, y \in \R^\dm$, $r, s \in \R^A$ it holds that $\abs{ f(x,r) - f(x,s) } \le L \max_{a \in A} \abs{r(a) - s(a)} $ and $\abs{f(x,r) - f(y, r)} \le K \norm{x - y}$,
    and the induction hypothesis that for all $l\in \N_0 \cap \br{ 0, n - 1 } $, $x, y \in \R^\dm$, $a \in A$ it holds that $\abs{ (u_l(x))(a) - (u_l(y))(a) } \le \br{ \eta K \sum_{i = 0}^{l - 1} (\eta L)^i } \norm{x - y}$ imply for all $x, y \in \R^\dm$, $a \in A$, $\gamma \in \coup_\dm(\tker^{\pr{a}}(x, \cdot), \tker^{\pr{a}}(y, \cdot))$ that
    \begin{align}
            \abs{ (u_n(x))(a) - (u_n(y))(a)} &= \abs*{ \int_{\R^\dm} f(\bfx, u_{n - 1}(\bfx)) \tker^{\pr{a}}(x, \dxx \bfx) - \int_{\R^\dm} f(\bfy, u_{n - 1}(\bfy)) \tker^{\pr{a}}(y, \dxx\bfy)} \nonumber \\
            &= \abs*{ \int_{\R^\dm \times \R^\dm} \hspace{-0.7cm} f(\bfx, u_{n - 1}(\bfx)) \gamma(\dxx(\bfx, \bfy)) - \int_{\R^\dm \times \R^\dm} \hspace{-0.7cm} f(\bfy, u_{n - 1}(\bfy)) \gamma(\dxx(\bfx, \bfy)) } \nonumber \\
            &\le \int_{\R^\dm \times \R^\dm}\hspace{-0.5cm} \abs*{ f(\bfx, u_{n - 1}(\bfx)) - f(\bfy, u_{n - 1}(\bfy)) } \gamma (\dxx(\bfx, \bfy)) \nonumber \\
            &\le \int_{\R^\dm \times \R^\dm} \hspace{-0.5cm} \abs*{ f(\bfx, u_{n - 1}(\bfx)) - f(\bfx, u_{n - 1}(\bfy)) } \gamma (\dxx(\bfx, \bfy)) \nonumber \\
            &\quad + \int_{\R^\dm \times \R^\dm}\hspace{-0.5cm} \abs*{ f(\bfx, u_{n - 1}(\bfy)) - f(\bfy, u_{n - 1}(\bfy)) } \gamma (\dxx(\bfx, \bfy)) \nonumber \\
            &\le L \int_{\R^\dm \times \R^\dm}  \max_{b \in A} \abs{ (u_{n - 1}(\bfx))(b) - (u_{n - 1}(\bfy))(b) } \gamma(\dxx(\bfx, \bfy)) \nonumber \\
            &\quad+ K \int_{\R^\dm \times \R^\dm} \hspace{-0.7cm} \norm{ \bfx - \bfy } \gamma(\dxx(\bfx, \bfy)) \nonumber \\
            &\le L \br*{ \eta K \SmallSum{i = 0}{n - 2} (\eta L)^i } \int_{\R^\dm \times \R^\dm} \hspace{-0.7cm} \norm{\bfx - \bfy} \gamma(\dxx(\bfx, \bfy)) + K \!\! \int_{\R^\dm \times \R^\dm} \hspace{-0.7cm} \norm{\bfx - \bfy} \gamma(\dxx(\bfx, \bfy)) \nonumber \\
            &= \br*{K \SmallSum{i = 0}{n - 1}(\eta L)^i} \int_{\R^\dm \times \R^\dm} \hspace{-0.7cm} \norm{\bfx - \bfy} \gamma(\dxx(\bfx, \bfy)).
    \end{align}
    This and the assumption that for all $x, y \in \R^\dm$, $a \in A$ it holds that $\wdist_\dm ( \tker^{\pr{a}}(x, \cdot), \tker^{\pr{a}}(y, \cdot) ) \le \eta  \norm{x - y} $ assure for all $x, y \in \R^\dm$, $a\in A$ that
    \begin{align}
        \abs{ (u_n(x))(a) - (u_n(y))(a) } &\le  \br*{ K \SmallSum{i = 0}{n - 1} (\eta L)^i } \inf_{\gamma \in \coup_\dm(\tker^{\pr{a}}(x,\cdot), \tker^{\pr{a}}(y,\cdot))} \int_{\R^\dm \times \R^\dm} \hspace{-0.7cm} \norm{\bfx - \bfy} \gamma(\dxx(\bfx, \bfy)) \\
        &= \br*{K \SmallSum{i = 0}{n - 1} (\eta L)^i}  \wdist_\dm \big( \tker^{\pr{a}}(x, \cdot), \tker_{\pr{a}}(y, \cdot) \big) \le  \br*{ \eta K \SmallSum{i = 0}{n - 1} (\eta L)^i } \norm{x - y}. \nonumber
    \end{align}
    This and induction prove \eqref{eq:lipschitz_continuity_of_fixed_point_iterates_kernel_w}.
    Furthermore, note that the assumption that $\eta L \in [0,1)$ ensures for all $n \in \N_0$ that $ \sum_{i = 0}^{n - 1} (\eta L)^i \le \sum_{i = 0}^{\infty} (\eta L)^i = \frac{1}{1 - \eta L}$.
    This yields that for all $n \in \N_0$, $x,y \in \R^\dm$, $a \in A$ it holds that 
    \begin{equation}
        \abs{ (u_n(x))(a) - (u_n(y))(a) } \le \tfrac{\eta K}{1 - \eta L} \|x - y\|.
    \end{equation}
    Combining this with the fact that $\lim_{n \to \infty} \big( \sup_{x \in \R^\dm, a \in A} \abs{\w(x)}^{-1} \abs{ (u_n(x))(a) - (u(x))(a) } \big) = 0$ establishes for all $x,y \in \R^\dm$, $a \in A$ that 
    \begin{equation}
        \abs{ (u(x))(a) - (u(y))(a) } = \lim_{n \rightarrow \infty} \abs{ (u_n(x))(a) - (u_n(y))(a) } \le \tfrac{\eta K}{1 - \eta L} \| x - y \|.
    \end{equation}
    This proves \cref{item:lipschitz_continuity_of_solutions_kernel_w_result}. The proof of \cref{lem:lipschitz_continuity_of_solutions_kernel_w} is thus complete.
\end{mproof}

\cfclear
\begin{lemma}\label{lem:stability_kernel_w_kernel}
    \cfconsiderloaded{lem:stability_kernel_w_kernel}
    Assume \cref{set:stability_w},
    let $\eta, c, L, K \in [0, \infty)$,
    assume $cL < 1$ and $\eta L < 1$,
    let $f \colon \R^\dm \times \R^{\abs{A}} \to \R$ be measurable,
    assume for all $x, y \in \R^\dm$, $r, s \in \R^{\abs{A}}$ that $\abs{ f(x, r) - f(x, s) } \le L \max_{a \in A} \abs{ r(a) - s(a)}$ and $\abs{ f(x, r) - f(y, r) } \le K \norm{x - y}$,
    let $\tker^{\pr{a}}_i \colon \R^\dm \times \Borel(\R^\dm) \to \br{0,1}$, $i \in \{1,2\}$, $a \in A$, be stochastic kernels,
    assume for all $i \in \{ 1,2 \}$, $x \in \R^\dm$, $a \in A$ that $ \int_{\R^\dm} \abs{\w(y)} \tker_i^{\pr{a}}(x, \dxx y) \le c \w(x)$
    and $\int_{\R^\dm} \norm{ y} \tker^{\pr{a}}_i(x, \dxx y) < \infty$,
    assume for all $i \in \cu{1,2}$ that
    \begin{equation}
        \sup_{(x, a) \in \R^\dm \times A} \abs{\w(x)}^{-1} \int_{\R^\dm} \abs{ f\pr{y, 0} } \tker_{i}^{\pr{a}}(x, \dxx y) < \infty,
    \end{equation}
    assume for all $x,y \in \R^\dm$, $a \in A$ that $\wdist_\dm\pr{\tker_2^{\pr{a}}(x, \cdot), \tker_2^{\pr{a}}(y, \cdot)} \le \eta \norm{x - y}$,
    Then
    \begin{enumerate}[label=(\roman *)]
        \item \label{item:stability_kernel_w_kernel_ex_uniq} there exist unique $u_i \in \W$, $i \in \cu{1,2}$, such that for every $i \in \cu{1,2}$, $x \in \R^\dm$, $a \in A$ it holds that
        \begin{equation}
            \int_{\R^\dm} \abs{f \pr{ y, u_i(y) }} \tker^{\pr{a}}_i \pr{x, \dxx y} < \infty \qquad \text{and} \qquad \pr{u_i\pr{x}}\pr{a} = \int_{\R^\dm} f\pr{ y, u_i(y) } \tker^{\pr{a}}_{i} \pr{x, \dxx y}
        \end{equation}
        and
        \item \label{item:stability_kernel_w_kernel_result} it holds that
        \begin{equation} \label{eq:stability_kernel_w_kernel}
            \br[\bigg]{\sup_{(x, a) \in \R^\dm \times A}\!\!\!\! \tfrac{\abs{ (u_1(x))(a) - (u_2(x))(a) }}{\abs{\w(x)}} } \le \tfrac{K}{(1 - \eta L) (1 - cL)} \br[\bigg]{ \sup_{(y, b) \in \R^\dm \times A} \!\!\!\! \tfrac{\wdist_\dm \pr{ \tker_{1}^{\pr{b}}(y, \cdot), \tker_{2}^{\pr{b}}(y, \cdot) }}{\abs{\w(y)}} }
        \end{equation}
    \end{enumerate}
    \cfout.
\end{lemma}

\begin{mproof}{\cref{lem:stability_kernel_w_kernel}}
    First, note that 
    the assumption that $\w$ is measurable,
    the assumption that $cL<1$,
    the assumption that for all $x \in \R^\dm$, $r,s \in \R^{\abs{A}}$ it holds that $\abs{f\pr{x,r} - f\pr{x,s}} \le L \max_{a \in A} \abs{ r\pr{a} - s\pr{a} }$,
    the assumption that for all $i \in \cu{1,2}$, $x \in \R^\dm$, $a \in A$ it holds that $\int_{\R^\dm} \abs{\w(y)} \tker_i^{\pr{a}} \pr{x, \dxx y} \le c \w(x)$, 
    the assumption that for all $i \in \cu{1,2}$ it holds that $\sup_{\pr{x,a} \in \R^\dm \times A} \abs{\w(x)}^{-1} \int_{\R^\dm} \abs{f\pr{y, 0}} \tker^{\pr{a}}_i\pr{x, \dxx y} < \infty$,
    and \cite[Lemma~2.2]{beck2023nonlinctrl} (applied for every $i \in \cu{1,2}$ with $c \with c$, $L \with L$, $\pr{\X, \mathcal{X}} \with \pr{\R^\dm, \cB\pr{\R^\dm}}$, $A\with A$, $\pr{\tker_a}_{a \in A} \with \pr{\tker_i^{\pr{a}}}_{a \in A}$, $f \with f$, $\w \with \br{ \R^\dm \ni x \mapsto \br{ A \ni a \mapsto \w(x) \in \pr{0,\infty} } \in \pr{0,\infty}^{\abs{A}} }$ in the notation of \cite[Lemma~2.2]{beck2023nonlinctrl}) prove \cref{item:stability_kernel_w_kernel_ex_uniq}.
    Next, observe that \cref{item:lipschitz_continuity_of_solutions_kernel_w_result} in \cref{lem:lipschitz_continuity_of_solutions_kernel_w} (applied with $\eta \with \eta$, $c \with c$, $L \with L$, $K \with K$, $f \with f$, $(\tker^{\pr{a}})_{a \in A} \with (\tker_2^{\pr{a}})_{a \in A}$ in the notation of \cref{lem:lipschitz_continuity_of_solutions_kernel_w}) yields for all $x, y \in \R^\dm$, $a \in A$ that
    \begin{equation}
        \abs{ (u_2(x))(a) - (u_2(y))(a) } \le \tfrac{\eta K}{1 - \eta L} \norm{x - y}.
    \end{equation}
    This,
    the triangle inequality,
    and the assumption that for all $x, y \in \R^\dm$, $r, s \in \R^{\abs{A}}$ it holds that $\abs{ f(x,r) - f(x,s) } \le L \max_{a \in A} \abs{ r(a) - s(a) }$ and $\abs{ f(x,r) - f(y, r) } \le K \norm{x - y}$
    demonstrate for all $x, y \in \R^\dm$ that 
    \begin{align}
        \begin{split}
            \abs*{ f(x, u_2(x)) - f(y, u_2(y)) } &\le \abs{ f(x, u_2(x)) - f(x, u_2(y)) } + \abs{f(x, u_2(y)) - f(y, u_2(y))}\\
            &\le L \max_{a \in A} \abs{ (u_2(x))(a) - (u_2(y))(a) } + K \norm{x - y} \\
            &\le \pr*{L \tfrac{\eta K}{1 - \eta L} + K} \norm{x - y} = \tfrac{K}{1 - \eta L} \norm{x - y}.
        \end{split}
    \end{align}
    Combining this,
    the fact that for all $x \in \R^\dm$, $a \in A$, $\gamma \in \coup_\dm(\tker_{1}^{\pr{a}}(x, \cdot), \tker_{2}^{\pr{a}}(x,\cdot))$, $B \in \Borel(\R^\dm)$ it holds that $\gamma(B \times \R^\dm) = \tker_{1}^{\pr{a}}(x, B)$ and $\gamma(\R^\dm \times B) = \tker_{2}^{\pr{a}}(x, B)$,
    the triangle inequality,
    the assumption that for all $x \in \R^\dm$, $r,s\in \R^{\abs{A}}$ it holds that $\abs{ f(x,r) - f(x,s) } \le L \max_{a \in A} \abs{r(a) - s(a)}$,
    and the assumption that for all $i \in \{1,2\}$, $x \in \R^\dm$ $a \in A$ it holds that $\int_{\R^\dm} \abs{\w(\bfx)} \tker^{\pr{a}}_i(x, \dxx \bfx) \le c \w(x)$
    with \cref{item:stability_kernel_w_kernel_ex_uniq}
    assures for all $x \in \R^\dm$, $a \in A$, $\gamma \in \coup_\dm(\tker_{1}^{\pr{a}}(x, \cdot), \tker_{2}^{\pr{a}}(x, \cdot))$ that
    \begin{align}
        \begin{split}
            \tfrac{\abs{ (u_1(x))(a) - (u_2(x))(a) }}{\abs{\w(x)}} &= \tfrac{1}{\abs{\w(x)}} \abs*{ \int_{\R^\dm} f\pr{\bfx, u_1(\bfx)} \tker_{1}^{\pr{a}}(x, \dxx \bfx) - \int_{\R^\dm} f\pr{\bfy, u_2(\bfy)} \tker_{2}^{\pr{a}}(x, \dxx \bfy) }  \\
            &= \tfrac{1}{\abs{\w(x)}} \abs*{ \int_{\R^\dm \times \R^\dm} \hspace{-0.7cm} f\pr{ \bfx, u_1(\bfx) } - f\pr{\bfy, u_2(\bfy)} \gamma (\dxx (\bfx, \bfy)) } \\
            &\le \tfrac{1}{\abs{\w(x)}} \int_{\R^\dm \times \R^\dm} \hspace{-0.5cm} \abs*{ f\pr{ \bfx, u_1(\bfx) } - f\pr{\bfy, u_2(\bfy)} } \gamma (\dxx (\bfx, \bfy)) \\
            &\le \tfrac{1}{\abs{\w(x)}}  \pr[\bigg]{\int_{\R^\dm \times \R^\dm} \hspace{-0.5cm} \abs{f\pr{ \bfx, u_1(\bfx) } - f\pr{\bfx, u_2(\bfx)}} \gamma(\dxx(\bfx, \bfy)) \\
            &\quad+ \int_{\R^\dm \times \R^\dm} \hspace{-0.5cm} \abs{f\pr{ \bfx, u_2(\bfx) } - f\pr{\bfy, u_2(\bfy)}} \gamma(\dxx(\bfx, \bfy))}  \\
            &\le \tfrac{1}{\abs{\w(x)}} \pr[\bigg]{L \int_{\R^\dm \times \R^\dm} \max_{b \in A} \abs{ (u_1(\bfx))(b) - (u_2(\bfx)(b) } \gamma(\dxx(\bfx, \bfy)) \\
            &\quad+ \tfrac{K}{1 - \eta L} \int_{\R^\dm \times \R^\dm} \hspace{-0.5cm} \norm{\bfx - \bfy} \gamma (\dxx (\bfx, \bfy))} \\
            &\le \tfrac{L}{\abs{\w(x)}} \br[\bigg]{ \sup_{(y,b) \in \R^\dm \times A} \!\!\!\! \tfrac{\abs{ (u_1(y))(b) - (u_2(y))(b) } }{ \abs{\w(y)} }} \int_{\R^\dm \times \R^\dm} \hspace{-0.5cm} \abs{\w(\bfx)} \gamma (\dxx (\bfx, \bfy)) \\
            &\quad+ \tfrac{K}{(1 - \eta L) \abs{\w(x)}} \int_{\R^\dm \times \R^\dm} \hspace{-0.5cm} \norm{\bfx - \bfy} \gamma (\dxx (\bfx, \bfy))\\
            &\le cL \br[\bigg]{ \sup_{(y,b) \in \R^\dm \times A} \!\!\!\! \tfrac{\abs{ (u_1(y))(b) - (u_2(y))(b) } }{ \abs{\w(y)} }} + \tfrac{K}{(1 - \eta L) \abs{\w(x)}} \int_{\R^\dm \times \R^\dm} \hspace{-0.5cm} \norm{\bfx - \bfy} \gamma (\dxx (\bfx, \bfy)).
        \end{split}
    \end{align}
    Hence, it holds for all $x \in \R^\dm$, $a \in A$ that
    \begin{align}
        \begin{split}
            &\tfrac{\abs{ (u_1(x))(a) - (u_2(x))(a) }}{\abs{\w(x)}} - cL \br[\bigg]{ \sup_{(y,b) \in \R^\dm \times A} \!\!\!\! \tfrac{\abs{ (u_1(y))(b) - (u_2(y))(b) } }{ \abs{\w(y)} }}\\
            &\quad\le \tfrac{K}{(1 - \eta L) \abs{\w(x)}} \pr*{\inf_{\gamma \in \coup_\dm( \tker_{1}^{\pr{a}}(x, \cdot), \tker_{2}^{\pr{a}}(x, \cdot ))} \int_{\R^\dm \times \R^\dm} \hspace{-0.5cm} \norm{\bfx - \bfy} \gamma (\dxx (\bfx, \bfy))} \\
            &\quad= \tfrac{K}{1 - \eta L} \tfrac{\wdist_\dm( \tker_{1}^{\pr{a}}(x, \cdot), \tker_{2}^{\pr{a}}(x, \cdot) )}{\abs{\w(x)}}.
        \end{split}
    \end{align}
    Combining this,
    the assumption that $cL < 1$,
    and \cref{item:stability_kernel_w_kernel_ex_uniq}
    shows that 
    \begin{equation}
        \br[\bigg]{\sup_{(x, a) \in \R^\dm \times A}\!\!\!\! \tfrac{\abs{ (u_1(x))(a) - (u_2(x))(a) }}{\w(x)} } \le \tfrac{K}{(1 - \eta L) (1 - cL)} \br[\bigg]{ \sup_{(y, b) \in \R^\dm \times A} \!\!\!\! \tfrac{\wdist_\dm \pr{ \tker_{1}^{\pr{b}}(y, \cdot), \tker_{2}^{b}(y, \cdot) }}{\abs{\w(y)}} }.
    \end{equation}
    This proves \cref{item:stability_kernel_w_kernel_result}. The proof of \cref{lem:stability_kernel_w_kernel} is thus complete.
\end{mproof}

\cfclear
\begin{propo} \label{prop:stability_kernel_w}
    \cfconsiderloaded{prop:stability_kernel_w}
    Assume \cref{set:stability_w},
    let $\eta, c, K \in [0, \infty)$,
    let $L_i \in [0, \infty)$, $i \in \{1,2\}$,
    assume for all $i \in \{1,2\}$ that $cL_i < 1$ and $\eta L_2 < 1$,
    let $f_i \colon \R^\dm \times \R^{\abs{A}} \to \R$, $i \in \{1, 2\}$, be measurable, 
    assume for all $i \in \{1,2\}$, $x, y \in \R^\dm$, $r, s \in \R^{\abs{A}}$ that $\abs{ f_i(x,r) - f_i(x,s) } \le L_i \max_{a \in A} \abs{ r(a) - s(a)}$ and $\abs{f_2(x,r) - f_2(y, r)} \le K \norm{x - y}$,
    let $\tker^{\pr{a}}_i \colon \R^\dm \times \Borel(\R^\dm) \to \br{0,1}$, $i \in \{1,2\}$, $a \in A$, be stochastic kernels,
    assume for all $i \in \{1,2\}$, $x \in \R^\dm$, $a \in A$ that $\int_{\R^\dm} \abs{\w(y)} \tker_i^{\pr{a}}(x, \dxx y) \le c \w(x)$,
    assume for all $i \in \cu{1,2}$ that
    \begin{equation}
        \begin{split}
            \sup_{(x, a) \in \R^\dm \times A} \abs{\w(x)}^{-1} \int_{\R^\dm} \abs{ f_i \pr{y, 0} } \tker_i^{\pr{a}}(x, \dxx y) < \infty,
        \end{split}
    \end{equation}
    assume for all $i \in \cu{1,2}$ that $\sup_{ (x, a) \in \R^\dm \times A } \abs{\w(x)}^{-1} \int_{\R^\dm} \norm{y} \tker_i^{\pr{a}}(x, \dxx y) < \infty $,
    assume for all $x, y \in \R^\dm$, $a \in A$ that $\wdist_\dm \pr{ \tker_2^{\pr{a}}(x, \cdot), \tker_2^{\pr{a}}(y, \cdot) } \le \eta \norm{x-y}$.
    Then
    \begin{enumerate}[label=(\roman *)]
        \item \label{item:stability_kernel_w_ex_uniq} there exist unique $u_i \in \W$, $i\in \cu{1,2}$, such that for all $i \in \cu{1,2}$, $x \in \R^\dm$, $a \in A$ it holds that
        \begin{equation}
            \begin{split}
                \int_{\R^\dm} \abs{ f_i\pr{ y, u_i(y)} } \tker_i^{\pr{a}}\pr{x, \dxx y} < \infty \qquad \text{and} \qquad \pr{u_i\pr{x}}\pr{a} = \int_{\R^\dm} f_i \pr{y, u_i\pr{y}}\tker_i^{\pr{a}}\pr{x, \dxx y}
            \end{split}
        \end{equation}
        and
        
        \item \label{item:stability_kernel_w_result} it holds that
        \begin{equation} \label{eq:stability_kernel_w}
            \begin{split}   
                \br[\bigg]{ \sup_{(x,a) \in \R^\dm \times A} \!\!\!\!\! \tfrac{\abs{ (u_1(x))(a) - (u_2(x))(a) }}{\abs{ \w(x) }} } &\le \tfrac{c}{1 - c(\min\{L_1, L_2))} \br[\bigg]{ \sup_{ (y,r) \in \R^\dm \times \R^{\abs{A}} } \!\!\!\!\! \tfrac{ \abs{ f_1(y,r) - f_2(y,r) } }{ \abs{\w(y)} } } \\
                &\quad+ \tfrac{K}{(1 - \eta L_2) (1 - cL_2)} \br[\bigg]{ \sup_{(y, b) \in \R^\dm \times A} \!\!\!\! \tfrac{\wdist_\dm \pr{ \tker_{1}^{\pr{b}}(y, \cdot), \tker_{2}^{\pr{b}}(y, \cdot) }}{\abs{\w(y)}} }
            \end{split} 
        \end{equation}
    \end{enumerate}
    \cfout.
\end{propo}

\begin{mproof}{\cref{prop:stability_kernel_w}}
    First, observe that
    the assumption that $\w$ is measurable,
    the assumption that for all $i \in \cu{1,2}$ it holds that $cL_i < 1$,
    the assumption that for all $i \in \cu{1,2}$, $x \in \R^\dm$, $r,s \in \R^{\abs{A}}$ it holds that $\abs{f_i\pr{ x, r} - f_i\pr{x,s}} \le L_i \max_{a \in A} \abs{ r\pr{a} - s\pr{a}}$,
    the assumption that for all $i \in \cu{1,2}$, $x \in \R^\dm$, $a \in A$ it holds that $\int_{\R^\dm} \abs{\w(y)} \tker_i^{\pr{a}}\pr{x, \dxx y} \le c \w(x)$,
    the assumption that for all $i \in \cu{1,2}$ it holds that $\sup_{\pr{x,a} \in \R^\dm \times A} \abs{\w(x)}^{-1} \int_{\R^\dm} \abs{ f_i \pr{y, u_i\pr{y}}} \tker_i^{\pr{a}}\pr{x, \dxx y} < \infty $,
    and \cite[Lemma~2.2]{beck2023nonlinctrl} (applied for every $i \in \cu{1,2}$ with $c \with c$, $L \with L_i$, $\pr{\X, \mathcal{X}} \with \pr{\R^\dm, \cB\pr{\R^\dm}}$, $A\with A$, $\pr{\tker_a}_{a \in A} \with \pr{\tker_i^{\pr{a}}}_{a \in A}$, $f \with f_i$, $\w \with \br{ \R^\dm \ni x \mapsto \br{A \ni a \mapsto \w(x) \in \pr{0,\infty}} \in \pr{0,\infty}^{\abs{A}} }$ in the notation of \cite[Lemma~2.2]{beck2023nonlinctrl}) establishes \cref{item:stability_kernel_w_ex_uniq}.
    Next, note that the assumption that for all $i \in \cu{1,2}$ it holds that $\sup_{\pr{x,a} \in \R^\dm \times A} \abs{\w(x)}^{- 1} \int_{\R^\dm} \norm{y} \tker_i^{\pr{a}}\pr{x, \dxx y} \le c \w(x)$
    and the fact that for all $x \in \R^\dm$ it holds that $\w(x) > 0$
    imply that for all $i \in \cu{1,2}$, $x \in \R^\dm$, $a \in A$ it holds that
    \begin{equation} \label{eq:stability_kernel_w_finite_first_moment}
        \begin{split}
            \int_{\R^\dm} \norm{y} \tker_i^{\pr{a}} \pr{x, \dxx y} &= \tfrac{\abs{\w(x)}}{\abs{\w(x)}} \int_{\R^\dm} \norm{y} \tker_i^{\pr{a}} \pr{x, \dxx y} \\
            &\le \abs{\w(x)} \br[\bigg]{\sup_{\pr{z,b} \in \R^\dm \times A} \!\!\!\!  \abs{\w(z)}^{-1} \!\! \int_{\R^\dm} \norm{y} \tker_i^{\pr{b}} \pr{z, \dxx y}} < \infty.
        \end{split}
    \end{equation}
    Moreover, note that Fubini's theorem,
    the triangle inequality, 
    and the assumption that for all $x, y \in \R^\dm$, $r \in \R^{\abs{A}}$ it holds that $\abs{ f_2(x,r) - f_2(y,r) } \le K \norm{x - y}$
    ensure for all $x \in \R^\dm$, $a \in A$ that
    \begin{equation} \label{eq:stability_kernel_w_cross_intregrability}
        \begin{split}
            \tfrac{1}{\abs{\w(x)}} \int_{\R^\dm} \abs{f_2\pr{\bfx, 0}} \tker_1^{\pr{a}}(x, \dxx \bfx) &= \tfrac{1}{\abs{\w(x)}} \int_{\R^\dm \times \R^\dm} \hspace{-0.5cm}\abs{ f_2\pr{\bfx, 0} } \pr[\big]{ \tker_1^{\pr{a}}(x, \dxx \bfx) \otimes \tker_2^{\pr{a}}(x, \dxx \bfy) } \\
            &\le \tfrac{1}{\abs{\w(x)}} \int_{\R^\dm \times \R^\dm} \hspace{-0.5cm} \abs{ f_2\pr{\bfx, 0} - f_2\pr{\bfy, 0} } \pr[\big]{ \tker_1^{\pr{a}}(x, \dxx \bfx) \otimes \tker_2^{\pr{a}}(x, \dxx \bfy) }  \\
            &\quad+ \tfrac{1}{\abs{\w(x)}} \int_{\R^\dm \times \R^\dm} \hspace{-0.5cm} \abs{ f_2\pr{\bfy, 0} } \pr[\big]{ \tker_1^{\pr{a}}(x, \dxx \bfx) \otimes \tker_2^{\pr{a}}(x, \dxx \bfy) } \\
            &\le \tfrac{K}{\abs{\w(x)}} \int_{\R^\dm \times \R^\dm} \hspace{-0.5cm}\norm{ \bfx - \bfy} \pr[\big]{ \tker_1^{\pr{a}}(x, \dxx \bfx) \otimes \tker_2^{\pr{a}}(x, \dxx \bfy) }\\
            &\quad+ \tfrac{1}{\abs{\w(x)}} \int_{\R^\dm \times \R^\dm} \hspace{-0.5cm} \abs{ f_2\pr{\bfy, 0} } \pr[\big]{ \tker_1^{\pr{a}}(x, \dxx \bfx) \otimes \tker_2^{\pr{a}}(x, \dxx \bfy) }\\
            &\le \tfrac{K}{\abs{\w(x)}} \int_{\R^\dm \times \R^\dm} \hspace{-0.5cm} \norm{\bfx} \pr[\big]{ \tker_1^{\pr{a}}(x, \dxx \bfx) \otimes \tker_2^{\pr{a}}(x, \dxx \bfy) }\\
            &\quad+ \tfrac{K}{\abs{\w(x)}} \int_{\R^\dm \times \R^\dm} \hspace{-0.5cm} \norm{\bfy} \pr[\big]{ \tker_1^{\pr{a}}(x, \dxx \bfx) \otimes \tker_2^{\pr{a}}(x, \dxx \bfy) } \\
            &\quad+ \tfrac{1}{\abs{\w(x)}} \int_{\R^\dm \times \R^\dm} \hspace{-0.5cm} \abs{ f_2\pr{\bfy, 0} } \pr[\big]{ \tker_1^{\pr{a}}(x, \dxx \bfx) \otimes \tker_2^{\pr{a}}(x, \dxx \bfy) }\\
            &= \tfrac{K}{\abs{\w(x)}} \int_{\R^\dm} \norm{\bfx} \tker_1^{\pr{a}}(x, \dxx \bfx) + \tfrac{K}{\abs{\w(x)}} \int_{\R^\dm} \norm{\bfy} \tker_2^{\pr{a}}(x, \dxx \bfy)\\
            &\quad+ \tfrac{1}{\abs{\w(x)}} \int_{\R^\dm} \abs{f_2\pr{\bfy, 0}} \tker_2^{\pr{a}} (x, \dxx \bfy)\\
            &\le K \pr*{ \br[\bigg]{\sup_{(y,b) \in \R^\dm \times A} \!\!\!\! \tfrac{\int_{\R^\dm} \norm{\bfx} \tker_1^{\pr{b}}(y, \dxx \bfx) }{ \abs{\w(y)} } } + \br[\bigg]{\sup_{(y,b) \in \R^\dm \times A} \!\!\!\! \tfrac{\int_{\R^\dm} \norm{\bfy} \tker_2^{\pr{b}} (y, \dxx \bfy)}{ \abs{\w(y)} } }  }\\
            &\quad+ \br[\bigg]{ \sup_{(y,b) \in \R^\dm \times A} \!\!\!\! \tfrac{ \int_{\R^\dm} \abs{f_2 \pr{ \bfy, 0 } } \tker_2^{\pr{b}}(y, \dxx \bfy) }{ \abs{\w(y)} } }.
        \end{split}
    \end{equation}
    This,
    the assumption that for all $i \in \{1, 2\}$ it holds that $\sup_{(y,b) \in \R^\dm \times A} \abs{\w(y)}^{-1} \int_{\R^\dm} \abs{f_i\pr{\bfy, 0}} \tker_i^{\pr{b}}(y, \dxx \bfy) <\infty$,
    and the assumption that for all $i \in \{1,2\}$ it holds that $\sup_{(y,b) \in \R^\dm \times A} \abs{\w(y)}^{-1} \int_{\R^\dm} \norm{\bfy} \tker_i^{\pr{b}} (y, \dxx \bfy) <\infty$,
    show that $\sup_{(y, b) \in \R^\dm \times A} \abs{\w(y)}^{-1} \int_{\R^\dm} \abs{f_2\pr{\bfy, 0}} \tker_1^{\pr{b}}(y, \dxx \bfy) < \infty$.
    Combining this
    with \cite[Lemma~2.2]{beck2023nonlinctrl} (applied with $c \with c$, $L \with L_2$, $(\X, \mathcal{X}) \with (\R^\dm, \Borel(\R^\dm))$, $A \with A$, $(\tker_a)_{a \in A} \with (\tker^{\pr{a}}_1)_{a \in A}$, $f \with f_2$, $\w \with \br{ \R^\dm \ni x \mapsto \br{ A \ni a \mapsto \w(x) \in \pr{0, \infty}  } \in \pr{0, \infty}^{\abs{A}}}$ in the notation of \cite[Lemma~2.2]{beck2023nonlinctrl})
    proves that there exists a unique $u \in \W$ such that for all $x \in \R^\dm$, $a \in A$ it holds that
    $\int_{\R^\dm} \abs{ f_2 \pr{y, u(y)} } \tker_1^{\pr{a}}(x, \dxx y) < \infty$ and $(u(x))(a) = \int_{\R^\dm}  f_2\pr{y, u(y)}  \tker_1^{\pr{a}}(x, \dxx y)$.
    This,
    \cref{eq:stability_kernel_w_finite_first_moment},
    the triangle inequality,
    \cref{item:stability_kernel_w_nonlinearity_result} in \cref{lem:stability_kernel_w_nonlinearity} (applied with
    $c \with c$,
    $L_1 \with L_1$,
    $L_2 \with L_2$,
    $f_1 \with f_1$,
    $f_2 \with f_2$,
    $\pr{\tker^{\pr{a}}}_{a \in A} \with \pr{\tker^{\pr{a}}_1}_{a \in A}$
    in the notation of \cref{lem:stability_kernel_w_nonlinearity}),
    and \cref{item:stability_kernel_w_kernel_result} \cref{lem:stability_kernel_w_kernel} (applied with
    $\eta \with \eta$,
    $c \with c$,
    $L \with L_2$,
    $K \with K$,
    $f \with f_2$,
    $\pr{ \tker^{\pr{a}}_1 }_{a \in A} \with \pr{\tker^{\pr{a}}_1}_{a \in A}$,
    $\pr{ \tker^{\pr{a}}_2 }_{a \in A} \with \pr{\tker^{\pr{a}}_2}_{a \in A}$
    in the notation of \cref{lem:stability_kernel_w_kernel})
    demonstrate that
    \begin{equation}
        \begin{split}
            \br[\bigg]{ \sup_{(x,a) \in \R^\dm \times A} \!\!\!\!\!\! \tfrac{\abs{ (u_1(x))(a) - (u_2(x))(a) }}{\abs{ \w(x) }} } &\le \br[\bigg]{ \sup_{(x,a) \in \R^\dm \times A} \!\!\!\!\! \tfrac{\abs{ (u_1(x))(a) - (u(x))(a) }}{\abs{ \w(x) }} } + \br[\bigg]{ \sup_{(x,a) \in \R^\dm \times A} \!\!\!\!\! \tfrac{\abs{ (u(x))(a) - (u_2(x))(a) }}{\abs{ \w(x) }} } \\
            &\le \tfrac{c}{1 - c(\min\cu{L_1, L_2})} \br[\bigg]{ \sup_{ (y,r) \in \R^\dm \times \R^{\abs{A}} } \!\!\!\!\! \tfrac{ \abs{ f_1(y,r) - f_2(y,r) } }{ \abs{\w(y)} } } \\
            &\quad+ \tfrac{K}{(1 - \eta L_2) (1 - cL_2)} \br[\bigg]{ \sup_{(y, b) \in \R^\dm \times A} \!\!\!\!\! \tfrac{\wdist_\dm \pr{ \tker_{1}^{\pr{b}}(y, \cdot), \tker_{2}^{\pr{b}}(y, \cdot) }}{\abs{\w(y)}} }.
        \end{split}
    \end{equation}
    This proves \cref{item:stability_kernel_w_result}. The proof of \cref{prop:stability_kernel_w} is thus complete.
\end{mproof}

\section{Artificial neural network (ANN) calculus}
\label{sec:ann_calculus} 

In this section we recall necessary notions and statements on ANNs which we need for the remainder of this work. The content of this section is well-known in the literature and the specific definitions and statements mostly consist of slightly extended or modified excerpts from Grohs et al. \cite[Section~2]{GrohsHornungJentzen2019published}, Grohs et al. \cite[Section~3]{grohs2019deep_published}, and Ackermann et al. \cite[Section~2]{ackermann2023deep}.

In particular,
\cref{def:ANN} is a slightly extended version of \cite[Definition~2.1]{GrohsHornungJentzen2019published},
\cref{def:multi} is \cite[Definition~2.2]{GrohsHornungJentzen2019published},
\cref{def:realization} is \cite[Definition~2.3]{GrohsHornungJentzen2019published},
\cref{def:compositions_of_anns} is \cite[Definition~2.5]{GrohsHornungJentzen2019published},
\cref{lem:compositions_of_anns} is a version of \cite[Proposition~2.6]{GrohsHornungJentzen2019published},
\cref{lem:associativity_compositions_of_anns} is \cite[Lemma~2.8]{GrohsHornungJentzen2019published},
\cref{def:powers_of_anns} is \cite[Definition~2.9]{GrohsHornungJentzen2019published},
\cref{def:extensions_of_anns} is \cite[Definition~2.11]{GrohsHornungJentzen2019published},
\cref{lem:certain_extensions_of_anns} is a slightly extended version of \cite[Lemma~2.13]{GrohsHornungJentzen2019published} where the proof of the extension is provided for the convenience of the reader,
\cref{lem:extensions_of_anns_involving_identities} is \cite[Lemma~2.14]{GrohsHornungJentzen2019published},
\cref{def:parallelizations_of_anns_same_depth} is a slightly modified version of \cite[Definition~2.17]{GrohsHornungJentzen2019published},
\cref{lem:parallelizations_of_anns_same_depth} is a collection of \cite[Lemma~2.18, Proposition~2.19, and Proposition~2.20]{GrohsHornungJentzen2019published},
\cref{def:parallelizations_of_anns_with_different_layer_structure} is \cite[Definition~2.22]{GrohsHornungJentzen2019published},
\cref{cor:parallelizations_of_anns_with_different_layer_structure_realization_involving_identities} is a slightly extended version of \cite[Corollary~2.23]{GrohsHornungJentzen2019published} where the proof of the extension is provided for the convenience of the reader,
\cref{def:sum_transpose_anns} is \cite[Definitions~2.13 and 2.15]{ackermann2023deep},
\cref{def:sum_of_anns_same_length} is \cite[Definition~2.16]{ackermann2023deep},
\cref{def:sum_of_anns_diff_length} is \cite[Definition~2.17]{ackermann2023deep},
\cref{def:scalar_multiplications_of_anns} is \cite[Definition~2.18]{ackermann2023deep},
\cref{lem:scalar_multiplications_of_anns} is \cite[Lemma~3.14]{grohs2019deep_published},
\cref{lem:linear_combinations_anns_same_length} is a slightly modified version of \cite[Lemma~2.19]{ackermann2023deep} where the proof is provided for the convenience of the reader,
and \cref{lem:linear_combination_anns_different_length} is a slightly modified version of \cite[Lemma~2.20]{ackermann2023deep} where the proof is provided for the convenience of the reader.

\subsection{Structured description of ANNs}
\label{subsec:ann_description}

\begin{definition}[ANNs] \label{def:ANN}
    We denote by $\bN$ the set given by
    \begin{equation} \label{eq:ANN_set}
        \textstyle
        \bN =
        \textstyle
        \bigcup_{L \in \N} \bigcup_{l_0, l_1, \dots, l_L \in \N} \pr[\big]{\bigtimes_{k = 1}^{L}   \pr{\R^{l_k \times l_{k - 1}} \times \R^{l_k}}}
    \end{equation}
    and we denote by
    $\paramANN \colon \bN \to \N$,
    $\cL \colon \bN \to \N$,
    $\cI \colon \bN \to \N$,
    $\cO \colon \bN \rightarrow \N$,
    $\cH \colon \bN \rightarrow \N_0 = \N \cup \{ 0 \}$,
    $\cD \colon \bN \rightarrow \cup_{l = 2}^\infty \N^l$,
    and $\dimANNlevel_n \colon \bN \rightarrow \N_0$, $n\in \N_0$,
    the functions which satisfy for all $L \in \N$, $l_0, l_1, \dots, l_L \in \N$, $\Phi \in \pr[\big]{\bigtimes_{k = 1}^{L} ( \R^{l_k \times l_{k - 1}} \times \R^{l_k} )} $, $n \in \N_0$ that $ \paramANN(\Phi) = \sum_{k = 1}^{L} l_k(l_{k - 1} + 1) $, $\cL(\Phi) = L$, $\cI(\Phi) = l_0$, $\cO(\Phi) = l_L$, $\cH(\Phi) = L - 1$, $\cD(\Phi) = (l_0, l_1, \dots, l_L)$, and
    \begin{equation} \label{eq:ANN_nth_layer}
        \dimANNlevel_n(\Phi) = \begin{cases}
            l_n
            &
            \colon
            n \le L
            \\
            0
            &
            \colon
            n > L.
        \end{cases}
    \end{equation}
\end{definition}

\cfclear
\begin{definition}[ANN] \label{def:ANN2}
    We say that $\Phi$ is an ANN if and only if it holds that $\Phi \in \bN$
    \cfout[.]
\end{definition}

\cfclear
\begin{lemma}\label{lem:ANN_size_estimate}
    \cfconsiderloaded{lem:ANN_size_estimate}
    For every $\dm \in \N$, $x \in \pr{x_1, \dots, x_\dm} \in \R^\dm$ let $\fnorm{x} \in \R$ satisfy $\fnorm{x} = \max_{i \in \cu{1,\dots, \dm}} \abs{x_i}$.
    Then it holds for all $\Phi \in \bN$ that $\max\{ \cL(\Phi), \fnorm{\cD(\Phi)} \} \le \paramANN(\Phi)$ and $\paramANN(\Phi) \le 2 \cL(\Phi) \fnorm{\cD(\Phi)}^2$
    \cfout[.]
\end{lemma}

\begin{mproof}{\cref{lem:ANN_size_estimate}}
    First, note that \cite[Lemma~2.4]{ackermann2023deep} demonstrates for all $\Phi \in \bN$ that $\cL(\Phi) + \fnorm{\cD(\Phi)} \le \paramANN(\Phi)$.
    Hence, it holds for all $\Phi \in \bN$ that $\max\{ \cL(\Phi), \fnorm{\cD(\Phi)} \} \le \paramANN(\Phi)$.
    Moreover, observe that for all $L \in \N$, $l_0, l_1, \dots, l_L \in \N$, $\Phi \in \pr{\bigtimes_{k = 1}^L \pr{ \R^{l_k \times l_{k - 1}} \times \R^{l_k} }}$ it holds that
    \begin{equation}
        \begin{split}
            \paramANN\pr{\Phi} = \sum_{k = 1}^L l_k\pr{l_{k - 1} + 1} \le \sum_{k = 1}^{L} \fnorm{\cD\pr{\Phi}} \pr{ \fnorm{\cD\pr{ \Phi }} + 1} \le 2L \fnorm{\cD \pr{ \Phi}}^2 = 2 \cL \pr{\Phi} \fnorm{\cD \pr{\Phi}}^2.
        \end{split}
    \end{equation}
    The proof of \cref{lem:ANN_size_estimate} is thus complete.
\end{mproof}

\begin{definition}[Multidimensional version] \label{def:multi}
    Let $\dm \in \N$, $\activation \in C(\R, \R)$.
    Then we denote by $\activationDim{\dm} \colon \R^\dm \to \R^\dm$ the function which satisfies for all $x = (x_1, x_2, \dots, x_\dm) \in \R^\dm$ that
    \begin{equation}
        \activationDim{\dm} (x) = \big( \activation(x_1), \activation(x_2), \dots, \activation(x_\dm) \big).
    \end{equation} 
\end{definition}

\cfclear
\begin{definition}[Realization associated to an ANN]\label{def:realization}
    \cfconsiderloaded{def:realization}
    Let $\activation \in C(\R, \R)$.
    Then we denote by $\cR_\activation \colon \bN \to  \cup_{k,l \in \N} C(\R^k, \R^l)$ the function which satisfies for all $L \in \N$, $l_0, l_1, \dots, l_L \in \N$, $\Phi = ( (W_1, B_1), (W_2, B_2), \dots, (W_L, B_L) ) \in \pr[\big]{\bigtimes_{k = 1}^L  (\R^{l_k \times l_{k - 1}} \times \R^{l_k})}$, $x_0 \in \R^{l_0}$, $x_1 \in \R^{l_1}$, \dots, $x_{L - 1} \in \R^{l_{L - 1}}$ with $x_k = \activationDim{l_k} (W_k x_{k - 1} + B_k)$, $k \in \N_0 \cap [1, L-1]$, that
    \begin{equation}
        \cR_\activation(\Phi) \in C(\R^{l_0}, \R^{l_L}), \quad \text{and} \quad \big(\cR_\activation(\Phi)\big)(x_0) = W_L x_{L - 1} + B_L.
    \end{equation}
    \cfout.
\end{definition}

\subsection{Compositions, affine linear transformations, powers, and extensions of ANNs}

\cfclear
\begin{definition}[Compositions of ANNs] \label{def:compositions_of_anns}
    \cfconsiderloaded{def:compositions_of_anns}
    We denote by $ \ANNcomp \colon \{ (\Phi_1, \Phi_2) \in \bN \times \bN : \cI(\Phi_1) = \cO(\Phi_2) \} \rightarrow \bN$ the function which satisfies for all $K, L \in \N$, $k_0, k_1, \dots, k_K, l_0, l_1, \dots, l_L \in \N$,
    $\Phi_1 = \pr{(W_1, B_1), \dots, (W_K, B_K)} \in  \pr[\big]{\bigtimes_{j = 1}^{K}  \pr{\R^{k_j \times k_{j - 1}}  \times \R^{k_j}} }$,
    $\Phi_2 = \pr{(\mathscr{W}_1, \mathscr{B}_1), \dots, (\mathscr{W}_L, \mathscr{B}_L)} \in \pr[\big]{\bigtimes_{j = 1}^{L} \pr{\R^{ l_j \times l_{j - 1} } \times \R^{l_j}}}$
    with $ k_0 = \cI(\Phi_1) = \cO(\Phi_2) = l_L$ that
    \begin{align} \label{eq:definition_compositions_of_anns}
        \begin{split}
            &\Phi_1 \bullet \Phi_2 = \\
            &\begin{cases}
                \pr[\big]{(W_1 \mathscr{W}_1, W_1\mathscr{B}_1 + B_1 )}
                &\colon K = 1 = L\\[1ex]
                \pr[\big]{(W_1 \mathscr{W}_1, W_1  \mathscr{B}_1 + B_1), (W_2, B_2), \dots, (W_K, B_K)} 
                &\colon K > 1 = L  \\[1ex]
                \pr[\big]{(\mathscr{W}_1, \mathscr{B}_1),(\mathscr{W}_2, \mathscr{B}_2), \dots (\mathscr{W}_{L - 1}, \mathscr{B}_{L - 1}), (W_1 \mathscr{W}_L, W_1 \mathscr{B}_L + B_1)}
                &\colon K = 1 < L \\[1ex]
                \begin{aligned}
                    &\pr[\big]{(\mathscr{W}_1, \mathscr{B}_1), (\mathscr{W}_2, \mathscr{B}_2), \dots, (\mathscr{W}_{L - 1}, \mathscr{B}_{L - 1}), (W_1 \mathscr{W}_L, W_1 \mathscr{B}_L + B_1), \\ &\hspace{0.2cm}(W_2, B_2), \dots, (W_K, B_K)}
                \end{aligned}
                &\colon K > 1 < L
            \end{cases}
        \end{split}
    \end{align}
    \cfout.
\end{definition}

\cfclear
\begin{lemma}[Compositions of ANNs] \label{lem:compositions_of_anns}
    \cfconsiderloaded{\cref{lem:compositions_of_anns}}
    Let $K, L \in \N$, $k_0, k_1, \dots, k_K, l_0, l_1, \dots, l_L \in \N$,
    let $\Phi_1 = \pr{(W_1, B_1), \dots, (W_K, B_K)}  \in  \pr[\big]{\bigtimes_{j = 1}^{K}  \pr{\R^{k_j \times k_{j - 1}}  \times \R^{k_j}}}$, $ \Phi_2 =  \pr{(\mathscr{W}_1, \mathscr{B}_1), \dots, (\mathscr{W}_L, \mathscr{B}_L)} \in \pr[\big]{ \bigtimes_{j = 1}^{L}  \pr{\R^{ l_j \times l_{j - 1} } \times \R^{l_j} }} $ satisfy $ k_0 = \cI(\Phi_1) = \cO(\Phi_2) = l_L$.
    Then
    \begin{enumerate}[label=(\roman *)]
        \item \label{item:compositions_of_anns_layer_structure} it holds that $ \cD ( \Phi_1 \bullet \Phi_2 ) = (l_0, \dots, l_{L - 1}, k_1, \dots, k_K) \in \N^{L + K} $,
        
        \item \label{item:compositions_of_anns_depth} it holds that $ \cL ( \Phi_1 \bullet \Phi_2) = \cL(\Phi_1) + \cL(\Phi_2) - 1 $,
        
        \item \label{item:compositions_of_anns_hidden_layers} it holds that $ \cH (\Phi_1 \bullet \Phi_2) = \cH(\Phi_1) + \cH(\Phi_2)$,
        
        \item \label{item:compositions_of_anns_size} it holds that 
        \begin{align}
            \begin{split}
                \paramANN ( \Phi_1 \bullet \Phi_2 ) &= \paramANN(\Phi_1) + \paramANN ( \Phi_2 ) + k_1 ( l_{L - 1} + 1 ) - k_1(k_0 + 1) - l_L ( l_{L-  1} + 1 )\\
                &= \paramANN(\Phi_1) + \paramANN(\Phi_2) + k_1 l_{L - 1} - k_1 k_0 - l_Ll_{L - 1} - l_L \\
                &\le \paramANN(\Phi_1) + \paramANN(\Phi_2) + k_1 l_{L - 1},
            \end{split}
        \end{align}
        and

        \item \label{item:compositions_of_anns_realizations} it holds for all $\activation \in C(\R, \R)$ that $ \cR_\activation( \Phi_1 \ANNcomp \Phi_2 ) \in C ( \R^{\cI(\Phi_2)  }, \R^{ \cO (\Phi_1) } ) $ and $ \cR_\activation ( \Phi_1 \bullet \Phi_2 ) = \cR_\activation ( \Phi_1 ) \circ \cR_\activation ( \Phi_2 )$
    \end{enumerate}
    %
    %
    \cfout.
\end{lemma}

\cfclear
\begin{lemma}[Associativity of compositions of ANNs] \label{lem:associativity_compositions_of_anns}
    \cfconsiderloaded{\lem:associativity_compositions_of_anns}
    Let $\Phi_1, \Phi_2, \Phi_3 \in \bN$ satisfy that $\cI (\Phi_1) = \cO (\Phi_2)$ and $\cI (\Phi_2) = \cO (\Phi_3)$. Then it holds that 
    \begin{equation}
        \left( \Phi_1 \ANNcomp \Phi_2 \right) \ANNcomp \Phi_3 = \Phi_1 \ANNcomp \left( \Phi_2 \ANNcomp \Phi_3 \right)
    \end{equation}
    \cfload.
\end{lemma}

\cfclear
\begin{definition}[Identity matrix] \label{def:identity_matrix}
    \cfconsiderloaded{def:identity_matrix}
    Let $\dm \in \N$. Then we denote by $\id_\dm \in \R^{\dm \times \dm}$ the identity matrix in $\R^{\dm \times \dm}$
    \cfout. 
\end{definition}

\cfclear
\begin{definition}[Powers of ANNs] \label{def:powers_of_anns}
    %
    %
    %
    %
    We denote by $(\cdot)^{\bullet n} \colon \{ \Phi \in \bN : \cI(\Phi) = \cO(\Phi) \} \rightarrow \bN$, $n \in \N_0$, the functions which satisfy for all $n \in \N_0$, $\Phi \in \bN$ with $\cI(\Phi) = \cO(\Phi)$ that
    \begin{equation} \label{eq:def_powers_of_anns}
        \Phi^{\bullet n } = \begin{cases}
            \left( \id_{\cO(\Phi)}, (0, 0, \dots, 0) \right) \in \left( \R^{\cO(\Phi) \times \cO(\Phi)} \times \R^{\cO(\Phi)} \right)
            &\colon
            n = 0
            \\
            \Phi \ANNcomp (\Phi^{\bullet (n - 1)})
            &\colon
            n \in \N
        \end{cases}
    \end{equation}
    \cfout.
\end{definition}

\cfclear
\begin{definition}[Extensions of ANNs] \label{def:extensions_of_anns}
    \cfconsiderloaded{def:extensions_of_anns}
    Let $L\in \N$, let $\Psi \in \bN$ with $\cI (\Psi) = \cO (\Psi)$. Then we denote by $\longerANN{L, \Psi} \colon \{ \Phi \in \bN : ( \cL (\Phi) \le L \text{ and } \cO (\Phi) = \cI (\Psi)  ) \} \rightarrow \bN$ satisfy for all $\Phi \in \bN$ with $ \cL (\Phi) \le L $ and $ \cI (\Psi) = \cO (\Phi) $ that $ \mathcal{E}_{L, \Psi} (\Phi) = \compANN{ (\power{\Psi}{(L - \cL (\Phi))}) }{\Phi}$
    \cfout.
\end{definition}

\cfclear
\begin{lemma}[Certain powers and extensions of ANNs] \label{lem:certain_extensions_of_anns}
    Let $\dm, l \in \N$, let $\Psi \in \bN$ satisfy $ \cD (\Psi) = (\dm , l , \dm) \in \N^3 $. Then
    \begin{enumerate}[label=(\roman *)]
        \item \label{item:certain_extensions_of_anns_item_power} it holds for all $n \in \N_0$ that $ \cL (\power{\Psi}{n}) = n + 1 $, $ \cD (\Psi^{\bullet n}) \in \N^{n + 2} $, and
        \begin{equation}
            \cD (\Psi^{ \bullet n }) = \begin{cases}
                (\dm, \dm) &: n = 0 \\
                (\dm, l, l, \dots, l, \dm) &: n \in \N,
            \end{cases}
        \end{equation}

        \item \label{item:certain_extensions_of_anns_item_extension} it holds for all $\Phi \in \bN$, $L \in \N \cap [\cL (\Phi), \infty)$ with $\cO (\Phi) = \dm$ that $ \cL (\mathcal{E}_{L, \Psi} (\Phi) ) = L $ and
        \begin{equation}
            \hspace{-0.42cm}\mathcal{P} ( \mathcal{E}_{L, \Psi} (\Phi) ) \le \begin{cases}
                \mathcal{P} (\Phi) &: \cL (\Phi) = L \\
                \left[ \left( \max \{ 1, \frac{l}{\dm} \} \right) \mathcal{P} (\Phi) + \left( \left( L - \cL (\Phi) - 1 \right) l + \dm \right)(l + 1)  \right] &: \cL (\Phi) < L,
            \end{cases}
        \end{equation} and 

        \item \label{item:certain_extensions_of_anns_item_extension_to_same_length} it holds for all $\Phi \in \bN$ with $\cO (\Phi) = \dm$ that $\longerANN{\cL(\Phi), \Psi} (\Phi) = \Phi$
    \end{enumerate}
    \cfout.
\end{lemma}

\begin{mproof}{\cref{lem:certain_extensions_of_anns}}
    First, note that \cite[Lemma~2.13]{GrohsHornungJentzen2019published} proves \cref{item:certain_extensions_of_anns_item_power,item:certain_extensions_of_anns_item_extension}.
    Next, observe that it holds that $\longerANN{\cL(\Phi), \Psi}(\Phi) = (\Psi^{\bullet (\cL(\Psi) - \cL(\Phi))}) \bullet \Phi = \Psi^{\bullet 0} \bullet \Phi$.
    This,
    the assumption that $\cO (\Phi) = \cI (\Psi) = \dm$,
    the fact that $\Psi^{\bullet 0} = ( \id_\dm, (0, 0, \dots, 0) ) \in \R^{\dm \times \dm} \times \R^\dm$,
    and \cref{lem:compositions_of_anns} demonstrate that $\longerANN{\cL (\Phi), \Psi} (\Phi) = \Psi^{\bullet 0} \bullet \Phi = \Phi$.
    This proves \cref{item:certain_extensions_of_anns_item_extension_to_same_length}.
    The proof of \cref{lem:certain_extensions_of_anns} is thus complete.
\end{mproof}

\cfclear
\begin{lemma}[Extensions of ANNs involving identities] \label{lem:extensions_of_anns_involving_identities}
    Let $\activation \in C(\R, \R)$,
    let $\mathbb{I} \in \bN$ satisfy for all $x \in \R^{ \cI (\mathbb{I}) }$ that $\cI (\mathbb{I}) = \cO (\mathbb{I})$ and $ \big(\cR_\activation(\mathbb{I})\big)(x) = x $.
    Then
    \begin{enumerate}[label=(\roman *)]
        \item \label{item:extensions_of_anns_involving_identities_item_powers_of_identity} it holds for all $n \in \N_0$, $ x \in \R^{ \cI (\mathbb{I}) } $ that $ \cR_\activation ( \power{\mathbb{I}}{n} ) \in C \left( \R^{ \cI (\mathbb{I}) }, \R^{ \cI (\mathbb{I}) } \right) $ and $ \big(\cR_\activation (\mathbb{I}^{\bullet n})\big)( x ) = x$ and
        
        \item \label{item:extensions_of_anns_involving_identities_item_extension_by_identity} it holds for all $\Phi \in \bN$, $L \in \N \cap [\cL  (\Phi), \infty)$, $x \in \R^{ \cI (\Phi) }$ with $ \cI (\mathbb{I}) = \cO (\Phi)$ that $ \cR_\activation (\longerANN{L, \mathbb{I}}(\Phi) ) \in C (\R^{ \cI (\Phi) }, \R^{ \cO (\Phi)  } )$ and $ \big(\cR_\activation (\mathcal{E}_{ L, \mathbb{I}})  (\Phi)  \big)(x) = \big(\cR_\activation(\Phi)\big)(x) $
    \end{enumerate}
    \cfout.
\end{lemma}

\subsection{Parallelizations of ANNs}
\label{subsec:parallelizations_of_anns}

\cfclear
\begin{definition}[Parallelizations of ANNs with the same depth] \label{def:parallelizations_of_anns_same_depth}
    \cfconsiderloaded{def:parallelizations_of_anns_same_depth}
    Let $n \in \N$. Then we denote by $\paraANN{n} \colon \{ (\Phi_1, \dots, \Phi_n) \in \bN^n : \cL (\Phi_1) = \cL (\Phi_2) = \dots = \cL (\Phi_n)  \} \rightarrow \bN $ the function which satisfies for all $L \in \N$, $ (l_{1,0}, l_{1,1}, \dots,  l_{1, L}), \dots, ( l_{n, 0}, l_{n, 1}, \dots, l_{n, L} ) \in \N^{L + 1} $,
    \begin{align}
        \begin{split}
            \Phi_1 &= ( (W_{1,1}, B_{1,1}), (W_{1,2}, B_{1,2}), \dots, (W_{1,L}, B_{1,L})  ) \in \textstyle \pr[\big]{\bigtimes_{k = 1}^{L}  ( \R^{ l_{1,k} \times l_{1, k- 1} }  \times \R^{l_{1,k}} ) } ,\\
            \Phi_2 &= ((W_{2,1}, B_{2,1}), (W_{2,2}, B_{2,2}), \dots, (W_{2,L}, B_{2,L})) \in \textstyle \pr[\big]{\bigtimes_{k = 1}^{L}  ( \R^{ l_{2,k} \times l_{2, k- 1} }  \times \R^{l_{2,k}} )  }, \\
            &\vdots\\
            \Phi_n &= ((W_{n,1}, B_{n,1}), (W_{n,2}, B_{n,2}), \dots, (W_{n,L}, B_{n,L})) \in \textstyle \pr[\big]{\bigtimes_{k = 1}^{L}  ( \R^{ l_{n,k} \times l_{n, k- 1} }  \times \R^{l_{n,k}} )}  ,
        \end{split}
    \end{align}
    \noindent
    $\pr*{ (\mathscr{W}_1, \mathscr{B}_1), (\mathscr{W}_2, \mathscr{B}_2), \dots, ( \mathscr{W}_L, \mathscr{B}_L )}  \in  \pr[\big]{\bigtimes_{j = 1}^{L} ( \R^{  (\sum_{k = 1}^{n} l_{k, j} )  \times  ( \sum_{k = 1}^{n} l_{k, j - 1} )}  \times \R^{ (\sum_{k = 1}^{n} l_{k, j} ) } )}  $ with 
    \begin{equation} \label{eq:parallelizations_of_anns_same_depth_def_matrix}
        \forall i \in \N \cap \in [1, L]: ( \mathscr{W}_i, \mathscr{B}_i ) = \begin{pmatrix}
            \begin{pmatrix}
                W_{1, i} & 0 & 0 & \hdots & 0\\
                0 & W_{2, i} & 0 & \hdots & 0\\
                0 & 0 & W_{3, i} & \hdots & 0\\
                \vdots & \vdots & \vdots & \ddots & \vdots \\
                0 & 0 & 0 & \hdots & W_{n, i}
            \end{pmatrix},
            \begin{pmatrix}
                B_{1, i} \\ B_{2, i} \\ B_{3, i} \\ \vdots \\ B_{n, i}
            \end{pmatrix}
        \end{pmatrix},
    \end{equation}
    that
    \begin{align}\label{eq:parallelizations_of_anns_same_depth_def_p_n}
        \begin{split}
            \paraANN{n} ( \Phi_1, \Phi_2, \dots \Phi_n  ) = \pr*{(\mathscr{W}_1, \mathscr{B}_1), (\mathscr{W}_2, \mathscr{B}_2), \dots, ( \mathscr{W}_L, \mathscr{B}_L ) }
        \end{split}
    \end{align}
    \cfout.
\end{definition}

\cfclear
\begin{lemma}\label{lem:parallelizations_of_anns_same_depth} 
    Let $n, L \in \N$, let $ (l_{1,0}, l_{1,1}, \dots,  l_{1, L}), \dots, ( l_{n, 0}, l_{n, 1}, \dots, l_{n, L} ) \in \N^{L + 1} $, and let 
    \begin{align}
        \begin{split}
            \Phi_1 &= ( (W_{1,1}, B_{1,1}), (W_{1,2}, B_{1,2}), \dots, (W_{1,L}, B_{1,L})  ) \in \textstyle \pr[\big]{ \bigtimes_{k = 1}^{L}  ( \R^{ l_{1,k} \times l_{1, k- 1} }  \times \R^{l_{1,k}} )} ,\\
            \Phi_2 &= ((W_{2,1}, B_{2,1}), (W_{2,2}, B_{2,2}), \dots, (W_{2,L}, B_{2,L})) \in \textstyle \pr[\big]{\bigtimes_{k = 1}^{L}  ( \R^{ l_{2,k} \times l_{2, k- 1} }  \times \R^{l_{2,k}} )}, \\
            &\vdots\\
            \Phi_n &= ((W_{n,1}, B_{n,1}), (W_{n,2}, B_{n,2}), \dots, (W_{n,L}, B_{n,L})) \in \textstyle \pr[\big]{\bigtimes_{k = 1}^{L}  ( \R^{ l_{n,k} \times l_{n, k- 1} }  \times \R^{l_{n,k}} )}.
        \end{split}
    \end{align}
    \cfload
    \noindent Then
    \begin{enumerate}[label=(\roman *)] 
        \item \label{item:parallelizations_of_anns_same_depth_item_element_in_set_of_anns} it holds that
        \begin{equation}
            \paraANN{n} ( \Phi_1, \dots, \Phi_n ) \in \textstyle \pr[\big]{\bigtimes_{j = 1}^{L} ( \R^{  (\sum_{k = 1}^{n} l_{k, j} )  \times  ( \sum_{k = 1}^{n} l_{k, j - 1} )}  \times \R^{ (\sum_{k = 1}^{n} l_{k, j} ) } )} \subseteq \bN,
        \end{equation}

        \item \label{item:parallelizations_of_anns_same_depth_item_realization} it holds that
        \begin{equation}
            \cR_\activation  ( \paraANN{n} ( \Phi_1, \dots, \Phi_n ) ) \in C ( \R^{ (\sum_{k = 1}^{n} \cI (\Phi_k) ) }, \R^{ ( \sum_{k = 1}^{n} \cO (\Phi_k) ) } ),
        \end{equation}

        \item \label{item:parallelizations_of_anns_same_depth_item_realization_evaluated} it holds for all $ x_1 \in \R^{\cI (\Phi_1)} $, $x_2 \in \R^{ \cI (\Phi_2) }$, \dots, $ x_n \in \R^{ \cI ( \Phi_n ) } $ that
        \begin{multline}
            \big(\cR_\activation ( \paraANN{n} ( \Phi_1, \dots, \Phi_n ) )\big)( x_1, x_2, \dots, x_n )\\
            = \left( \big(\cR_\activation (\Phi_1)\big)(x_1), \big(\cR_\activation (\Phi_2)\big)(x_2), \dots, \big(\cR_\activation (\Phi_n)\big)(x_n) \right) \in \R^{ ( \sum_{k = 1}^{n} \cO (\Phi_k) ) },
        \end{multline}
        and

        \item \label{item:parallelizations_of_anns_same_depth_item_layer_structure} it holds that
        \begin{equation}
            \cD  ( \paraANN{n} ( \Phi_1, \dots, \Phi_n ) ) ) = \bigg( \SmallSum{k = 1}{n} l_{k, 0}, \SmallSum{k = 1}{n} l_{k, 1}, \dots, \SmallSum{k = 1}{n} l_{k, L} \bigg) 
        \end{equation}

    \end{enumerate}
    \cfout.
\end{lemma}

\cfclear
\begin{definition}[Parallelizations of ANNs with different layer structure] \label{def:parallelizations_of_anns_with_different_layer_structure}
    \cfconsiderloaded{def:parallelizations_of_anns_with_different_layer_structure}
    Let $n\in\N$, $ \Psi = (\Psi_1, \Psi_2, \dots, \Psi_n) \in \bN^n$ satisfy for all $i \in \N \cap [1, n]$ that $\cL (\Psi_i) = 2$, $\cI (\Psi_i) = \cO (\Psi_i)$. Then we denote by $ \paraLANN{n}{\Psi} \colon \{ (\Phi_1, \Phi_2, \dots, \Phi_n) \in \bN^n : \cI (\Psi_i) = \cO (\Phi_i), i \in \N \cap [1, n]\} \rightarrow \bN$ the function which satisfies for all $ \Phi = (\Phi_1, \dots, \Phi_n) \in \bN^n$ with $ \cI (\Psi_i) = \cO (\Phi_i) $, $i \in \N \cap [1,n]$, that
    \begin{equation} \label{eq:parallelizations_of_anns_with_different_layer_structure_definition}
        \paraLANN{n}{\Psi}(\Phi) = \paraANN{n} \Big( \longerANN{\max_{i \in \N \cap  [1,n]  } \cL(\Phi_i), \Psi_1}(\Phi_1), \dots, \longerANN{\max_{i \in \N \cap [1,n] } \cL(\Phi_i), \Psi_n}(\Phi_n) \Big)
    \end{equation}
    \cfout.
\end{definition}

\cfclear
\begin{corollary} \label{cor:parallelizations_of_anns_with_different_layer_structure_realization_involving_identities}
    Let $n \in \N$, $\activation \in C(\R, \R)$, $\mathbb{I} = (\mathbb{I}_1, \mathbb{I}_2, \dots, \mathbb{I}_n), \Phi = (\Phi_1, \Phi_2, \dots, \Phi_n) \in \bN^n$ satisfy for all $i \in \N \cap [1, n] $, $x \in \R^{\cO (\Phi_i)}$ that $ \cL (\mathbb{I}_i) = 2 $, $\cI (\mathbb{I}_i) = \cO (\mathbb{I}_i) = \cO (\Phi_i)$, and $\big(\cR_\activation (\mathbb{I}_i)\big)(x) = x$. Then
    \begin{enumerate}[label=(\roman *)]
        \item \label{item:parallelizations_of_anns_with_different_layer_structure_realization_involving_identities_realization} it holds that $ \cR_\activation \big( \paraLANN{n}{\mathbb{I}}(\Phi)  \big) \in C \big( \R^{ (\sum_{i = 1}^{n} \cI (\Phi_i) ) }, \R^{ (\sum_{i = 1}^{n} \cO (\Phi_i) ) } \big) $ and
        
        \item \label{item:parallelizations_of_anns_with_different_layer_structure_realization_involving_identities_realization_evaluated} it holds for all $x_1 \in \R^{ \cI (\Phi_1) }, \dots, x_n \in \R^{ \cI (\Phi_n) }$ that
        
        \begin{equation}
            \big(\cR_\activation(\paraLANN{n}{\mathbb{I}}(\Phi) )\big)(x_1, x_2, \dots, x_n) = \begin{pmatrix}
                \big(\cR_\activation(\Phi_1)\big)(x_1)\\
                \big(\cR_\activation(\Phi_2)\big)(x_2)\\
                \vdots\\
                \big(\cR_\activation(\Phi_n)\big)(x_n)
            \end{pmatrix} \in \R^{ (\sum_{i = 1}^{n} \cO (\Phi_i) )  },
        \end{equation}
        and
        \item \label{item:parallelizations_of_anns_with_different_layer_structure_realization_involving_identities_architecture} it holds that
        \begin{equation}
            \begin{split}
                \cD \pr*{ \paraLANN{n}{\I}\pr{\Phi} } &= \bigg( \SmallSum{k = 1}{n} \dimANNlevel_0\pr[\big]{ \longerANN{\max_{i \in \N \cap \br{1,n}} \cL(\Phi_i), \Psi_k}\pr{\Phi_k} },
                \SmallSum{k = 1}{n} \dimANNlevel_1\pr[\big]{ \longerANN{\max_{i \in \N \cap \br{1,n}} \cL(\Phi_i), \Psi_k}\pr{\Phi_k} }, \dots \\
                &\quad \dots, \SmallSum{k = 1}{n} \dimANNlevel_{\max_{i \in \N \cap \br{1,n}} \cL\pr{\Phi_i} } \pr[\big]{ \longerANN{\max_{i \in \N \cap \br{1,n}} \cL(\Phi_i), \Psi_k}\pr{\Phi_k} } \bigg)
            \end{split}
        \end{equation}
    \end{enumerate}
    \cfout.
\end{corollary}

\begin{mproof}{\cref{cor:parallelizations_of_anns_with_different_layer_structure_realization_involving_identities}}
    First, note that \cite[Corollary~2.23]{GrohsHornungJentzen2019published} proves
    \cref{item:parallelizations_of_anns_with_different_layer_structure_realization_involving_identities_realization,item:parallelizations_of_anns_with_different_layer_structure_realization_involving_identities_realization_evaluated}.
    Moreover, observe that the assumption that for all $i \in \N \cap \br{1, n}$ it holds that $\cL\pr{\I_i} = 2$,
    the assumption that for all $i \in \N \cap \br{1,n}$ it holds that $\cI\pr{\I_i} = \cO\pr{\I_i} = \cO\pr{\Phi_i}$,
    and \cref{item:certain_extensions_of_anns_item_extension} in \cref{lem:certain_extensions_of_anns}
    ensure for all $k \in \N \cap \br{1,n}$ it holds that $\cL \pr{ \longerANN{\max_{i \in \N \cap \br{1,n}} \cL\pr{\Phi_i}, \Psi_k  } \pr{\Phi_k} } = \max_{i \in \N \cap \br{1,n}} \cL\pr{\Phi_i}$.
    Combining this,
    \eqref{eq:ANN_nth_layer},
    \eqref{eq:parallelizations_of_anns_with_different_layer_structure_definition},
    and \cref{item:parallelizations_of_anns_same_depth_item_layer_structure} in \cref{lem:parallelizations_of_anns_same_depth}
    yields that
    \begin{equation}
        \begin{split}
            \cD \pr*{ \paraLANN{n}{\I}\pr{\Phi} } &= \bigg( \SmallSum{k = 1}{n} \dimANNlevel_0\pr[\big]{ \longerANN{\max_{i \in \N \cap \br{1,n}} \cL(\Phi_i), \Psi_k}\pr{\Phi_k} },
            \SmallSum{k = 1}{n} \dimANNlevel_1\pr[\big]{ \longerANN{\max_{i \in \N \cap \br{1,n}} \cL(\Phi_i), \Psi_k}\pr{\Phi_k} }, \dots \\
            &\quad \dots, \SmallSum{k = 1}{n} \dimANNlevel_{\max_{i \in \N \cap \br{1,n}} \cL\pr{\Phi_i} } \pr[\big]{ \longerANN{\max_{i \in \N \cap \br{1,n}} \cL(\Phi_i), \Psi_k}\pr{\Phi_k} } \bigg).
        \end{split}
    \end{equation}
    This proves \cref{item:parallelizations_of_anns_with_different_layer_structure_realization_involving_identities_architecture}.
    The proof of \cref{cor:parallelizations_of_anns_with_different_layer_structure_realization_involving_identities} is thus complete.
\end{mproof}

\subsection{Linear combinations of ANNs}
\label{subsec:ann_calculus_linear_combinations}

\cfclear
\begin{definition} \label{def:sum_transpose_anns}
    \cfconsiderloaded{def:sum_transpose_anns}
    Let $m, n\in \N$. Then we denote by $\fS_{m,n} \in \R^{m \times (m n)} \times \R^m$, $\fT_{m, n} \in \R^{(m n) \times m} \times \R^{m n} $ the ANNs given by
    \begin{align}
        \begin{split}
            \fS_{m,n} &= \Big( \begin{pmatrix}
                \id_m &\id_m & \hdots & \id_m
            \end{pmatrix}, 0 \Big)  \qquad \text{and} \qquad
            \fT_{m, n} = \begin{pmatrix}
                \begin{pmatrix}
                    \id_m \\ \id_m \\ \vdots \\ \id_m
                \end{pmatrix}, 0
            \end{pmatrix}
        \end{split}
    \end{align}
    \cfout.
\end{definition}

\cfclear
\begin{definition}[Sums of ANNs with the same length] \label{def:sum_of_anns_same_length}
    \cfconsiderloaded{def:sum_of_anns_same_length}
    Let $u \in \Z$, $v \in \Z \cap [u, \infty)$, $\Phi_u, \Phi_{u + 1}, \dots, \Phi_v \in \bN$ satisfy for all $i \in \Z \cap [u, v]$ that $\cL(\Phi_i) = \cL (\Phi_u)$, $\cI(\Phi_i) = \cI(\Phi_u)$, $\cO(\Phi_i) = \cO(\Phi_u)$. Then we denote by $\oSum_{i = u}^v \Phi_i$ the ANN given by
    \begin{equation}
        \OSum{i = u}{v} \Phi_i = \Big( \fS_{\cO(\Phi_u), v - u + 1} \ANNcomp \big[ \paraANN{v - u + 1} (\Phi_u, \Phi_{u + 1}, \dots, \Phi_v) \big] \ANNcomp \fT_{\cI(\Phi_u), v - u + 1} \Big) \in \bN
    \end{equation}
    \cfout.
\end{definition}

\cfclear
\begin{definition}[Sums of ANNs with different lengths] \label{def:sum_of_anns_diff_length}
    \cfconsiderloaded{def:sum_of_anns_diff_length}
    Let $u \in \Z$, $v \in \Z \cap [u, \infty)$, $\Phi_u, \Phi_{u + 1}, \dots, \Phi_{v}, \Psi \in \bN$ satisfy for all $i \in \Z \cap [u, v]$ that $\cI( \Phi_i) = \cI(\Phi_u)$, $\cO(\Phi_i) = \cI(\Psi) = \cO(\Psi)$, and $\cL(\Psi) = 2$. Then denote by $\bSum_{i = u, \Psi}^v \Phi_i$ the ANN given by
    \begin{equation}
        \BSum{i = u}{\Psi}{v} \Phi_i = \OSum{i = u}{v} \longerANN{\max_{j \in \Z \cap [u, v]} \cL (\Phi_j), \Psi} (\Phi_i)
    \end{equation}
    \cfout.
\end{definition}

\cfclear
\begin{definition}[Scalar multiplications of ANNs]\label{def:scalar_multiplications_of_anns}
    \cfconsiderloaded{def:scalar_multiplications_of_anns}
    We denote by $\scalarANN \colon \R \times \bN \rightarrow \bN$ the function which satisfies for all $\alpha \in \R$, $\Phi \in \bN$ that
    \begin{equation}
        \alpha \circledast \Phi = \left( \alpha \id_{\cO (\Phi)}, 0 \right) \ANNcomp \Phi
    \end{equation}
    \cfout.
\end{definition}

\cfclear
\begin{lemma}\label{lem:scalar_multiplications_of_anns} 
    Let $\alpha \in \R$, $\Phi \in \bN$. Then
    \begin{enumerate}[label=(\roman *)]
        \item \label{item:scalar_multiplications_of_anns_layer_structure} it holds that $\cD (\alpha \scalarANN \Phi) = \cD (\Phi)$,
        
        \item \label{item:scalar_multiplications_of_anns_realization} it holds that $\cR_\activation ( \alpha \circledast \Phi ) \in C (\R^{\cI (\Phi)}, \R^{\cO (\Phi)})$, and
        
        \item \label{item:scalar_multiplications_of_anns_realization_evaluated} it holds for all $x \in \R^{\cI (\Phi)}$ that $ \big(\cR_\activation (\alpha \circledast \Phi)\big)(x) = \alpha \big( \cR_\activation (\Phi) \big)(x) $
    \end{enumerate}
    \cfout.
\end{lemma}

\cfclear
\begin{lemma}[Linear combinations of ANNs with the same length]
    \label{lem:linear_combinations_anns_same_length}
    %
    %
    Let $u \in \Z$, $v \in \Z \cap [u, \infty)$, $h_u, h_{u + 1}, \dots, h_v \in \R$, $\activation \in C(\R, \R)$, $\Phi_u, \Phi_{u + 1}, \dots, \Phi_v, \Psi \in \bN$, satisfy for all $i \in \Z \cap [u, v]$ that $\cI(\Phi_i) = \cI (\Phi_u)$, $\cO(\Phi_i) = \cO (\Phi_u)$, $\cL(\Phi_i) = \cL(\Phi_u)$ and $ \Psi = \oSum_{i = u}^{v} ( h_i \scalarANN \Phi_i)$. Then
    \begin{enumerate}[label=(\roman *)]
        \item \label{item:linear_combinations_anns_same_length_layer_structure} it holds that
        \begin{equation} 
            \cD (\Psi) = \Big( \cI (\Phi_u), \SmallSum{i = u}{v}\dimANNlevel_1(\Phi_i) ,\dots, \SmallSum{i = u}{v} \dimANNlevel_{\cL(\Phi_u) - 1} (\Phi_i) , \cO(\Phi_u) \Big),
        \end{equation}
        
        \item \label{item:linear_combinations_anns_same_length_realization} it holds that $\cR_\activation (\Psi) \in C(\R^{\cI (\Phi_u)}, \R^{\cO(\Phi_u)})$, and
        
        \item \label{item:linear_combinations_anns_same_length_realization_evaluated}
        it holds for all $x \in \R^{\cI(\Phi_u)}$ that
        \begin{equation}
            \big( \cR_\activation(\Psi) \big)(x) = \sum_{i = u}^{v} h_i \big( \cR_\activation(\Phi_i) \big)(x)
        \end{equation}
    \end{enumerate}
    \cfout[.]
\end{lemma}

\begin{mproof}{\cref{lem:linear_combinations_anns_same_length}}
    First, note that \cref{lem:scalar_multiplications_of_anns} ensures for all $i \in \Z \cap [u, v]$ that $\cD( h_i \circledast \Phi_i ) = \cD (\Phi_i)$. This, the fact that $\cD (\fS_{\cO(\Phi_u), v - u + 1}) = ((v - u + 1)\cO (\Phi_u), \cO(\Phi_u))$, the fact that $\cD (\fT_{\cI(\Phi_u), v - u + 1}) = ( \cI(\Phi_u), (v - u + 1)\cI(\Phi_u) )$, \cref{lem:compositions_of_anns}, and \cref{item:parallelizations_of_anns_same_depth_item_layer_structure} in \cref{lem:parallelizations_of_anns_same_depth} yield 
    \begin{equation}
        \cD (\Psi) = \Big( \cI(\Phi_u), \SmallSum{i = u}{v} \dimANNlevel_1(\Phi_i), \SmallSum{i = u}{v} \dimANNlevel_{2}(\Phi_i), \dots, \SmallSum{i = u}{v} \dimANNlevel_{\cL (\Phi_u) - 1} (\Phi_i) ,\cO(\Phi_u) \Big).
    \end{equation}
    This proves \cref{item:linear_combinations_anns_same_length_layer_structure}. Moreover, observe that \cref{lem:compositions_of_anns} and \cref{item:scalar_multiplications_of_anns_realization_evaluated} in \cref{lem:scalar_multiplications_of_anns} establish for all $x \in \R^{\cI(\Phi_u)}$ that
    \begin{align}
        \begin{split}
            \big( \cR_\activation (\Psi) \big)(x) &= \big( \cR_\activation(\fS_{\cO(\Phi_u), v - u + 1}) \circ \paraANN{v - u + 1}( h_u \circledast \Phi_u, h_{u + 1} \circledast \Phi_{u + 1}, \dots, h_v \circledast \Phi_v) \circ \fT_{\cI(\Phi_u), v - u + 1}\big) (x) \\
            &= \big( \cR_\activation (\fS_{\cO(\Phi_u), v-u +1}) \circ \paraANN{v - u + 1}(h_u \circledast \Phi_u, h_{u + 1} \circledast \Phi_{u + 1}, \dots, h_v \circledast \Phi_v) \big) (x, x, \dots, x)\\
            &= \big( \cR_\activation(\fS_{\cO(\Phi_u), v- u +1}) \big) \begin{pmatrix}
                \big( \cR_\activation (h_u \circledast \Phi_u) \big)(x) \\
                \big( \cR_\activation (h_{u + 1} \circledast \Phi_{u + 1}) \big)(x) \\
                \vdots\\
                \big( \cR_\activation (h_v \circledast \Phi_v) \big)(x)
            \end{pmatrix} = \sum_{i = u}^{v} h_i \big( \cR_\activation(\Phi_i) \big)(x).
        \end{split} 
    \end{align}
    This proves \cref{item:linear_combinations_anns_same_length_realization} and \cref{item:linear_combinations_anns_same_length_realization_evaluated}. The proof of \cref{lem:linear_combinations_anns_same_length} is thus complete.
\end{mproof}

\cfclear
\begin{lemma}[Linear combinations of ANNs with different lengths] \label{lem:linear_combination_anns_different_length}
    %
    %
    Let $L\in \N$, $\activation \in C(\R, \R)$, $u \in \Z$, $v \in [u, \infty) \cup \Z$, $h_u, h_{u + 1}, \dots, h_v \in \R$, $\Phi_u, \Phi_{u + 1}, \dots, \Phi_v, \Psi, \fJ \in \bN$ satisfy for all $x \in \R^{\cI(\fJ)}$, $i \in \Z \cap [u, v]$ that $\cI(\Phi_i) = \cI(\Phi_u)$, $\cO(\Phi_i) = \cI (\fJ) = \cO (\fJ)$, $L = \max_{j \in \Z \cap [u, v]} \cL(\Phi_j)$, $\cL (\fJ) = 2$, $\big( \cR_\activation (\fJ) \big)(x) = x$, and $\Psi = \bSum_{i = u, \fJ}^v (h_i \scalarANN \Phi_i)$. Then
    \begin{enumerate}[label=(\roman *)]
        \item \label{item:linear_combinations_anns_different_length_layer_structure} it holds that
        \begin{equation}
            \cD (\Psi) = \Big( \cI(\Phi_u), \SmallSum{i = u}{v} \dimANNlevel_1(\longerANN{L, \fJ}(\Phi_i)), \dots, \SmallSum{i = u}{v} \dimANNlevel_{L - 1}( \longerANN{L, \fJ}(\Phi_i) ), \cO(\Phi_u) \Big),
        \end{equation}

        \item \label{item:linear_combinations_anns_different_length_realization} it holds that $\cR_\activation (\Psi) \in C(\R^{\cI(\Phi_u)}, \R^{\cO(\Phi_u)}  )$, and
        
        \item \label{item:linear_combinations_anns_different_length_realization_evaluated} it holds for all $x \in \R^{\cI(\Phi_u)}$ that
        \begin{equation}
            \big( \cR_\activation(\Psi) \big)(x) = \sum_{i = u}^{v} h_i \big( \cR_\activation(\Phi_i) \big)(x)
        \end{equation}
    \end{enumerate} 
    \cfout.
\end{lemma}

\begin{mproof}{\cref{lem:linear_combination_anns_different_length}}
    First, note that \cref{item:scalar_multiplications_of_anns_layer_structure} in \cref{lem:scalar_multiplications_of_anns} ensures for all $i \in \Z \cap [u, v]$ that $\cD (h_i \scalarANN \Phi_i) = \cD(\Phi_i)$.\ %
    Combining this,
    the assumption that it holds that $\cL(\fJ) = 2$,
    the assumption that it holds that $\cI(\fJ) = \cO(\fJ)$,
    the assumption that for all $x \in \R^{\cI(\fJ)}$ it holds that $\pr*{ \cR_\activation(\fJ) }(x) = x$, and
    \cref{item:certain_extensions_of_anns_item_power,item:certain_extensions_of_anns_item_extension,item:certain_extensions_of_anns_item_extension_to_same_length} in \cref{lem:certain_extensions_of_anns}
    yields for all $i \in \Z \cap [u, v]$ that $\cD (\longerANN{L, \fJ} ( h_i \scalarANN \Phi_i ) ) = \cD(\longerANN{L, \fJ}(\Phi_i))$.\ %
    This and \cref{item:linear_combinations_anns_same_length_layer_structure} in \cref{lem:linear_combinations_anns_same_length}
    demonstrate that
    \begin{align}
        \begin{split}
            \cD(\Psi) &= \cD \pr*{ \BSum{i = u}{\fJ}{v} (h_i \scalarANN \Phi_i) } = \cD \pr*{ \OSum{i = u}{v} \longerANN{L, \fJ}(h_i \scalarANN \Phi_i)} \\
            &= \pr[\big]{ \cI(\longerANN{L, \fJ} (h_u \scalarANN \Phi_u)), \smallsum_{i = u}^{v} \dimANNlevel_1(\longerANN{ L, \fJ }(h_i \scalarANN \Phi_i)), \dots\\
            &\quad \dots, \smallsum_{i = u}^{v} \dimANNlevel_{L - 1}( \longerANN{L, \fJ}(h_i \scalarANN \Phi_i) ), \cO(\longerANN{L, \fJ} (h_u \scalarANN \Phi_u)) } \\
            &= \pr[\big]{ \cI(\Phi_u), \smallsum_{i = u}^{v} \dimANNlevel_1(h_i \scalarANN \Phi_i), \dots, \smallsum_{i = u}^{v} \dimANNlevel_{L - 1}(h_i \scalarANN \Phi_i), \cO (\Phi_u) } \in \N^{L + 2}.
        \end{split}
    \end{align}
    Hence, it holds that $\cR_\activation(\Psi) \in C(\R^{\cI(\Phi_u)}, \R^{\cO(\Phi_u)})$.\ %
    This proves \cref{item:linear_combinations_anns_different_length_layer_structure,item:linear_combinations_anns_different_length_realization}.\ %
    Furthermore, observe that \cref{item:linear_combinations_anns_same_length_realization_evaluated} in \cref{lem:linear_combinations_anns_same_length},
    \cref{item:scalar_multiplications_of_anns_realization_evaluated} in \cref{lem:scalar_multiplications_of_anns}
    assures for all $x \in \R^{\cI(\Phi_u)}$ that
    \begin{align}
        \begin{split}
            \pr*{ \cR_\activation (\Psi) }(x) &= \sum_{i = u}^{v} \pr[\big]{ \cR_\activation(h_i \scalarANN \Phi_i) }(x) = \sum_{i = u}^{v} h_i \pr[\big]{\cR_\activation(\Phi_i)}(x).
        \end{split}
    \end{align}
    This proves \cref{item:linear_combinations_anns_different_length_realization_evaluated}.\ %
    The proof of \cref{lem:linear_combination_anns_different_length} is thus complete.
\end{mproof}

\section{ANN representations}
\label{sec:ann_representations}

In this section, we develop representation results for ANNs.
In \cref{subsec:ann_representations_max} we establish a representation for maxima through ANNs with leaky ReLU activation function. To this end, we refer to the notions found in \cite[Section~4.2.2]{jentzen2023mathematicalintroductiondeeplearning} for a representation for maxima through ReLU ANNs.
In \cref{subsec:ann_representations_MLFP} we develop representations for MLFP approximations through ANNs with general activation function.

\subsection{ANN representations for the maximum function on finite sets}
\label{subsec:ann_representations_max}

\cfclear
\begin{lemma}\label{lem:max_one_step_leaky}
    Let $\beta \in [0, \infty)\backslash \{1 \}$,
    let $\activation \in C(\R, \R)$ satisfy for all $x \in \R$ that $\phi(x) = \max\{ x, \beta x\}$.
    Then there exists $\Phi_n \in \bN$, $n \in \N \cap [2, \infty)$, such that
    \begin{enumerate}[label=(\roman *)]
        \item \label{item:max_one_step_leaky_realization} it holds for all $n \in \N \cap [2, \infty)$ that \begin{equation}
            \cR_\activation (\Phi_n) \in \begin{cases}
                C(\R^n, \R^{\nicefrac{n}{2}}) &\colon n \text{ is even} \\
                C(\R^n, \R^{\nicefrac{(n + 1)}{2}}) &\colon n \text{ is odd},
            \end{cases}
        \end{equation}

        \item \label{item:max_one_step_leaky_realization_evaluated} it holds for all $n \in \N \cap [2, \infty)$, $x = (x_1, x_2, \dots, x_n) \in \R^n$ that
        \begin{equation}
            \big( \cR_\activation (\Phi_n) \big)(x) = \begin{cases}
                ( \max\{x_1, x_2\}, \dots, \max\{x_{n - 1}, x_n\}) \in \R^{\nicefrac{n}{2}} &\colon n \text{ is even}\\
                ( \max\{x_1, x_2 \}, \dots, \max\{x_{n - 2}, x_{n - 1} \}, x_n) \in \R^{\nicefrac{(n + 1)}{2}} &\colon n \text{ is odd,}
            \end{cases}
        \end{equation}

        \item \label{item:max_one_step_leaky_architecture} it holds for all $n \in \N \cap [2, \infty)$ that
        \begin{equation}
            \cD (\Phi_n) = \begin{cases}
                (n, 2n, \tfrac{n}{2}) &\colon n \text{ is even}\\
                (n, 2n, \tfrac{n + 1}{2}) &\colon n \text{ is odd},
            \end{cases}
        \end{equation}
        and

        \item \label{item:max_one_step_leaky_depth} it holds for all $n \in \N \cap [2, \infty)$ that $\cL(\Phi_n) = 2$
    \end{enumerate}
    \cfout.
\end{lemma}

\begin{mproof}{\cref{lem:max_one_step_leaky}}
    Throughout this proof let $\gamma, \delta \in \R$ satisfy $\gamma = \tfrac{|1 - \beta|}{(1 - \beta)(1 - \beta^2)}$ and $ \delta = \tfrac{|1 - \beta|}{1 - \beta}$,
    let $T \in \R^{1 \times 2}$, $w_1 \in \R^{4 \times 2}$, $w_2 \in \R^{1 \times 4}$ satisfy
    \begin{align}
        \begin{split} \label{eq:max_one_step_leaky_def_matrix}
            T &= \begin{pmatrix}
                (1 + \beta)^{-1} & - (1+ \beta)^{-1}
            \end{pmatrix}, \qquad w_1 = \begin{pmatrix}
                \phantom{-}\delta & -\delta \\
                -\delta & \phantom{-}\delta \\
                \phantom{-}0 & \phantom{-}\delta \\
                \phantom{-}0 & -\delta
            \end{pmatrix}, \qquad \text{and}\\
            w_2 &= \begin{pmatrix}
                \gamma & \gamma \beta & \gamma(1 - \beta) & \gamma (\beta - 1)
            \end{pmatrix},
        \end{split}
    \end{align}
    let $\fJ, \bM \in \bN$ satisfy
    \begin{equation} \label{eq:max_one_step_leaky_def_nets_id_max}
        \fJ = \begin{pmatrix}
            \begin{pmatrix}
                \begin{pmatrix}
                    \phantom{-}1 \\ -1
                \end{pmatrix}, \begin{pmatrix}
                    0 \\ 0
                \end{pmatrix}
            \end{pmatrix},
            \begin{pmatrix}
                T, 0
            \end{pmatrix}

        \end{pmatrix}
        \qquad \text{and} \qquad 
            \bM = \begin{pmatrix}
                \begin{pmatrix}
                    w_1, 0
                \end{pmatrix}, \begin{pmatrix}
                    w_2, 0
                \end{pmatrix}
            \end{pmatrix},
    \end{equation}
    let $\Phi_n \in \bN$, $n \in \N \cap [2, \infty)$, satisfy for all $n \in \N \cap [2, \infty)$ that
    \begin{align} \label{eq:max_one_step_leaky_def_Phi_n}
        \Phi_n = \begin{cases}
            \paraANN{\tfrac{n}{2}}(\bM, \bM, \dots, \bM) &: n \text{ is even},\\
            \paraANN{\frac{n + 1}{2}}(\bM, \bM, \dots, \bM, \fJ) &: n \text{ is odd}.
        \end{cases}
    \end{align}
    First, note that \cite[Equation (4.99)]{ackermann2023deep} ensures for all  $x \in \R$ that $\max\{x, 0\} = \gamma \beta \activation (- \delta x) + \gamma \activation (\delta x)$
    This,
    \eqref{eq:max_one_step_leaky_def_matrix},
    and \eqref{eq:max_one_step_leaky_def_nets_id_max},
    imply that for all $x \in \R$, $y = (y_1, y_2) \in \R^2$ it holds that $\cL(\fJ) = \cL(\bM) =  2$, $\cD(\fJ) = (1,2,1) \in \N^3$, $\cD(\bM) = (2, 4, 1) \in \N^3$, $\cR_\activation (\fJ) \in C(\R, \R)$, $\cR_\activation (\bM) \in C(\R^2, \R)$,
    \begin{align}
        \begin{split}
            \big( \cR_\activation (\fJ)\big)(x) &= T \activationDim{2} \left( \begin{pmatrix}
                \phantom{-}1 \\-1
            \end{pmatrix}x \right) = (1 + \beta)^{-1} (\activation (x) - \activation ( - x)) \\
            &= (1 + \beta)^{-1} \big( \max\{x, \beta x\} - \max\{- x, -\beta x\} \big) = \tfrac{x + \beta x}{1 + \beta} = x, \qquad\text{ and}    
        \end{split}\\
        \begin{split}
            \big( \cR_\activation (\bM) \big)(y) &= w_2 \activationDim{4} \left( w_1 \begin{pmatrix}
                y_1 \\ y_2
            \end{pmatrix} \right) = w_2 \activationDim{4} \left( \begin{pmatrix}
                \phantom{-}\delta & -\delta\\
                -\delta & \phantom{-}\delta\\
                \phantom{-}0 & \phantom{-}\delta \\
                \phantom{-}0 & -\delta
            \end{pmatrix} \begin{pmatrix}
                y_1 \\ y_2
            \end{pmatrix}\right) \\
            &= \begin{pmatrix}
                \gamma & \gamma \beta & \gamma(1 - \beta) & \gamma (\beta - 1)
            \end{pmatrix} \begin{pmatrix}
                \activation(\delta (y_1 - y_2)) \\
                \activation(- \delta (y_1 - y_2)) \\
                \activation(\delta y_2) \\
                \activation(-\delta y_2)
            \end{pmatrix} \\
            &= \gamma \activation(\delta (y_1 - y_2)) + \gamma \beta \activation(- \delta (y_1 - y_2)) + \gamma \activation(\delta y_2) + \gamma \beta \activation(-\delta y_2)\\
            &\quad- \gamma \beta \activation(\delta y_2) - \gamma \activation(-\delta y_2) \\
            &= \max\{y_1 - y_2, 0 \} + \max\{y_2, 0 \} - \max\{-y_2, 0 \} =  \max\{y_1, y_2\}.
        \end{split}
    \end{align}
    Combining this,
    \eqref{eq:max_one_step_leaky_def_Phi_n},
    and \cref{lem:parallelizations_of_anns_same_depth} therefore demonstrates for all $n \in \N \cap [2, \infty)$, $x = (x_1, x_2, \dots, x_n) \in \R^n$ that $\cL (\Phi_n) = 2$, 
    \begin{align}
        \begin{split}
            \cD( \Phi_n ) &= \begin{cases}
                (n, 2n, \frac{n}{2}) &: n \text{ is even}\\
                (n, 2n, \frac{n + 1}{2}) &: n \text{ is odd},
            \end{cases} 
        \end{split}\\
        \begin{split}
            \cR_\activation (\Phi_n) &\in \begin{cases}
                C(\R^n, \R^{\nicefrac{n}{2}}) &: n \text{ is even}\\
                C(\R^n, \R^{\nicefrac{(n + 1)}{2}}) &: n \text{ is odd},
            \end{cases}
        \end{split}
    \end{align}
    and 
    \begin{align}
        \begin{split}
            \big( \cR_\activation (\Phi_n) \big)(x) &= \begin{cases}
                \big(\cR_\activation (\paraANN{\tfrac{n}{2}}(\bM, \bM, \dots, \bM))\big)(x_1, x_2, \dots, x_n) &: n \text{ is even}\\
                \big(\cR_\activation(\paraANN{\frac{n + 1}{2}}(\bM, \bM, \dots, \bM, \fJ))\big) (x_1, x_2, \dots, x_n) &: n \text{ is odd}
            \end{cases}\\
            &= \begin{cases}
                \begin{pmatrix}
                    \big( \cR_\activation(\bM) \big)(x_1, x_2)\\
                    \big( \cR_\activation(\bM) \big)(x_3, x_4)\\
                    \vdots\\
                    \big( \cR_\activation(\bM) \big)(x_{n - 1}, x_n)\\
                \end{pmatrix} &: n \text{ is even}\\
                \begin{pmatrix}
                    \big( \cR_\activation(\bM) \big)(x_1, x_2)\\
                    \big( \cR_\activation(\bM) \big)(x_3, x_4)\\
                    \vdots\\
                    \big( \cR_\activation(\bM) \big)(x_{n - 2}, x_{n - 1})\\
                    \big( \cR_\activation(\fJ) \big)(x_n)
                \end{pmatrix} &: n \text{ is odd}
            \end{cases} \\
            &= \begin{cases}
                \big( \max\{x_1, x_2\}, \max\{x_3, x_4\}, \dots, \max\{x_{n - 1}, x_n\} \big) &: n \text{ is even}\\
                \big(\max\{x_1, x_2\}, \max\{x_3, x_4\}, \dots, \max\{x_{n - 2}, x_{n - 1}\}, x_n \big) &: n \text{ is odd}.
            \end{cases}
        \end{split}
    \end{align}
    This proves \cref{item:max_one_step_leaky_realization,item:max_one_step_leaky_realization_evaluated,item:max_one_step_leaky_architecture,item:max_one_step_leaky_depth}. The proof of \cref{lem:max_one_step_leaky} is thus complete.
\end{mproof}

\cfclear
\begin{lemma}[Maximum of finite set of real numbers (leaky ReLU ANN)] \label{lem:max_of_finite_set_leaky}
    Let $\beta \in [0, \infty) \backslash \{1\}$,
    for every $\dm \in \N$, $x = \pr{x_1, \dots, x_\dm} \in \R^\dm$ let $\fnorm{x} \in \R$ satisfy $\fnorm{x} = \max_{i \in \cu{1,\dots, \dm}} \abs{x_i}$,
    let $\ceil{\cdot} \colon \R \to \Z$ satisfy for all $x \in \R$ that $\ceil{x} = \min \{ l \in \Z : l \ge x \}$,
    let $\activation \in C(\R, \R)$ satisfy for all $x \in \R$ that $\activation(x) = \max\{x, \beta x\}$.
    Then there exists $\Psi_m \in \bN$, $m \in \N \cap [2, \infty)$, such that
    \begin{enumerate}[label=(\roman *)]
        \item \label{item:max_of_finite_set_leaky_realization} it holds for all $m \in \N \cap [2, \infty)$ that $\cR_\activation (\Psi_m) \in C(\R^m, \R)$,
        \item \label{item:max_of_finite_set_leaky_realization_evaluated} it holds for all $m \in \N \cap [2, \infty)$, $x = (x_1, x_2, \dots, x_m) \in \R^m$ that
        \begin{equation}
            \big( \cR_\activation (\Psi_m) \big)(x) = \max_{i \in [1, m] \cap \N} x_i,
        \end{equation}
        and
        \item \label{item:max_of_finite_set_leaky_width_depth_size} it holds for all $m \in \N \cap [2, \infty)$ that
        \begin{equation}
            \fnorm{\cD(\Psi_m)} \le 2m \qquad \text{ and } \qquad  \cL(\Psi_m) = \ceil{\log_2(m)} + 1
        \end{equation}
    \end{enumerate}
    \cfout.
\end{lemma}

\begin{mproof}{\cref{lem:max_of_finite_set_leaky}}
    First, note that \cref{lem:max_one_step_leaky} assures that there exist $\Phi_n \in \bN$, $n \in \N \cap [2, \infty)$, such that \cref{item:max_one_step_leaky_realization,item:max_one_step_leaky_realization_evaluated,item:max_one_step_leaky_architecture,item:max_one_step_leaky_depth} in \cref{lem:max_one_step_leaky} hold true.
    Throughout the remainder of this proof let $\Psi_m \in \bN$, $m \in \N \in [2, \infty)$, satisfy for all $m \in \N \cap [3, \infty)$ that $\Psi_2 = \Phi_2$ and $\Psi_m = \Psi_{\ceil{\frac{m}{2}}} \bullet \Phi_m$.
    We prove \cref{item:max_of_finite_set_leaky_realization,item:max_of_finite_set_leaky_realization_evaluated,item:max_of_finite_set_leaky_width_depth_size} by induction on $m \in \N \cap [2, \infty)$.
    For the base case $m = 2$ observe that the fact that $\Psi_2 = \Phi_2$
    and \cref{item:max_one_step_leaky_realization,item:max_one_step_leaky_realization_evaluated,item:max_one_step_leaky_architecture} in \cref{lem:max_one_step_leaky}
    imply for all $x = (x_1, x_2) \in \R^2$ that $\cR_\activation (\Psi_2) \in C(\R^2, \R)$, $\big( \cR_\activation (\Psi_2) \big)(x_1, x_2) = \max\{x_1, x_2\}$, $\fnorm{\cD(\Psi_2)} = 4$, and $\cL(\Psi_2) = 2 = \ceil{\log_2(2)} + 1$.
    This proves \cref{item:max_of_finite_set_leaky_realization,item:max_of_finite_set_leaky_realization_evaluated,item:max_of_finite_set_leaky_width_depth_size} for the base case $m = 2$.
    For the induction step $ \N \cap [2, \infty) \ni (m - 1) \induct m \in \N \cap [3, \infty)$ let $m \in \N \cap [3, \infty)$ and assume that \cref{item:max_of_finite_set_leaky_realization,item:max_of_finite_set_leaky_realization_evaluated,item:max_of_finite_set_leaky_width_depth_size} hold true for all $l \in \N \cap [2, m - 1]$.
    Observe that the fact that for all $l \in \Z$ it holds that
    \begin{equation}
        \ceil[\big]{\tfrac{l}{2}} = \begin{cases}
            \tfrac{l}{2} &: l \text{ is even}\\
            \tfrac{l + 1}{2} &: l \text{ is odd}
        \end{cases}
    \end{equation}
    and \cref{item:max_one_step_leaky_architecture} in \cref{lem:max_one_step_leaky} demonstrate for all $l \in \N \cap [2, \infty)$ that
    \begin{equation} \label{eq:max_of_finite_set_leaky_architecture_Phi_n}
        \cD (\Phi_l) = \left. \begin{cases}
            (l, 2l, \nicefrac{l}{2}) &: l \text{ is even}\\
            (l, 2l, \nicefrac{(l + 1)}{2}) &: l \text{ is odd}
        \end{cases} \right\} = (l, 2l, \ceil{\tfrac{l}{2}}) \in \N^3.
    \end{equation}
    This, \cref{lem:compositions_of_anns}, and the induction hypothesis establish that 
    \begin{align}
        \begin{split} \label{eq:max_of_finite_set_leaky_in_out_realization}
            \cI (\Psi_m) &= \cI (\Phi_m) = m, \qquad \cO (\Psi_m) = \cO( \Psi_{\ceil{\frac{m}{2}}}) = 1,\\
            \text{and} \qquad \cR_\activation(\Psi_m) &= \cR_\activation(\Psi_{\ceil{\frac{m}{2}}}) \circ \cR_\activation(\Phi_m) \in C(\R^m, \R).
        \end{split}
    \end{align}
    Moreover, note that \cref{lem:compositions_of_anns},
    \cref{item:max_one_step_leaky_realization_evaluated} in \cref{lem:max_one_step_leaky},
    the fact that for all $l \in \N \cap [2, \infty)$ it holds that $\ceil{\tfrac{l}{2}} < l$,
    and the induction hypothesis that for all $l \in \N \cap [2, m - 1]$, $x = (x_1, \dots, x_l) \in \R^l$ it holds that $\big( \cR_\activation (\Psi_l)\big)(x) = \max_{i \in \N \cap [1, l]} x_i$
    show for all $x = (x_1, x_2, \dots, x_m) \in \R^m$ that 
    \begin{align} \label{eq:max_of_finite_set_leaky_realization_evaluated}
        \big( \cR_\activation (\Psi_m) \big)(x) &= \big( \cR_\activation (\Psi_{\ceil{\frac{m}{2}}}) \circ \cR_\activation (\Phi_{m}) \big)(x_1, x_2, \dots, x_m) \nonumber \\
        &= \begin{cases}
            \max\!\big\{\! \max\{x_1, x_2\}, \max\{x_3, x_4\}, \dots, \max\{x_{m - 1}, x_m\}\! \big\} &: m \text{ is even}\\
            \max\!\big\{\! \max\{x_1, x_2\}, \max\{x_3, x_4\}, \dots, \max\{x_{m - 2}, x_{m - 1}\}, x_m\! \big\} &: m \text{ is even}
        \end{cases} \nonumber \\
        &= \max\{ x_1, x_2, \dots, x_m \}.
    \end{align}
    In addition, observe that \cref{lem:compositions_of_anns},
    \eqref{eq:max_one_step_leaky_def_Phi_n},
    the fact that for all $l \in \N \cap [2, \infty)$ it holds that $\ceil{\tfrac{l}{2}} < l$,
    and the induction hypothesis assure that
    \begin{equation} \label{eq:max_of_finite_set_leaky_width_Psi_m}
        \fnorm{ \cD (\Psi_m )} \le \max\big\{  \fnorm{\cD (\Psi_{\ceil{\frac{m}{2}}})}, \fnorm{ \cD (\Phi_m)} \big\} \le \max\big\{ 2 \ceil{\tfrac{m}{2}}, 2m \big\} = 2m.
    \end{equation}
    Furthermore, note that the fact that for all $l\in \N \cap [2, \infty)$ it holds that $\ceil{\log_2(2 \ceil{\nicefrac{l}{2}})} = \ceil{\log_2(l)}$,
    \cref{lem:compositions_of_anns},
    \cref{item:max_one_step_leaky_depth} in \cref{lem:max_one_step_leaky},
    and the induction hypothesis imply that
    \begin{align}
        \begin{split}
            \cL(\Psi_m) &= \cL(\Psi_{\ceil{\frac{m}{2}}}) + \cL(\Phi_m) - 1 = \cL(\Psi_{\ceil{\frac{m}{2}}}) + 1 = \big( \ceil{\log_2(\ceil{\nicefrac{m}{2}})} + 1 \big) + 1\\
            &= \big( \ceil{\log_2(2 \ceil{\nicefrac{m}{2}})} \big) + 1 = \ceil{\log_2(m)} + 1.
        \end{split}
    \end{align}
    Combining this,
    \eqref{eq:max_of_finite_set_leaky_in_out_realization},
    \eqref{eq:max_of_finite_set_leaky_realization_evaluated},
    and \eqref{eq:max_of_finite_set_leaky_width_Psi_m} completes the induction step. Induction hence proves \cref{item:max_of_finite_set_leaky_realization,item:max_of_finite_set_leaky_realization_evaluated,item:max_of_finite_set_leaky_width_depth_size}. The proof of \cref{lem:max_of_finite_set_leaky} is thus complete.
\end{mproof}

\subsection{ANN representations for MLFP approximations}
\label{subsec:ann_representations_MLFP}

\cfclear
\begin{lemma}[ANN representations for MLFP Approximations] \label{lem:mlfp_as_ann_general}
    Let $\dm, \fd \in \N$,
    let $M \in \N \cap [2, \infty)$,
    let $\Theta = \cup_{n \in \N} \Z^n$,
    let $A$ be a nonempty finite set,
    let $\action \in A$,
    let $\activation \in C(\R, \R)$,
    let $\fJ, \bF \in \bN$ satisfy $\cD (\fJ) = (1, \fd, 1)$, $\cR_\activation(\fJ) = \idfunc_{\R}$, and $\cR_\activation (\bF) \in C(\R^{\dm + |A|}, \R) $,
    let $\bX_a^\theta \in \bN$, $\theta \in \Theta$, $a \in A$, satisfy for all $\theta \in \Theta$, $a \in A$ that $\cD (\bX_a^\theta) = \cD (\bX_\action^0)$ and $\cR_\activation(\bX_a^\theta) \in C(\R^\dm, \R^\dm)$,
    for every $\dm \in \N$, $x \in \pr{x_1, \dots, x_\dm} \in \R^\dm$ let $\fnorm{x} \in \R$ satisfy $\fnorm{x} = \max_{i \in \cu{1,\dots, \dm}} \abs{x_i}$,
    let $c = 2 \max\{ \fnorm{\cD(\bF)}, \fnorm{\cD(\bX_\action^0)},  \fd \dm + \fd \abs{A} \}$,
    let $\varmlfp_n^\theta \colon \R^\dm \to \R^{|A|}$, $\theta \in \Theta$, $n \in \N_0$, satisfy for all $\theta \in \Theta$, $n \in \N_0$, $x \in \R^\dm$, $a \in A$ that
    \begin{multline} \label{eq:mlfp_as_ann_general_def_mlfp_scheme}
        \big(\varmlfp_{n}^{\theta} (x) \big)(a) = \sum_{l = 0}^{n - 1} \frac{1}{M^{n - l}} \sum_{i = 1}^{M^{n - l}} \Big[ \big( \cR_\activation(\bF) \big) \big(  \big( \cR_\activation(\bX_a^{(\theta, l, i)}) \big)(x), \varmlfp_{l}^{(\theta, l, i)} \big( \big(\cR_\activation( \bX_a^{(\theta, l, i)} ) \big)(x)  \big) \big)\\
        -\mathbbm{1}_{\N}\pr{l} \bigl( \cR_\activation (\bF) \bigr) \bigl( \bigl(\cR_\activation (\bX_a^{(\theta, l, i)})\bigr)(x) , \varmlfp_{\max\{l - 1, 0\}}^{(\theta, -l, i)} \bigl( \bigl( \cR_\activation (\bX_a^{(\theta, l, i)}) \bigr)(x) \bigr)  \bigr) \Big].
    \end{multline}
    Then there exists $\Phi_n^\theta \in \bN$, $\theta \in \Theta$, $n \in \N_0$, such that
    \begin{enumerate}[label=(\roman *)]
        \item \label{item:mlfp_as_ann_general_realization} it holds for all $\theta \in \Theta$, $n \in \N_0$ that $\cR_\activation(\Phi_n^\theta) \in C(\R^\dm, \R^{|A|})$,

        \item \label{item:mlfp_as_ann_general_realization_evaluated} it holds for all $\theta \in \Theta$, $n \in \N_0$, $x \in \R^\dm$ that $\bigl( \cR_\activation (\Phi_n^\theta) \bigr) (x) = \varmlfp_n^\theta(x) \in \R^{|A|}$,

        \item \label{item:mlfp_as_ann_general_layer_structure} it holds for all $\theta \in \Theta$, $n \in \N_0$ that $ \cD(\Phi_n^\theta) = \cD(\Phi_n^0)$,
        
        \item \label{item:mlfp_as_ann_general_depth} it holds for all $\theta \in \Theta$, $n \in \N_0$ that
        \begin{equation}
            \cL(\Phi_n^\theta) \le n \big( \cL(\bF) + \cL(\bX_\action^0) - 1 \big) + 1,
        \end{equation}
        
        \item \label{item:mlfp_as_ann_general_width} it holds for all $\theta \in \Theta$, $n \in \N_0$ that
        \begin{equation}
            \fnorm{\cD (\Phi_n^\theta)} \le c (4 |A| M)^n,
        \end{equation}
        and
        \item \label{item:mlfp_as_ann_general_size} it holds for all $\theta \in \Theta$, $n \in \N_0$ that
        \begin{equation}
            \paramANN(\Phi_n^\theta) \le 2 \bigl( n \bigl( \cL(\bF) + \cL(\bX_\action^0) - 1 \bigr) + 1 \bigr) c^2 (4|A|M)^{2n}
        \end{equation}
    \end{enumerate}
    \cfout.
\end{lemma}

\begin{mproof}{\cref{lem:mlfp_as_ann_general}}
    Throughout this proof let $\I_k \in \bN$, $k \in \N$, satisfy for all $k \in \N$ that
    \begin{equation}
        \I_k = \paraANN{k} (\fJ, \dots, \fJ),    
    \end{equation}
    let $\mathbb{J} \in \bN^2$ satisfy $\mathbb{J} = (\I_\dm, \I_{\abs{A}})$.
    Next, observe that \cref{lem:parallelizations_of_anns_same_depth} ensures for all $k \in \N$, $x \in \R^k$ that $\cD (\I_k) = (k, k \fd, k) \in \N^3$, $\cL(\I_k) = 2$, $\cR_\activation (\I_k) \in C(\R^k, \R^k)$, and $\bigl( \cR_\activation (\I_k) \bigr)(x) = x$.
    Throughout the remainder of this proof let $\Lambda_n^{\theta, a }, \Phi_n^\theta, \Xi_n^\theta \in \bN$, $\theta \in \Theta$, $n \in \N_0$, $a \in A$, satisfy for all $\theta \in \Theta$, $n \in \N$, $a \in A$ that
    \begin{align}
        \begin{split} \label{eq:mlfp_as_ann_general_def_base_case_ANN} 
            \Lambda_0^{\theta, a} &= \left( \begin{pmatrix}
                0 & 0 & \cdots & 0
            \end{pmatrix}, 0 \right) \in \R^{1 \times \dm} \times \R^1, \qquad \Phi_0^\theta = \paraANN{\abs{A}} ( (\Lambda_0^{\theta, a})_{a \in A} ) \bullet \fT_{\dm, \abs{A}},\\
            \Xi_0^\theta &= \left( 0, 0 \right) \in \R^{({\abs{A}} + \dm) \times \dm} \times \R^{\abs{A} + \dm},
        \end{split}\\
        \phantom{a} \nonumber \\
        \begin{split} \label{eq:mlfp_as_ann_general_Lambda}
            \Lambda_n^{\theta, a} &= \Bigg[ \BSum{l = 0}{\fJ}{n-1} \Bigg[ \scalar{\frac{1}{M^{n - l}}}{ \left( \BSum{i = 1}{\fJ}{M^{n - l}} \left( \bF \bullet \Xi_l^{(\theta, l, i)} \bullet \bX_a^{(\theta, l, i)} \right) \right)} \Bigg] \Bigg] \\
            &\quad \bSum_{\,\fJ} \Bigg[ \BSum{l = 0}{\fJ}{n - 1} \Bigg[ \scalar{\frac{- \mathbbm{1}_\N (l)}{ M^{n - l} }}{ \left( \BSum{i = 1}{\fJ}{M^{n - l}} \left( \bF \bullet \Xi_{\max\{ l - 1, 0 \} }^{(\theta, -l , i)} \bullet \bX_a^{(\theta, l, i)} \right) \right) } \Bigg]  \Bigg],
        \end{split}\\
        \phantom{a} \nonumber \\
        \begin{split} \label{eq:mlfp_as_ann_general_Phi}
            \Phi_n^\theta &= \paraANN{\abs{A}} ( ( \Lambda_n^{\theta, a} )_{a \in A} ) \bullet \fT_{\dm, \abs{A}},
        \end{split}\\
        \phantom{a} \nonumber \\
        \begin{split} \label{eq:mlfp_as_ann_general_Xi}
            \Xi_n^{\theta} &= \paraLANN{2}{\mathbb{J}} \left( \I_\dm, \Phi_n^\theta \right) \bullet \fT_{\dm, 2},
        \end{split} 
    \end{align}   
    let $\Psi_{n, l}^{\theta, j, a} \in \bN$, $\theta \in \Theta$, $n \in \N$, $l \in \N_0 \cap [0, n - 1]$, $j \in \{0,1\}$, $a \in A$,
    satisfy for all $\theta \in \Theta$, $n \in \N$, $l \in \N_0 \cap [0, n - 1]$, $j \in \{0,1\}$, $a \in A$ that
    \begin{equation} \label{eq:mlfp_as_ann_general_Psi}
        \Psi_{n, l}^{\theta, j, a} = \BSum{i = 1}{\fJ}{M^{n - l}} \left( \bF \bullet \Xi_{\max\{ l - j, 0 \} }^{( \theta, (-1)^jl, i )} \bullet \bX_a^{(\theta, l, i)} \right) \in \bN,
    \end{equation}
    and let $\Gamma_n^{\theta, j, a} \in \bN$, $\theta \in \Theta$, $n \in \N$, $j \in \{0,1\}$, $a \in A$,
    satisfy for all $\theta \in \Theta$, $n \in \N$, $j \in \{0,1\}$, $a \in A$ that
    \begin{equation} \label{eq:mlfp_as_ann_general_Gamma}
        \Gamma_n^{\theta, j, a} = \BSum{l = 0}{\fJ}{n - 1} \left[ \scalar{\tfrac{ (-1)^j \mathbbm{1}_\N(l + 1 - j) }{ M^{n - l} }}{ \Psi_{n, l}^{\theta, j , a} } \right],
    \end{equation}
    let $L^j_{l} \in \N$, $l \in \N_0 \cap [0, \infty) $, $j \in \{0,1\}$, satisfy for all $l \in \N_0 \cap [0, \infty)$, $j \in \{0,1\}$ that
    \begin{equation}
        L^j_{l} = \cL \bigl( \bF \bullet \Xi_{\max \{l - j, 0\} }^{0} \bullet \bX_\action^0 \bigr),
    \end{equation}
    let $\fL_n^j \in \N$, $n \in \N$, $j \in \{0,1\}$, satisfy for all $n\in \N$, $j \in \{0,1\}$ that
    \begin{equation}
        \fL_n^j = \max_{l \in \N_0 \cap [0, n-1]} L_{l}^j,
    \end{equation}
    let $\mathbb{L}_n \in \N$, $n \in \N$, satisfy for all $n \in \N$ that
    \begin{equation} \label{eq:mlfp_as_ann_general_depth_nth_iterate} 
        \mathbb{L}_n = \max_{j \in \{0,1\}} \fL_n^j.
    \end{equation}
    Observe that for $\theta \in \Theta$, $n \in \N$, $a \in A$ it holds that
    \begin{equation}
        \Lambda_n^{\theta, a} = \Gamma_{n}^{\theta, 0, a} \,\bSum_{\,\fJ}\, \Gamma_n^{\theta, 1, a}.
    \end{equation}
    We prove \cref{item:mlfp_as_ann_general_realization,item:mlfp_as_ann_general_realization_evaluated,item:mlfp_as_ann_general_layer_structure,item:mlfp_as_ann_general_depth,item:mlfp_as_ann_general_width,item:mlfp_as_ann_general_size} by induction on $n \in \N_0$.
    For the base case $n = 0$ note that \eqref{eq:mlfp_as_ann_general_def_base_case_ANN} ensures for all $\theta \in \Theta$, $a \in A$ that $\Lambda_0^{\theta, a} = (0, 0) \in \R^{1 \times \dm} \times \R^1$.
    This establishes for all $\theta \in \Theta$, $a \in A$ that $\cD (\Lambda_0^{\theta, a}) = \cD (\Lambda_0^{0, \action}) = (\dm, 1)$.
    Combining this, the fact that for all $j, k \in \N$ it holds that $\cD (\fT_{k, j}) = (k, kj)$, \eqref{eq:mlfp_as_ann_general_def_base_case_ANN}, \cref{lem:compositions_of_anns},
    and \cref{lem:parallelizations_of_anns_same_depth}
    proves for all $\theta \in \Theta$ that $\cD (\Phi_0^\theta) = \cD (\Phi_0^0) = (\dm, \abs{A})$.
    Hence, \cref{item:mlfp_as_ann_general_layer_structure} holds for the base case $n = 0$.
    Moreover, note that \eqref{eq:mlfp_as_ann_general_def_base_case_ANN} ensures for all $\theta \in \Theta$ that $\cD(\Xi_0^\theta) = \cD(\Xi_0^0) = (\dm, \dm + \abs{A})$.
    Thus for all $\theta \in \Theta$, $a \in A$ it holds that
    \begin{equation}
        \cL (\Xi_0^\theta) = \cL(\Phi_0^\theta) = \cL (\Lambda_0^{\theta, a}) =  0 \cdot \bigl( \cL (\bF) + \cL (\bX_\action^0) - 1 \bigr) + 1.
    \end{equation}
    This proves \cref{item:mlfp_as_ann_general_depth} for the base case $n = 0$.
    The fact that for all $\theta \in \Theta$ it holds that $\cD(\Phi_0^\theta) = (\dm, \abs{A})$ implies for all $\theta \in \Theta$ that $\fnorm{ \cD(\Phi_0^\theta) } = \max \{\dm, \abs{A}\} \le \dm + \abs{A} \le c$ and $ \paramANN(\Phi_0^\theta) = \abs{A}(\dm + 1) \le c^2 \le 2c^2$.
    This establishes \cref{item:mlfp_as_ann_general_width,item:mlfp_as_ann_general_size} for the base case $n = 0$.
    Furthermore, note that \eqref{eq:mlfp_as_ann_general_def_base_case_ANN} implies for all $\theta \in \Theta$, $x \in \R^\dm$, $a \in A$ that $\cR_\activation (\Lambda_0^{\theta, a}) \in C(\R^\dm, \R)$ and $\bigl( \cR_\activation(\Lambda_0^{\theta, a}) \bigr)(x) = 0$.
    This,
    \eqref{eq:mlfp_as_ann_general_def_base_case_ANN},
    \cref{lem:compositions_of_anns},
    \cref{lem:parallelizations_of_anns_same_depth},
    and the fact that for all $\theta \in \Theta$, $x \in \R^\dm$, $a \in A$ it holds that $ \bigl( \varmlfp_0^{\theta}(x) \bigr)(a) = 0$
    demonstrate for all $\theta \in \Theta$, $x \in \R^\dm$ that $\cR_\activation (\Phi_0^\theta) \in C (\R^\dm, \R^{|A|})$ and
    \begin{align}
        \begin{split}
            \bigl(\cR_\activation (\Phi_0^\theta) \bigr)(x) &= \bigl(  \cR_\activation( \paraANN{\abs{A}}((\Lambda_0^{\theta, a})_{a \in A}) \bullet \fT_{\dm, \abs{A}} )  \bigr)(x) = \bigl( \cR_\activation ( \paraANN{\abs{A}}((\Lambda_0^{\theta, a})_{a \in A}) ) \bigr)\begin{pmatrix}
                x \\x \\ \vdots \\ x
            \end{pmatrix}\\
            &= \Bigl( \bigl( \cR_\activation (\Lambda_0^{\theta, a}) \bigr)(x) \Bigr)_{a \in A} = 0 = \varmlfp_0^\theta(x).
        \end{split}
    \end{align}
    Hence \cref{item:mlfp_as_ann_general_realization,item:mlfp_as_ann_general_realization_evaluated} hold true for the base case $n = 0$.
    Thus \cref{item:mlfp_as_ann_general_realization,item:mlfp_as_ann_general_realization_evaluated,item:mlfp_as_ann_general_layer_structure,item:mlfp_as_ann_general_depth,item:mlfp_as_ann_general_width,item:mlfp_as_ann_general_size} hold for the base case $n = 0$.
    For the induction step $\N_0 \ni (n - 1) \induct n \in \N$ let $n \in \N$ and assume that \cref{item:mlfp_as_ann_general_realization,item:mlfp_as_ann_general_realization_evaluated,item:mlfp_as_ann_general_layer_structure,item:mlfp_as_ann_general_depth,item:mlfp_as_ann_general_width,item:mlfp_as_ann_general_size} hold for all $l \in \N_0 \cap [0, n - 1]$.
    First, note that \eqref{eq:mlfp_as_ann_general_Xi},
    \cref{lem:compositions_of_anns},
    \cref{cor:parallelizations_of_anns_with_different_layer_structure_realization_involving_identities},
    and the induction hypothesis imply for all $\theta \in \Theta$, $l \in \N_0 \cap [0, n - 1]$, $x \in \R^\dm$ that
    \begin{align}
        \begin{split} \label{eq:mlfp_as_ann_general_Xi_realization_evaluated} 
            \bigl( \cR_\activation (\Xi_l^\theta)\bigr) (x) &= \bigl( \cR_\activation ( \paraLANN{2, \mathbb{J}} ( \I_\dm, \Phi_l^\theta ) \bullet \fT_{\dm, 2} ) \bigr)(x) = \bigl( \cR_\activation (\paraLANN{2, \mathbb{J}} ( \I_\dm, \Phi_l^\theta ) ) \circ \cR_\activation (\fT_{\dm, 2}) \bigr) (x)\\
            &= \bigl(\cR_\activation (\paraLANN{2, \mathbb{J}} ( \I_\dm, \Phi_l^\theta ) )\bigr)(x, x) = \begin{pmatrix}
                \bigl( \cR_\activation (\I_\dm) \bigr)(x) \\ \bigl( \cR_\activation (\Phi_l^\theta) \bigr)(x)
            \end{pmatrix} = \begin{pmatrix}
                x \\ \varmlfp_l^\theta(x)
            \end{pmatrix} \in \R^{\dm + \abs{A}}.
        \end{split}
    \end{align}
    This and \cref{lem:compositions_of_anns} demonstrate for all $\theta, \vartheta \in \Theta$, $l \in \N_0 \cap [0, n - 1]$, $x \in \R^\dm$, $a \in A$ that 
    \begin{align}
        \begin{split} \label{eq:mlfp_as_ann_general_summands_realization_evaluated}
            \bigl( \cR_\activation (\bF \bullet \Xi_l^\vartheta \bullet \bX_a^\theta) \bigr)(x) &= \bigl( \cR_\activation(\bF) \circ \cR_\activation(\Xi_l^\vartheta) \bigr) \bigl( \bigl( \cR_\activation(\bX_a^\theta) \bigr)(x) \bigr)\\
            &= \bigl( \cR_\activation (\bF) \bigr) \Bigl( \bigl( \cR_\activation(\bX_a^\theta) \bigr)(x) , \varmlfp_l^\vartheta \bigl( \bigl( \cR_\activation(\bX_a^\theta) \bigr)(x) \bigr)  \Bigr).
        \end{split}
    \end{align}
    Moreover, observe that \eqref{eq:mlfp_as_ann_general_Xi} and the induction hypothesis
    establish for all $\theta \in \Theta$, $l \in \N_0 \cap [0, n -1]$ that $\cD (\Xi_l^\theta) = \cD(\Xi_l^0) $.
    Combining this,
    \cref{lem:compositions_of_anns},
    and the assumption that for all $\theta \in \Theta$, $a \in A$ it holds that $\cD (\bX_a^\theta) = \cD (\bX_\action^0)$
    proves for all $\theta, \vartheta \in \Theta$, $l \in \N_0 \cap [0, n - 1]$, $j \in \{0, 1\}$, $a \in A$ that 
    \begin{equation} 
        \cD (\bF \bullet \Xi_{\max\{l - j, 0\}}^\vartheta \bullet \bX_a^\theta) = \cD ( \bF \bullet \Xi_{\max \{l - j, 0\}}^0 \bullet \bX_\action^0 ). 
    \end{equation}
    This implies for all $\theta, \vartheta \in \Theta$, $l \in \N_0 \cap [0, n - 1]$, $j \in \{0, 1\}$, $a \in A$ that 
    \begin{equation}
        \cL \bigl( \bF \bullet \Xi_{\max\{ l - j, 0 \}}^\vartheta \bullet \bX_a^\theta \bigr) = \cL \bigl(\bF \bullet \Xi_{\max \{l - j, 0\}}^0 \bullet \bX_\action^0 \bigr) = L_{l}^j.
    \end{equation}
    This,
    \cref{lem:linear_combination_anns_different_length} (applied for every $\theta \in \Theta$, $l \in \N_0 \cap [0, n - 1]$, $j \in \{0, 1\}$, $a \in A$ with
    \begin{align}
        \begin{split}
            &u \with 1, \qquad v \with M^{n - l}, \qquad \fJ \with \fJ, \qquad (h_k)_{k \in \Z \cap [u, v]} \with (1)_{i \in \N \cap [1, M^{n - l}]},\\
            &L \with L^j_{l}, \qquad \activation \with \activation, \qquad \Psi \with \Psi_{n, l}^{\theta, j, a},\\
            &(\Phi_k)_{k \in \Z \cap [u, v]} \with (\bF \bullet \Xi_{\max\{l - j, 0\}}^{(\theta, (-1)^jl, i)} \bullet \bX_a^{(\theta, l, i)})_{i \in \N \cap [1, M^{n - l}]}
        \end{split}
    \end{align}
    in the notation of \cref{lem:linear_combination_anns_different_length}),
    and \cref{item:certain_extensions_of_anns_item_extension_to_same_length} in \cref{lem:certain_extensions_of_anns} yield for all $\theta \in \Theta$, $l \in \N_0 \cap [0, n - 1]$, $j \in \{0, 1\}$, $x \in \R^\dm$, $a \in A$ that
    \begin{align}
        \begin{split} \label{eq:mlfp_as_ann_general_Psi_layer_structure}
            &\cD (\Psi_{n, l}^{\theta, j, a})\\
            &= \pr[\bigg]{\dm, \SmallSum{i = 1}{M^{n - l}} \dimANNlevel_1 (\longerANN{L^j_{l}, \fJ} (\bF \bullet \Xi_{\max\{l - j, 0\}}^{(\theta, (-1)^jl, i)} \bullet \bX_a^{(\theta, l, i)}) ), \SmallSum{i = 1}{M^{n - l}} \dimANNlevel_2  (\longerANN{L^j_{l}, \fJ} (\bF \bullet \Xi_{\max\{l - j, 0\}}^{(\theta, (-1)^jl, i)} \bullet \bX_a^{(\theta, l, i)}) ), \dots\\
            &\quad\dots,\SmallSum{i = 1}{M^{n - l}} \dimANNlevel_{L_{l}^j - 1} (\longerANN{L^j_{l}, \fJ} (\bF \bullet \Xi_{\max\{l - j, 0\}}^{(\theta, (-1)^jl, i)} \bullet \bX_a^{(\theta, l, i)}) ), 1} \\
            &= \pr[\bigg]{\dm, M^{n - l} \dimANNlevel_1 (\bF \bullet \Xi^0_{\max\{l - j, 0\}} \bullet \bX_\action^0 ),  M^{n - l} \dimANNlevel_2 (\bF \bullet \Xi^0_{\max\{l - j, 0\}} \bullet \bX_\action^0 ), \dots\\
            &\quad\dots,  M^{n - l} \dimANNlevel_{L^j_{l} - 1} (\bF \bullet \Xi^0_{\max\{l - j, 0\}} \bullet \bX_\action^0 ), 1} \in \N^{L^j_{l} + 1}
        \end{split}
    \end{align}
    and
    \begin{align}
        \begin{split} \label{eq:mlfp_as_ann_general_Psi_realization_evaluated}
            \bigl(\cR_\activation(\Psi_{n, l}^{\theta, j, a})\bigr)(x) &= \sum_{i = 1}^{M^{n - l}} \Bigl( \cR_\activation  \bigl( \bF \bullet \Xi_{\max\{l - j, 0\}}^{(\theta, (-1)^jl, i)} \bullet \bX_a^{(\theta, l, i)}  \bigr) \Bigr)(x)\\
            &= \sum_{i = 1}^{M^{n - l}} \Bigl( \cR_\activation (\bF) \Bigr) \Bigl( \bigl( \cR_\activation (\bX_a^{(\theta, l, i)}) \bigr)(x) , \varmlfp_{\max\{l - j, 0\}}^{(\theta, (-1)^jl, i)} \bigl(  \bigl( \cR_\activation (\bX_a^{(\theta, l, i)}) \bigr)(x)  \bigr)  \Bigr).
        \end{split}
    \end{align}
    Furthermore, observe that \eqref{eq:mlfp_as_ann_general_Psi_layer_structure} and \cref{lem:scalar_multiplications_of_anns} assure for all $\theta \in \Theta$, $l \in \N_0 \cap [0, n - 1]$, $j \in \{0, 1\}$, $a \in A$ that
    \begin{equation}
        \cD \Bigl( \scalar{\tfrac{(-1)^j\mathbbm{1}_{\N}(l + 1 - j)}{M^{n - l}}}{\Psi_{n, l}^{\theta, j, a}} \Bigr) = \cD (\Psi_{n, l}^{\theta, j, a}) = \cD (\Psi_{n, l}^{0, j, \action}).
    \end{equation}
    Hence, it holds for all $\theta \in \Theta$, $l \in \N_0 \cap [0, n - 1]$, $j \in \{0, 1\}$, $a \in A$ that $\cL (\Psi_{n, l}^{\Theta, j, a}) = \cL (\Psi_{n, l}^{0, j, \action}) = L^j_{l}$. This ensures for all $\theta \in \Theta$, $a \in A$, $j \in \{0, 1\}$ that $\max_{l \in \N_0 \cap [0, n - 1] } \cL(\Psi_{n, l}^{\theta, j, a}) = \max_{l \in \N_0 \cap [0, n - 1]} L_{l}^j = \fL_n^j $.
    Moreover, note that this,
    \eqref{eq:mlfp_as_ann_general_Gamma},
    and \cref{lem:linear_combination_anns_different_length} (applied for every $\theta \in \Theta$, $j \in \{0, 1\}$, $a \in A$ with
    \begin{align}
        \begin{split}
            &u \with 0, \qquad v \with n - 1, \qquad \fJ \with \fJ, \qquad (h_k)_{k \in \Z \cap [u, v]} \with \Bigl( \tfrac{(-1)^j \mathbbm{1}_{\N}(l + 1 - j)}{M^{n - l}} \Bigr)_{l \in \N_0 \cap [0, n - 1]}, \\
            &L \with \fL_n^j, \qquad \activation \with \activation, \qquad \Psi \with \Gamma_n^{\theta, j , a},\qquad (\Phi_k)_{k \in \Z \cap [u, v]} \with ( \Phi_{n, l}^{\theta, j, a} )_{l \in \N_0 \cap [0, n - 1]}
        \end{split}
    \end{align}
    in the notation of \cref{lem:linear_combination_anns_different_length}) establish for all $\theta \in \Theta$, $j \in \{0 , 1\}$, $x \in \R^\dm$, $a \in A$ that
    \begin{align}
        \begin{split} \label{eq:mlfp_as_ann_general_Gamma_layer_structure}
            \cD ( \Gamma_n^{\theta, j, a} ) &=  \pr*{\dm, \SmallSum{l = 0}{n - 1} \dimANNlevel_1( \longerANN{\fL_n^j, \fJ} (\Psi_{n, l}^{\theta, j, a})  ), \SmallSum{l = 0}{n - 1} \dimANNlevel_2( \longerANN{\fL_n^j, \fJ} (\Psi_{n, l}^{\theta, j, a}) ), \dots, \SmallSum{l = 0}{n - 1} \dimANNlevel_{\fL_n^j - 1}(\longerANN{\fL_n^j, \fJ} (\Psi_{n, l}^{\theta, l, a}) ), 1} \\
            &= \pr*{d, \SmallSum{l = 0}{n - 1} \dimANNlevel_1( \longerANN{\fL_n^j, \fJ} (\Psi_{n, l}^{0, j, \action}) ), \SmallSum{l = 0}{n - 1} \dimANNlevel_2( \longerANN{\fL_n^j, \fJ} (\Psi_{n, l}^{0, j, \action}) ), \dots, \SmallSum{l = 0}{n - 1} \dimANNlevel_{\fL_n^j - 1}( \longerANN{\fL_n^j, \fJ} (\Psi_{n, l}^{0, j, \action}) ), 1} \\
            &= \cD (\Gamma_n^{0, j, \action}) \in \N^{\fL_n^j + 1}
        \end{split}
    \end{align}
    and 
    \begin{align}
        \begin{split} \label{eq:mlfp_as_ann_general_Gamma_realization_evaluated}
            \pr*{\cR_\activation (\Gamma_n^{\theta, j, a})} (x) &= \sum_{l = 0}^{n - 1} \tfrac{(-1)^j\mathbbm{1}_{\N}(l + 1 - j)}{M^{n - l}} \pr[\big]{\cR_\activation (\Psi_{n, l}^{\theta, j, a})} (x)\\
            &= \sum_{l = 0}^{n - 1} \tfrac{(-1)^j\mathbbm{1}_{\N}(l + 1 - j)}{M^{n - l}} \sum_{i = 1}^{M^{n - l}}  \pr*{\cR_\activation (\bF) \circ \cR_\activation (\Xi_{\max\{l - j, 0\}}^{(\theta, (-1)^jl, i)}) \circ \cR_\activation (\bX_a^{(\theta, l, i)})} (x).
        \end{split}
    \end{align}
    Note that \eqref{eq:mlfp_as_ann_general_Gamma_layer_structure} implies for all $\theta \in \Theta$, $a \in A$ that
    \begin{equation}
        \max\{ \cL(\Gamma_n^{\theta, 0, a}), \cL(\Gamma_n^{\theta, 1, a}) \} = \max\{ \cL (\Gamma_n^{0, 0, \action}), \cL (\Gamma_n^{0, 1, \action}) \} = \max\{ \fL_n^0, \fL_n^1 \} = \mathbb{L}_n.
    \end{equation}
    Combining this, \eqref{eq:mlfp_as_ann_general_Gamma}, \eqref{eq:mlfp_as_ann_general_Xi_realization_evaluated}, \eqref{eq:mlfp_as_ann_general_summands_realization_evaluated}, \eqref{eq:mlfp_as_ann_general_Gamma_layer_structure}, \eqref{eq:mlfp_as_ann_general_Gamma_realization_evaluated}, and \cref{lem:linear_combination_anns_different_length} (applied for every $\theta \in \Theta$, $a \in A$ with $u \with 1$, $v \with 2$, $\fJ \with \fJ$, $h_1 \with 1$, $h_2 \with 1$, $L \with \mathbb{L}_n$, $\activation \with \activation$, $\Phi_1 \with \Gamma_n^{\theta, 0, a}$, $\Phi_2 \with \Gamma_n^{\theta, 1, a}$, $\Psi \with \Lambda_n^{\theta, a}$ in the notation of \cref{lem:linear_combination_anns_different_length}) proves for all $\theta \in \Theta$, $a \in A$ that
    \begin{align}
        \begin{split} \label{eq:mlfp_as_ann_general_Lambda_layer_structure}
            \cD (\Lambda_n^{\theta, a}) &= \cD ( \Gamma_n^{\theta, 0, a} \,\bSum_{\, \fJ}\, \Gamma_n^{\theta, 1, a} ) \\
            &= \pr[\bigg]{\dm, \dimANNlevel_1 ( \longerANN{\mathbb{L}_n, \fJ} ( \Gamma_n^{\theta, 0, a} )  ) + \dimANNlevel_1 ( \longerANN{\mathbb{L}_n, \fJ} (\Gamma_n^{\theta, 1, a}) ), \dimANNlevel_2 ( \longerANN{\mathbb{L}_n, \fJ} ( \Gamma_n^{\theta, 0, a} )  ) + \dimANNlevel_2 ( \longerANN{\mathbb{L}_n, \fJ} (\Gamma_n^{\theta, 1, a}) ), \dots\\
            &\quad \dots, \dimANNlevel_{\mathbb{L}_n - 1} ( \longerANN{\mathbb{L}_n, \fJ} ( \Gamma_n^{\theta, 0, a} )  ) + \dimANNlevel_{\mathbb{L}_n - 1} ( \longerANN{\mathbb{L}_n, \fJ} (\Gamma_n^{\theta, 1, a}) ), 1} \\
            &= \pr[\bigg]{\dm, \dimANNlevel_1 ( \longerANN{\mathbb{L}_n, \fJ} ( \Gamma_n^{0, 0, \action} )  ) + \dimANNlevel_1 ( \longerANN{\mathbb{L}_n, \fJ} (\Gamma_n^{0, 1, \action}) ), \dimANNlevel_2 ( \longerANN{\mathbb{L}_n, \fJ} ( \Gamma_n^{0, 0, \action} )  ) + \dimANNlevel_2 ( \longerANN{\mathbb{L}_n, \fJ} (\Gamma_n^{0, 1, \action}) ), \dots\\
            &\quad \dots, \dimANNlevel_{\mathbb{L}_n - 1} ( \longerANN{\mathbb{L}_n, \fJ} ( \Gamma_n^{0, 0, \action} )  ) + \dimANNlevel_{\mathbb{L}_n - 1} ( \longerANN{\mathbb{L}_n, \fJ} (\Gamma_n^{0, 1, \action}) ), 1} \\
            &= \cD (\Lambda_n^{0, \action}) \in \N^{\mathbb{L}_n + 1}.
        \end{split}
    \end{align}
    and 
    \begin{align}
        \begin{split} \label{eq:mlfp_as_ann_general_Lambda_realization_evaluated}
            \bigl( \cR_\activation (\Lambda_n^{\theta, a}) \bigr)(x) &= \bigl(\cR_\activation (\Gamma_n^{\theta, 0, a}) \bigr)(x) + \bigl(\cR_\activation (\Gamma_n^{\theta, 1, a})\bigr)(x) \\
            &= \bigg[\SmallSum{l = 0}{n - 1} \tfrac{1}{M^{n - l}} \bigl(\cR_\activation (\Psi_{n, l}^{\theta, 0, a})\bigr)(x) \bigg] + \bigg[ \SmallSum{l = 0}{n - 1} \tfrac{-\mathbbm{1}_{\N}(l)}{M^{n- l}} \bigl( \cR_\activation (\Psi_{n, \max\{l - 1, 0\}}^{\theta, 1, a}) \bigr)(x) \bigg]\\
            &= \sum_{l = 0}^{n - 1} \tfrac{1}{M^{n - l}} \sum_{i = 1}^{M^{n - l}}\Bigl( \cR_\activation(\bF) \circ \cR_\activation (\Xi_l^{(\theta, l, i)}) \circ \cR_\activation(\bX_a^{(\theta, l, i)}) \Bigr) (x) \\
            &\quad- \sum_{l = 0}^{n - 1} \tfrac{\mathbbm{1}_{\N}(l)}{M^{n - l}} \sum_{i = 1}^{M^{n - l}} \bigl( \cR_\activation (\bF) \circ \cR_\activation (\Xi_{\max\{ l - 1, 0\}}^{(\theta, -l, i)}) \circ \cR_\activation( \bX_a^{(\theta, l, i)}) \bigr) (x)\\
            &= \bigl(\varmlfp_n^{\theta}(x)\bigr)(a).
        \end{split}
    \end{align}
    This,
    \eqref{eq:mlfp_as_ann_general_Lambda_layer_structure},
    \eqref{eq:mlfp_as_ann_general_Phi},
    \cref{lem:compositions_of_anns},
    \cref{lem:parallelizations_of_anns_same_depth}
    demonstrate for all $\theta \in \Theta$, $x \in \R^\dm$ that
    \begin{align}
        \begin{split} \label{eq:mlfp_as_ann_general_Phi_layer_structure}
            \cD (\Phi_n^\theta) &= \cD (\paraANN{\abs{A}} (( \Lambda_n^{\theta, a}  )_{a \in A}) \bullet \fT_{\dm, \abs{A}} ) \\
            &= \Bigl( \dm, \SmallSum{a \in A}{} \dimANNlevel_1 (\Lambda_n^{\theta, a}), \SmallSum{a \in A}{} \dimANNlevel_2 (\Lambda_n^{\theta, a}), \dots, \SmallSum{a \in A}{} \dimANNlevel_{\mathbb{L}_n - 1} (\Lambda_n^{\theta, a}), \SmallSum{a \in A}{} \dimANNlevel_{\mathbb{L}_n} (\Lambda_n^{\theta, a}) \Bigr) \\
            &= \Bigl( \dm, \abs{A} \dimANNlevel_1(\Lambda_n^{0, \action}), \abs{A}  \dimANNlevel_1(\Lambda_n^{0, \action}), \dots, \abs{A} \dimANNlevel_{\mathbb{L}_n - 1}(\Lambda_n^{0, \action}), \abs{A} \Bigr) \in \N^{\mathbb{L}_n + 1}
        \end{split}
    \end{align}
    and
    \begin{align}
        \begin{split}
            \bigl( \cR_\activation (\Phi_n^\theta) \bigr)(x) &= \bigl( \cR_\activation (\paraANN{\abs{A}} (( \Lambda_n^{\theta, a}  )_{a \in A}) \bullet \fT_{\dm, \abs{A}} ) \bigr)(x) = \Bigl(\cR_\activation \bigl(\paraANN{\abs{A}} (( \Lambda_n^{\theta, a}  )_{a \in A})  \bigr)\Bigr) \begin{pmatrix}
                x \\ x \\ \vdots \\ x
            \end{pmatrix} \\
            &= \Bigl( \bigl( \cR_\activation (\Lambda_n^{\theta, a}) \bigr)(x) \Bigr)_{a \in A} = \Bigl( \bigl(\varmlfp_n^\theta (x)\bigr)(a) \Bigr)_{a \in A} = \varmlfp_n^\theta (x) \in \R^{|A|}.
        \end{split}
    \end{align}
    This and induction prove \cref{item:mlfp_as_ann_general_realization,item:mlfp_as_ann_general_realization_evaluated,item:mlfp_as_ann_general_layer_structure}.
    Moreover, note that \eqref{eq:mlfp_as_ann_general_def_base_case_ANN},
    \eqref{eq:mlfp_as_ann_general_Xi},
    the fact that for all $j, k \in \N$ it holds that $\cL(\fT_{k, j}) = 1$,
    \cref{lem:compositions_of_anns},
    \cref{item:certain_extensions_of_anns_item_extension_to_same_length} in \cref{lem:certain_extensions_of_anns},
    \cref{def:parallelizations_of_anns_with_different_layer_structure},
    \cref{cor:parallelizations_of_anns_with_different_layer_structure_realization_involving_identities},
    and the induction hypothesis yield for all $\theta \in \Theta$, $l \in \N \cap [1, n - 1]$, that
    \begin{equation}
        \cL (\Xi_0^\theta) = \cL (\Xi_0^0) = 1
    \end{equation}
    and
    \begin{align}
        \begin{split}
            \cL (\Xi_l^\theta) &= \cL (\Xi_l^0) = \cL \bigl( \paraLANN{2}{\mathbb{J}} (\I_\dm, \Phi_l^0) \bullet \fT_{\dm, 2} \bigr) = \cL \bigl( \paraLANN{2}{\mathbb{J}} (\I_\dm, \Phi_l^0) \bigr)\\
            &= \cL \bigl( \paraANN{2} ( \longerANN{ \max\{ \cL(\I_\dm), \cL(\Phi_l^0) \}, \I_\dm }(\I_\dm), \longerANN{ \max\{\cL(\I_\dm), \cL (\Phi_l^0)\}, \I_{\abs{A}} } (\Phi_l^0) ) \bigr) \\
            &= \max\{ \cL (\I_\dm), \cL(\Phi_l^0) \} = \max \{2, \cL (\Phi_l^0)\} \\
            &\le \max\{2, l (\cL(\bF) + \cL(\bX_\action^0) - 1) + 1\} = l \bigl( \cL (\bF) + \cL (\bX_\action^0) - 1 \bigr) + 1.
        \end{split}
    \end{align}
    This ensures for all $\theta \in \Theta$, $l \in \N_0 \cap [0, n - 1]$ that 
    \begin{equation}
        \cL (\Xi_l^\theta) = \cL (\Xi_l^0) \le l \bigl( \cL(\bF) + \cL(\bX_\action^0) - 1 \bigr) + 1.
    \end{equation}
    This, \eqref{eq:mlfp_as_ann_general_depth_nth_iterate}, \eqref{eq:mlfp_as_ann_general_Phi_layer_structure}, and \cref{lem:compositions_of_anns} establish for all $\theta \in \Theta$ that 
    \begin{align}
        \begin{split}
            \cL (\Phi_n^\theta) &= \mathbb{L}_n\\
            &= \max_{j \in \{0,1\}} \fL_n^j \\ 
            &= \max_{ \substack{l \in \{0, 1, \dots, n- 1\},\\ j \in \{0,1\}} } \cL ( \bF \bullet \Xi_{\max\{l - j, 0\}}^0 \bullet \bX_\action^0 ) \\
            &= \max_{ \substack{l \in \{0, 1, \dots, n- 1\},\\ j \in \{0,1\}} } \big[ \cL (\bF) + \cL (\bX_\action^0) + \cL (\Xi_{ \max \{l - j, 0\}}^0) - 2 \big] \\
            &\le \max_{ \substack{l \in \{0, 1, \dots, n- 1\},\\ j \in \{0,1\}} } \big[ \cL (\bF) + \cL (\bX_\action^0) - 2 + \max\{ l - j, 0 \} \bigl( \cL (\bF) + \cL (\bX_\action^0) - 1 \bigr) + 1 \big] \\
            &= \cL (\bF) + \cL (\bX_\action^0) - 2 + (n - 1)\bigl( \cL (\bF) + \cL (\bX_\action^0) - 1 \bigr) + 1\\
            &\le n \bigl( \cL(\bF) + \cL (\bX_\action^0) - 1 \bigr) + 1.
        \end{split}
    \end{align}
    This and induction prove \cref{item:mlfp_as_ann_general_depth}.
    Next, observe that \eqref{eq:mlfp_as_ann_general_Psi_layer_structure},
    \eqref{eq:mlfp_as_ann_general_Gamma_layer_structure},
    \eqref{eq:mlfp_as_ann_general_Lambda_layer_structure},
    and \cref{item:extensions_of_anns_involving_identities_item_extension_by_identity} in \cref{lem:extensions_of_anns_involving_identities}
    demonstrate for all $\theta \in \Theta$, $a \in A$ that
    \begin{align}
        \begin{split} \label{eq:mlfp_as_ann_general_Lambda_width}
            \fnorm{ \cD (\Lambda_n^{\theta, a}) } &= \fnorm{ \cD (\Lambda_n^{0, \action})} \le \fnorm{ \cD (\longerANN{\mathbb{L}_n, \fJ} (\Gamma_n^{0, 0, \action}) ) } + \fnorm{ \cD (\longerANN{\mathbb{L}_n, \fJ} (\Gamma_n^{0, 1, \action}) ) } \\
            &\le \max\big\{ \fnorm{ \cD (\fJ) }, \fnorm{ \cD(\Gamma_n^{0,0,\action})}\big\} + \max\big\{\fnorm{ \cD (\fJ) }, \fnorm{ \cD(\Gamma_n^{0,1,\action}) }\big\} \\
            &\le \max\! \bigg\{\! \fd, \SmallSum{l = 0}{n - 1} \fnorm{ \cD (\longerANN{\fL_n^0, \fJ} ( \Psi_{n, l}^{0,0,\action} )) } \!\bigg\} + \max \! \bigg\{ \! \fd, \SmallSum{l = 0}{n - 1} \fnorm{ \cD( \longerANN{ \fL_n^1, \fJ } ( \Psi_{n, l}^{0, 1, \action} ) )  }\! \bigg\} \\
            &\le \max\! \bigg\{\! \fd, \SmallSum{l = 0}{n - 1} \max\!\big\{ \fd, \fnorm{ \cD (\Psi_{n, l}^{0,0,\action}) } \big\} \! \bigg\} + \max \! \bigg\{\! \fd, \SmallSum{l = 0}{n - 1} \max\!\big\{ \fd, \fnorm{ \cD (\Psi_{n, l}^{0,1,\action}) } \big\}  \!\bigg\} \\
            &\le \SmallSum{l = 0}{n - 1} \max \big\{ \fd, \fnorm{ \cD (\Psi_{n, l}^{0,0,\action}) } \big\} + \SmallSum{l = 0}{n - 1} \max \big\{ \fd, \fnorm{ \cD (\Psi_{n, l}^{0,1,\action}) } \big\} \\
            &\le \SmallSum{l = 0}{n - 1} \max\big\{ \fd, M^{n - l} \max \big\{ \dm, \fnorm{ \cD( \bF \bullet \Xi_l^0 \bullet \bX_\action^0 ) } \big\} \big\} \\
            &\quad+ \SmallSum{l = 0}{n - 1} \max \big\{ \fd, M^{n - l} \max \big\{ \dm, \fnorm{ \cD( \bF \bullet \Xi_{\max\{l - 1, 0\}}^0 \bullet \bX_\action^0 ) } \big\} \big\} \\
            &\le \SmallSum{l = 0}{n - 1} M^{n - l} \max \big\{ \fd, \fnorm{ \cD (\bF \bullet \Xi_l^0 \bullet \bX_\action^0)  } \big\}\\
            &\quad+ \SmallSum{l = 0}{n - 1} M^{n - l} \max \big\{ \fd, \fnorm{ \cD (\bF \bullet \Xi_{ \max\{l - 1, 0\} }^0 \bullet \bX_\action^0)  } \big\}\\
            &\le \SmallSum{l = 0}{n - 1} M^{n - l} \max \big\{ \fd, \fnorm{ \cD (\bF) }, \fnorm{\cD(\bX_\action^0)}, \fnorm{\cD( \Xi_l^0) } \big\}\\
            &\quad+ \SmallSum{l = 0}{n - 1} M^{n - l} \max \big\{ \fd, \fnorm{ \cD (\bF) }, \fnorm{\cD(\bX_\action^0)}, \fnorm{ \cD(\Xi_{ \max\{ l - 1, 0 \} }^0) } \big\}.
        \end{split}
    \end{align}
    Furthermore, note that the fact that for all $j, k \in \N$ it holds that $\cD (\fT_{k, j}) = (k, jk)$,
    \eqref{eq:mlfp_as_ann_general_def_base_case_ANN},
    \eqref{eq:mlfp_as_ann_general_Xi},
    \cref{lem:compositions_of_anns},
    \cref{lem:parallelizations_of_anns_same_depth},
    and \cref{def:parallelizations_of_anns_with_different_layer_structure} ensure for all $l \in \N \cap [1, n - 1]$ that $\fnorm{\cD (\Xi_0^0)} = \dm + \abs{A} \le \fd\dm + \fd \abs{A} + \fnorm{\cD(\Phi_0^0)}$ and
    \begin{align}
        \begin{split}
            \fnorm{ \cD ( \Xi_l^0 )} &= \fnorm{ \cD ( \paraLANN{2}{ \mathbb{J}}(\I_\dm, \Phi_l^0 ) \bullet \fT_{\dm, 2} ) }\\
            &\le \max \cu*{2\dm, \fnorm{ \cD\pr{\paraLANN{2}{\mathbb{J}} (\I_\dm, \Phi_l^0)} }} \\
            &\le \max \cu*{2\dm, \fnorm[\big]{ \cD\pr{\longerANN{\max\{ \cL(\I_\dm), \cL(\Phi_l^0) \}, \I_\dm  } (\I_\dm)} } + \fnorm[\big]{ \cD\pr{\longerANN{ \max\{\cL(\I_\dm), \cL(\Phi_l^0)\}, \I_{\abs{A}} }(\Phi_l^0)} }}  \\
            &\le \max \cu*{2\dm, \fd \dm +  \max \cu*{\fd \abs{A}, \fnorm{\cD (\Phi_l^0)}} } \\
            &= \fd \dm+ \max \cu*{ \fd \abs{A}, \fnorm{\cD(\Phi_l^0)}} \le \fd \dm + \fd \abs{A} + \fnorm{\cD(\Phi_l^0)}.
        \end{split}
    \end{align}
    This implies for all $l \in \N_0 \cap [0, n - 1]$ that $\fnorm{ \cD (\Xi_l^0)} \le \fd \dm + \fd |A| + \fnorm{ \cD(\Phi_l^0)}$.
    Combining this,
    \eqref{eq:mlfp_as_ann_general_Lambda_width},
    the assumption that $M \ge 2$,
    and the induction hypothesis yields that
    \begin{align}
        \begin{split}
            \fnorm{ \cD (\Lambda_n^{0, \action}) } &\le \SmallSum{l = 0}{n - 1} M^{n - l} \max\big\{ \fd, \fnorm{\cD (\bF)}, \fnorm{ \cD (\bX_\action^0)}, \fnorm{ \cD (\Xi_l^0)} \big\}\\
            &\quad+ \SmallSum{l = 0}{n - 1} M^{n - l} \max\big\{ \fd, \fnorm{\cD (\bF)}, \fnorm{ \cD (\bX_\action^0)}, \fnorm{ \cD (\Xi_{\max\{l - 1, 0\}}^0)} \big\} \\
            &\le \SmallSum{l = 0}{n - 1} M^{n - l} \max\big\{ \fd, \fnorm{\cD (\bF)}, \fnorm{ \cD (\bX_\action^0)}, \fd \dm + \fd |A| + \fnorm{\cD(\Phi_l^0)} \big\}\\
            &\quad+ \SmallSum{l = 0}{n - 1} M^{n - l} \max\big\{ \fd, \fnorm{\cD (\bF)}, \fnorm{ \cD (\bX_\action^0)}, \fd \dm + \fd |A| + \fnorm{\cD(\Phi_{\max\{ l - 1, 0 \}}^0)} \big\} \\
            &\le \SmallSum{l = 0}{n - 1} M^{n - l} \big[ \max\big\{ \fnorm{\cD (\bF)}, \fnorm{ \cD (\bX_\action^0)}, \fd \dm + \fd |A| \big\} + \fnorm{\cD(\Phi_l^0)} \big] \\
            &\quad+ \SmallSum{l = 0}{n - 1} M^{n - l} \big[ \max\big\{ \fnorm{\cD (\bF)}, \fnorm{ \cD (\bX_\action^0)}, \fd \dm + \fd |A| \big\} + \fnorm{\cD(\Phi_{\max\{ l - 1, 0 \}}^0)} \big] \\
            &= c \bigg[\SmallSum{l = 0}{n - 1} M^{n - l}\bigg] + \bigg[ \SmallSum{l = 0}{n - 1} M^{n - l} \big[ \fnorm{ \cD (\Phi_l^0)} + \fnorm{ \cD (\Phi_{\max\{ l - 1, 0 \}}^0)} \big]  \bigg]\\
            &\le c \Bigg[ \bigg[ \SmallSum{l = 0}{n - 1} M^{n - l} \bigg] + \bigg[ \SmallSum{ l = 0 }{n - 1} M^{n - l} (4|A|M)^l \bigg] + \bigg[\SmallSum{l = 0}{n - 1} M^{n - l} (4|A|M)^{\max\{l - 1, 0\}}\bigg] \Bigg]\\
            &= c \Bigg[ M \tfrac{M^n - 1}{M - 1} + M^n \bigg[ \SmallSum{l = 0}{n - 1} (4|A|)^l \bigg] + M^n + M^{n - 1} \bigg[ \SmallSum{l = 1}{n - 1} (4|A|)^{l - 1} \bigg] \Bigg] \\
            &\le c \bigg[ 3M^n + M^n \tfrac{(4|A|)^n - 1}{4|A| - 1} + M^{n - 1} \tfrac{(4|A|)^{n - 1} - 1}{4|A| - 1}\bigg]\\
            &= c M^n \bigg[ 3 + \tfrac{(4|A|)^n - 1}{4|A| - 1} + \tfrac{1}{M} \tfrac{(4|A|)^{n - 1} - 1}{4|A| - 1} \bigg] \\
            &\le c M^n \bigg[ \tfrac{(4|A|)^n}{4|A| - 1} + \tfrac{(4|A|)^{n - 1}}{4|A| - 1} + \tfrac{3(4|A| - 1) - 2}{4|A| - 1}\bigg]\\
            &= c(4|A|M)^n \bigg[ \tfrac{1}{4|A| - 1} + \tfrac{1}{4|A|(4|A| - 1)} + \tfrac{3(4|A| - 1) - 2}{ (4|A|)^n (4|A| - 1)}\bigg]\\
            &\le c(4|A|M)^n \bigg[ \tfrac{4|A| + 1}{4|A|(4|A| - 1)} + \tfrac{3(4|A| - 1) - 2}{4|A|(4|A| - 1)} \bigg] \\
            &= c(4|A|M)^n \bigg[ \tfrac{4(4|A| - 1)}{4|A|(4|A| - 1)} \bigg] = c(4|A|M)^n \tfrac{1}{|A|}.
        \end{split}
    \end{align} 
    Combining this and \eqref{eq:mlfp_as_ann_general_Phi} demonstrates that
    \begin{equation}
        \fnorm{ \cD (\Phi_n^0)} \le |A| \fnorm{\cD (\Lambda_n^{0, \action})} \le c (4|A|M)^n.
    \end{equation}
    This and induction prove \cref{item:mlfp_as_ann_general_width}. This, \cref{item:mlfp_as_ann_general_depth}, and the fact that for all $L \in \N$, $\zeta \in \bN$, $l_0, l_1, \dots, l_L \in \N$ with $\cD (\zeta) = (l_0, l_1, \dots, l_L) \in \N^{L + 1}$ it holds that $\paramANN(\zeta) = \sum_{k = 1}^{L} l_k(l_{k - 1} + 1)$ establish that
    \begin{align}
        \begin{split}
            \paramANN(\Phi_n^0) &\le \SmallSum{k = 1}{\cL(\Phi_n^0)} \fnorm{\cD (\Phi_n^0)} \bigl( \fnorm{ \cD (\Phi_n^0)} + 1 \bigr) \le 2\cL(\Phi_n^0) \fnorm{ \cD (\Phi_n^0) }^2\\
            &\le 2 \big[ n \bigl( \cL(\bF) + \cL(\bX_\action^0) - 1 \bigr) + 1 \big] c^2 (4|A|M)^{2n}.
        \end{split}
    \end{align}  
    This and induction prove \cref{item:mlfp_as_ann_general_size}. Hence, \cref{item:mlfp_as_ann_general_realization,item:mlfp_as_ann_general_realization_evaluated,item:mlfp_as_ann_general_layer_structure,item:mlfp_as_ann_general_depth,item:mlfp_as_ann_general_width,item:mlfp_as_ann_general_size} hold true. The proof of \cref{lem:mlfp_as_ann_general} is thus complete.
\end{mproof}

\section{ANN approximations for functional fixed-point equations}
\label{sec:ANN_approximations_fixed_point_equations}

\subsection{ANN approximations for functional fixed-point equations with general activation functions}
\label{subsec:general_fp_equation}

\cfclear
\begin{theo} \label{thm:main_simple}
    Let $\Theta = \cup_{n \in \N} \Z^n$, let $A$ be a nonempty finite set,
    let $(\Omega, \mathcal{F}, \mathbb{P})$ be a probability space,
    let $\genConst \in [1, \infty)$, 
    for every $\dm \in \N$ let $ \lambda_\dm, \eta_\dm, K_\dm, L_\dm \in \pr{0,\infty}$,
    for every $\dm \in \N$, $\eps \in \err$ let $\lipFANN_{\dm, \eps} \in \pr{0,\infty}$,
    assume $\sup_{\dm \in \N} \lambda_\dm L_\dm < 1$, $\sup_{\dm \in \N} \eta_\dm L_\dm < 1$, $\sup_{(\dm, \eps) \in \N \times \err} \lambda_\dm \lipFANN_{\dm, \eps} < 1$,
    for every $\dm \in \N$ let $\w_\dm \colon \R^\dm \to \pr{0, \infty}$ be measurable,
    for every $\dm \in \N$, $x \in \pr{x_1, \dots, x_\dm} \in \R^\dm$ let $\fnorm{x} \in \R$ satisfy $\fnorm{x} = \max_{i \in \cu{ 1,\dots, \dm}} \abs{x_i}$,
    for every $\dm \in \N$ let $f_\dm \colon \R^\dm \times \R^{\abs{A}} \to \R$ be measurable,
    for every $\dm \in \N$ let $\tker_\dm^{\pr{a}} \colon \R^\dm \times \Borel(\R^\dm) \to \br{0,1}$, $a \in A$, be stochastic kernels,
    let $\action \in A$, 
    let $\ol{\omega} \in \Omega$,
    let $\activation \in C(\R, \R)$,
    let $\fJ \in \bN$ satisfy $\cR_\activation(\fJ) = \idfunc_\R$ and $\cL(\fJ) = 2$,
    let $(\bF_{\dm, \eps})_{(\dm, \eps) \in \N \times \err} \subseteq \bN$ satisfy for all $\dm \in \N$, $\eps \in \N$ that $\cR_\activation(\bF_{\dm, \eps}) \in C(\R^{\dm + \abs{A}}, \R)$,
    let $\bX^{\theta, a}_{\dm, \eps} \colon \Omega \to \bN$, $\theta \in \Theta$, $\dm \in \N$, $\eps \in \err$, $a \in A$, satisfy for all $\theta \in \Theta$, $\dm \in \N$, $\eps \in \err$, $a \in A$, $\omega \in \Omega$ that $\cD (\bX_{\dm, \eps}^{\theta, a}(\omega)) = \cD (\bX_{\dm, \eps}^{0, \action}(\ol{\omega}))$ and $\cR_\activation( \bX_{\dm, \eps}^{0, \action}(\omega) ) \in C(\R^\dm, \R^\dm)$,
    for every $\dm \in \N$, $\eps \in \err$ let $\tkerANN_{\dm, \eps}^{\pr{a}} \colon \R^\dm \times \Borel(\R^\dm) \to \br{0,1}$, $a \in A$, satisfy for all $\theta \in \Theta$, $\dm \in \N$, $\eps \in \err$, $x \in \R^\dm$, $a \in A$, $Y \in \Borel(\R^\dm)$ that $ \tkerANN_{\dm, \eps}^{\pr{a}}(x, Y) = \mathbb{P} \br[\big]{ \pr[\big]{ \cR_\activation (\bX_{\dm, \eps}^{\theta, a}  ) }(x) \in Y }$,
    for every $\dm \in \N$, $\eps \in \err$ let $(\mathcal{F}^\theta_{\dm, \eps})_{\theta \in \Theta}$ be independent sub-$\sigma$-algebras of $\mathcal{F}$,
    assume for all $\dm \in \N$, $\eps \in \err$ that $\pr{\cR_{\activation} \pr{\bX^{\theta, a}_{\dm, \eps}}}_{a \in A} \colon \R^\dm \times \Omega \to \pr{\R^\dm}^{\abs{A}}$, $\theta \in \Theta$, are i.i.d.\ random fields,
    assume for all $\theta \in \Theta$, $\dm \in \N$, $\eps \in \err$ that $\pr{\cR_{\activation} \pr{\bX^{\theta, a}_{\dm, \eps}}}_{a \in A}$ is $(\Borel(\R^\dm) \otimes \mathcal{F}^\theta_{\dm, \eps}) / \Borel((\R^\dm)^{\abs{A}})$-measurable,
    assume for all $\dm \in \N$, $\eps \in \err$ that $\max \{ \lambda_\dm, K_\dm \} \le \genConst \dm^\genConst$ and
    \begin{align} \label{eq:main_simple_poly_growth_assumption_ANN}
        \max\cu[\big]{ \cL (\bF_{\dm, \eps}), \cL(\bX_{\dm, \eps}^{0, \action}(\ol{\omega})), \fnorm{\cD (\bF_{\dm, \eps}) }, \fnorm{ \cD (\bX_{\dm, \eps}^{0, \action}(\ol{\omega})) } } \le \genConst \dm^{\genConst}\eps^{-\genConst},
    \end{align}
    assume for all $\dm \in \N$, $\eps \in \err$, $x,y \in \R^\dm$, $r, s \in \R^{\abs{A}}$, $a \in A$ that 
    \begin{align}
        \begin{split} \label{eq:main_simple_approximation_assumption}
            &\abs[\big]{ \pr[\big]{ \cR_\activation (\bF_{\dm, \eps}) }(x,r) -  f_\dm (x,r) } \le \eps \genConst \dm^{\genConst} \abs{\w_\dm(x)}, \qquad \wdist_\dm\pr[\big]{ \tkerANN_{\dm, \eps}^{\pr{a}}(x, \cdot), \tker_\dm^{\pr{a}}(x, \cdot) } \le \eps \genConst \dm^{\genConst} \abs{\w_\dm(x)},
        \end{split} \\
        \begin{split} \label{eq:main_simple_lipschitz}
            & \abs*{ f_\dm (x,r) - f_\dm(x,s) } \le L_\dm \max_{b \in A} \abs*{ r(b) - s(b) }, \qquad \abs*{f_\dm(x, r) - f_\dm(y,r)} \le K_\dm \norm{ x - y}, \\
            &\abs*{ \pr[\big]{ \cR_\activation (\bF_{\dm, \eps}) }(x,r) - \pr[\big]{ \cR_\activation (\bF_{\dm, \eps}) }(x,s) } \le \lipFANN_{\dm, \eps} \max_{b \in A} \abs*{ r(b) - s(b) },\quad \text{and} \\
            &\wdist_\dm \pr[\big]{ \tker_{\dm}^{\pr{a}}(x, \cdot), \tker_{\dm}^{\pr{a}}(y, \cdot) } \le \eta_\dm \norm{x - y},
        \end{split}
    \end{align}
    assume for all $\dm \in \N$, $\eps \in \err$, $a \in A$ that
    \begin{equation}
        \begin{split} \label{eq:main_simple_first_moment_kernels}
            \sup_{x \in \R^\dm} \abs{\w_\dm\pr{x}}^{-1} \int_{\R^\dm} \norm{y} \tker_\dm^{\pr{a}}\pr{x, \dxx y} < \infty \qquad \text{and} \qquad \sup_{x \in \R^\dm} \abs{\w_\dm \pr{x}}^{-1 } \int_{\R^\dm} \norm{y} \tkerANN_{\dm, \eps}^{\pr{a}} \pr{x, \dxx y} < \infty,
        \end{split}
    \end{equation} 
    assume for all $\dm \in \N$, $x \in \R^\dm$, $\eps \in \err$, $a \in A$ that 
    \begin{align}
        \begin{split} \label{eq:main_simple_regularity_perturbed_inputs}
            &\pr*{ \int_{\R^\dm} \abs{\w_\dm(y)}^2 \tkerANN_{\dm, \eps}^{\pr{a}}(x, \dxx y) }^{\!\!\nicefrac{1}{2}} \le \lambda_\dm \abs{ \w_\dm(x) }, \\
            &\pr*{ \int_{\R^\dm} \abs*{ \pr[\big]{ \cR_\activation(\bF_{\dm, \eps}) } (y, 0) }^2 \tkerANN_{\dm, \eps}^{\pr{a}}(x, \dxx y)    }^{\!\!\nicefrac{1}{2}} \le \genConst \dm^\genConst \abs{\w_\dm(x)},
        \end{split} \\
        \begin{split} \label{eq:main_simple_regularity_target_inputs}
            &\int_{\R^\dm} \abs{ \w_\dm (y) } \tker_\dm^{\pr{a}}(x, \dxx y) \le \lambda_\dm \abs{\w_\dm (x)}, \quad \text{and} \quad \sup_{\substack{z \in \R^\dm,\\ b \in A}} \abs{\w_\dm(z)}^{-1} \!\! \int_{\R^\dm} \abs{ f_\dm (y, 0) } \tker_\dm^{\pr{b}}(z, \dxx y) < \infty.
        \end{split}
    \end{align}
    Then there exists $c \in \R$ such that for every $\dm \in \N$, $\eps \in \err$ and every probability measure $\mu \colon \Borel\pr{\R^\dm} \to \br{0,1}$ with $\pr[]{ \int_{\R^\dm} \abs{\w_\dm(x)}^2 \mu (\dxx x) }^{\nicefrac{1}{2}} \le \genConst \dm^\genConst$ it holds that
    \begin{enumerate}[label=(\roman *)]
        \item \label{item:main_simple_ex_uniq_sol} there exists a unique measurable $u \colon \R^\dm \to \R^{\abs{A}}$ which satisfies for all $x \in \R^\dm$, $a \in A$ that 
        $\sup_{(z, b) \in \R^\dm \times A} \abs{\w_\dm(z)}^{-1} \abs{(u(z))(b)} < \infty$,
        $\int_{\R^\dm} \abs{f_\dm(y, u(y))} \tker_\dm^{\pr{a}}(x, \dxx y) < \infty$, and
        \begin{equation}
            \pr*{ u (x) } \pr{a} = \int_{\R^\dm} f_\dm \pr{ y, u(y) } \tker^{\pr{a}}_\dm(x, \dxx y) 
        \end{equation}
        and
        \item \label{item:main_simple_ex_ann} there exists $\Psi \in \bN$ which satisfies that $\cH \pr{ \Psi } > 0$, $\paramANN(\Psi) \le c \dm^c \eps^{-c}$, $\cR_\activation(\Psi) \in C(\R^\dm, \R^{\abs{A}})$, and
        \begin{equation} \label{eq:main_simple_existence_approximating_ann}
            \pr*{ \int_{\R^\dm} \max_{a \in A} \abs*{ (u(x))(a) - \pr*{\pr[\big]{\cR_\activation(\Psi)}(x)}(a) }^2 \mu(\dxx x) }^{\!\!\nicefrac{1}{2}} \le \eps
        \end{equation}
    \end{enumerate}
    \cfout.
\end{theo}

\begin{mproof}{\cref{thm:main_simple}}
    Throughout this proof let $\fd \in \N$ satisfy $\cD (\fJ) = (1, \fd, 1)$,
    let $\fK \in \R$ satisfy $\fK = \max\cu{\fc, \fd\pr{1 + \abs{A}}}$,
    let $M \in \N \cap [2, \infty)$ satisfy
    \begin{equation} \label{eq:main_simple_def_M}
        M > \pr*{\frac{  1 + \pr{\sup_{(\dm, \eps) \in \N \times \err} \lambda_\dm \lipFANN_{\dm, \eps}}(2\abs{A} - 1) }{ 1 - \pr{ \sup_{(\dm, \eps) \in \N \times \err} \lambda_\dm \lipFANN_{\dm, \eps} } }}^2,
    \end{equation}
    let $\alpha_{\dm, \eps}, \gamma_{\dm, \eps} \in [0, \infty)$, $\dm \in \N$, $\eps \in \err$, satisfy for all $\dm \in \N$, $\eps \in \err$ that
    \begin{align}
        \begin{split} \label{eq:main_simple_def_alpha_gamma_dm_eps}
            \alpha_{\dm, \eps} &= \tfrac{\lambda_\dm \lipFANN_{\dm, \eps} (1 + |A|M^{\nicefrac{-1}{2}}) + M^{\nicefrac{-1}{2}} + \sqrt{ \big( \lambda_\dm \lipFANN_{\dm, \eps} (1 + |A|M^{\nicefrac{-1}{2}}) + M^{\nicefrac{-1}{2}} \big)^2 + 4M^{\nicefrac{-1}{2}} \lambda_\dm \lipFANN_{\dm, \eps} (|A| - 1)}}{2} \qquad \text{and} \\
            \gamma_{\dm, \eps} &= \tfrac{3}{2} \max\Big\{ \tfrac{\genConst}{1 - \lambda_\dm \lipFANN_{\dm, \eps}}, \tfrac{ \abs{A} \genConst}{\abs{A} \lambda_\dm \lipFANN_{\dm, \eps} + 1} \Big\},
        \end{split}
    \end{align}
    let $\mlfppower, \mlfpconst \in \br{0, \infty}$ satisfy $\mlfppower = \sup_{(\dm, \eps) \in \N \times \err} \alpha_{\dm, \eps}$ and $\mlfpconst = \sup_{(\dm, \eps) \in \N \times \err} \gamma_{\dm, \eps}$,
    let $\scrn \colon \err \to \N \cup \{+\infty\}$ satisfy for all $\eps \in \err$ that $\scrn(\eps) = \min( \{ n \in \N : \mlfpconst \mlfppower^n \le \eps\} \cup \{+ \infty\} )$.\ %
    First, observe that the assumption that $\sup_{(\dm, \eps) \in \N \times \err} \lambda_\dm \lipFANN_{\dm, \eps} < 1$
    and \eqref{eq:main_simple_def_M} 
    ensure that $\mlfppower \in [M^{-\nicefrac{1}{2}}, 1)$.
    Next, note that the assumption that $\sup_{(\dm, \eps) \in \N \times \err} \lambda_\dm \lipFANN_{\dm, \eps} < 1$
    and the assumption that $\genConst \in [1, \infty)$ imply that $\mlfpconst \in [\nicefrac{3}{2}, \infty)$.
    This and the fact that $\mlfppower \in \pr{0,1}$ show for all $\eps \in \err$ that $\scrn(\eps) \in \N$.
    Furthermore, let $\mlfp_{\dm, \eps, n}^\theta \colon \R^\dm \times \Omega \to \R^{\abs{A}}$, $\theta \in \Theta$, $\dm \in \N$, $\eps \in \err$, $n \in \N_0$ satisfy for all $\theta \in \Theta$, $\dm \in \N$, $\eps \in \err$, $n \in \N_0$, $x \in \R^\dm$, $a \in A$ that 
    \begin{align}
        \begin{split}
            &\pr*{\mlfp_{\dm, \eps, n}^\theta(x)}(a)\\
            &= \SmallSum{l = 0}{ n -1} \tfrac{1}{M^{n- l}} \br[\bigg]{ \SmallSum{i = 1}{M^{n - l}} \br[\Big]{ \pr[\Big]{ \cR_\activation(\bF_{\dm, \eps}) } \pr[\Big]{  \pr[\big]{ \cR_\activation ( \bX_{\dm, \eps}^{(\theta, l, i), a} ) }(x)  , \mlfp_{\dm, \eps, l}^{(\theta, l, i)} \pr[\big]{ \pr[\big]{ \cR_\activation ( \bX_{\dm, \eps}^{(\theta, l, i), a} ) }(x) } }  \\
            &\quad- \mathbbm{1}_\N(l) \pr[\Big]{ \cR_\activation(\bF_{\dm, \eps}) } \pr[\Big]{  \pr[\big]{ \cR_\activation ( \bX_{\dm, \eps}^{(\theta, l, i), a} ) }(x)  , \mlfp_{\dm, \eps, \max\{l - 1, 0\}}^{(\theta, -l, i)} \pr[\big]{ \pr[\big]{ \cR_\activation ( \bX_{\dm, \eps}^{(\theta, l, i), a} ) }(x) } } }  }. 
        \end{split}
    \end{align}
    Note that
    the assumption that for all $\dm \in \N$ it holds that $\lambda_\dm L_\dm < 1$,
    the assumption that for all $\dm \in \N$ it holds that $\sup_{(x,a) \in \R^\dm \times A} \abs{\w_\dm(x)}^{-1} \!\! \int_{\R^\dm} \abs{ f_\dm (y, 0) } \tker_\dm^{\pr{a}}(x, \dxx y) < \infty$,
    the assumption that for all $\dm \in \N$, $x \in \R^\dm$, $a \in A$ it holds that $\int_{\R^\dm} \abs{ \w_\dm (y) } \tker_\dm^{\pr{a}}(x, \dxx y) \le \lambda_\dm \abs{\w_\dm (x)}$,
    and \cite[Lemma~2.2]{beck2023nonlinctrl} (applied for every $\dm \in \N$ with $c\with \lambda_\dm$, $L \with L_\dm$, $(\X, \mathcal{X}) \with (\R^\dm, \Borel(\R^\dm))$, $A \with A$, $(\tker_a)_{a \in A} \with (\tker_\dm^{\pr{a}})_{a \in A}$, $f \with f_\dm$, $\w \with \br{\R^\dm \ni x \mapsto \br{A \ni a \mapsto \w_\dm(x) \in \pr{0,\infty}} \in \pr{0,\infty}^{\abs{A}}}$ in the notation of \cite[Lemma~2.2]{beck2023nonlinctrl})
    yield that for every $\dm \in \N$ there exists a unique measurable $u_\dm \colon \R^\dm \to \R^{\abs{A}}$ which satisfies for all $x \in \R^\dm$, $a \in A$ that $\sup_{(z,b) \in \R^\dm \times A} \abs{\w_\dm(z)}^{-1} \abs{(u_\dm(z))(b)} < \infty$, $\int_{\R^\dm} \abs{ f_\dm \pr{y, u_\dm(y)} } \tker_\dm^{\pr{a}}(x, \dxx y) < \infty$, and 
    \begin{equation} \label{eq:main_simple_existence_exact_solution}
        (u_\dm(x))(a) = \int_{\R^\dm} f_\dm \pr{y, u_\dm(y)} \tker_\dm^{\pr{a}}(x, \dxx y).
    \end{equation}
    This proves \cref{item:main_simple_ex_uniq_sol}.
    Next, observe that
    the assumptions that for all $\dm \in \N$, $\eps \in \err$ it holds that $\lambda_\dm \lipFANN_{\dm, \eps} < 1$,
    the assumptions that for all $\dm \in \N$, $\eps \in \err$, $x \in \R^\dm$, $a \in A$ it holds that $\pr[\big]{ \int_{\R^\dm} \abs{\w_\dm(y)}^2 \tkerANN_{\dm, \eps}^{\pr{a}}(x, \dxx y) }^{\nicefrac{1}{2}} \le \lambda_\dm \abs{ \w_\dm(x) }$ and
    $\pr[\big]{ \int_{\R^\dm} \abs*{ \pr[\big]{ \cR_\activation(\bF_{\dm, \eps}) } (y, 0) }^2 \tkerANN_{\dm, \eps}^{\pr{a}}(x, \dxx y)    }^{\nicefrac{1}{2}} \le \genConst \dm^\genConst \abs{\w_\dm(x)}$,
    and \cite[Lemma~2.2]{beck2023nonlinctrl} (applied for every $\dm \in \N$, $\eps \in \err$ with $c \with \lambda_\dm$, $L \with \lipFANN_{\dm, \eps}$, $(\X, \mathcal{X}) \with (\R^\dm, \Borel(\R^\dm))$, $A \with A$, $(\tker_a)_{a \in A} \with (\tkerANN_{\dm, \eps}^{\pr{a}})_{a \in A}$, $f \with \cR_\activation( \bF_{\dm, \eps})$, $\w \with \br{ \R^\dm \ni x \mapsto \br{ A \ni a \mapsto \w_\dm(x) \in \pr{0,\infty} } \in \pr{0, \infty}^{\abs{A}} }$ in the notation of \cite[Lemma~2.2]{beck2023nonlinctrl})
    imply that for every $\dm \in \N$, $\eps \in \err$ there exists a unique measurable $v_{\dm, \eps} \colon \R^\dm \to \R^{\abs{A}}$ which satisfies for all $x \in \R^\dm$, $a \in A$ that $\sup_{(z,b) \in \R^\dm \times A} \abs{\w_\dm(z)}^{-1} \abs{(v_{\dm, \eps}(z))(b)} < \infty$, $\int_{\R^\dm} \abs*{ \pr{\cR_\activation (\bF_{\dm, \eps})} \pr{y, v_{\dm, \eps}(y)} } \tkerANN^{\pr{a}}_{\dm, \eps}(x, \dxx y) < \infty$, and
    \begin{equation} \label{eq:main_simple_existence_perturbed_solution_1}
        (v_{\dm, \eps}(x))(a) = \int_{\R^\dm} \pr{\cR_\activation (\bF_{\dm, \eps})} \pr{y, v_{\dm, \eps}(y)} \tkerANN_{\dm, \eps}^{\pr{a}}(x, \dxx y).
    \end{equation}
    This and 
    the assumption that for all $\theta \in \Theta$, $\dm \in \N$, $\eps \in \err$, $x \in \R^\dm$, $a \in A$, $Y \in \Borel(\R^\dm)$ it holds that $ \tkerANN_{\dm, \eps}^{\pr{a}}(x, Y) = \mathbb{P} \br[\big]{ \pr[\big]{ \cR_\activation (\bX_{\dm, \eps}^{\theta, a}  ) }(x) \in Y }$
    ensure for all $\dm \in \N$, $\eps \in \err$, $x \in \R^\dm$, $a \in A$ that 
    \begin{align}
        \begin{split} \label{eq:main_simple_existence_perturbed_solution_2}
            (v_{\dm, \eps}(x))(a) &= \int_{\R^\dm} \pr{\cR_\activation (\bF_{\dm, \eps})} \pr{y, v_{\dm, \eps}(y)} \tkerANN_{\dm, \eps}^{\pr{a}}(x, \dxx y) \\
            &= \E \br*{ \pr[\big]{\cR_\activation(\bF_{\dm, \eps})} \pr*{ \pr*{\cR_\activation(\bX^{0,a}_{\dm, \eps})}(x) , v_{\dm, \eps}\pr*{\pr*{\cR_\activation(\bX^{0,a}_{\dm, \eps})}(x)} } }.
        \end{split}
    \end{align}
    Moreover, observe that the triangle inequality demonstrates for all $\dm, n \in \N$, $\eps \in \err$, $x \in \R^\dm$ that
    \begin{align}
        \begin{split}
            &\pr[\Big]{\E \br[\Big]{ \max_{a \in A} \abs*{ (u_\dm(x))(a) - \pr*{ \mlfp_{\dm, \eps, n}^{0}(x) }(a) }^2  }}^{\!\!\nicefrac{1}{2}} \\
            &\le \pr[\bigg]{ \E \br[\bigg]{ \pr[\Big]{ \max_{a \in A} \abs*{ (u_\dm(x))(a) - (v_{\dm, \eps}(x))(a) } + \max_{b \in A} \abs*{ (v_{\dm, \eps}(x))(b) - \pr*{ \mlfp_{\dm, \eps, n}^0(x) }(b) }  }^2 } }^{\!\!\nicefrac{1}{2}} \\
            &\le \pr[\Big]{ \E \br[\Big]{ \max_{a \in A} \abs*{ (u_\dm(x))(a) - (v_{\dm, \eps}(x))(a) }^2 } }^{\!\!\nicefrac{1}{2}} + \pr[\Big]{ \E \br[\Big]{\max_{a \in A} \abs*{ (v_{\dm, \eps}(x))(a) - \pr*{ \mlfp_{\dm, \eps, n}^0(x) }(a) }^2} }^{\!\!\nicefrac{1}{2}}\\
            &= \max_{a \in A} \abs*{ (u_\dm(x))(a) - (v_{\dm, \eps}(x))(a) } + \pr*{ \E\br[\Big]{\max_{a \in A} \abs*{ (v_{\dm, \eps}(x))(a) - \pr*{ \mlfp_{\dm, \eps, n}^0(x) }(a) }^2} }^{\!\!\nicefrac{1}{2}}.
        \end{split}
    \end{align}
    Hence, it holds for all $\dm, n \in \N$, $\eps \in \err$, $x \in \R^\dm$ that
    \begin{align}
        \begin{split} \label{eq:main_simple_error_decomposition}
            &\E \br[\Big]{ \max_{a \in A} \abs*{ (u_\dm(x))(a) - \pr*{ \mlfp_{\dm, \eps, n}^{0}(x) }(a) }^2  }\\
            &\le \pr*{ \max_{a \in A} \abs*{ (u_\dm(x))(a) - (v_{\dm, \eps}(x))(a) } + \pr[\Big]{ \E \br[\Big]{\max_{a \in A} \abs*{ (v_{\dm, \eps}(x))(a) - \pr*{ \mlfp_{\dm, \eps, n}^0(x) }(a) }^2} }^{\!\!\nicefrac{1}{2}} }^{\!\!2} \\
            &= \abs*{\w_\dm(x)}^2 \pr*{ \max_{a \in A} \tfrac{\abs{ (u_\dm(x))(a) - (v_{\dm, \eps}(x))(a) }}{\abs{\w_\dm(x)}}  + \pr*{\tfrac{ \E \br*{ \max_{a \in A} \abs*{ (v_{\dm, \eps}(x))(a) - \pr*{ \mlfp_{\dm, \eps, n}^0(x) }(a) }^2 } }{\abs{\w_\dm(x)}^2}}^{\!\!\nicefrac{1}{2}}   }^{\!\!2}.
        \end{split}
    \end{align}
    Furthermore, note that
    the assumptions that for all $\dm \in \N$, $\eps \in \err$, $x,y \in \R^\dm$, $r,s \in \R^{\abs{A}}$ it holds that
    $\abs*{ f_\dm (x,r) - f_\dm(x,s) } \le L_\dm \max_{b \in A} \abs*{ r(b) - s(b) }$,
    $\abs*{f_\dm(x, r) - f_\dm(y,r)} \le K_\dm \norm{ x - y}$,
    $\abs*{ \pr[\big]{ \cR_\activation (\bF_{\dm, \eps}) }(x,r) - \pr[\big]{ \cR_\activation (\bF_{\dm, \eps}) }(x,s) } \le \lipFANN_{\dm, \eps} \max_{b \in A} \abs*{ r(b) - s(b) }$, and
    $\wdist_\dm \pr[\big]{ \tker_{\dm}^{\pr{a}}(x, \cdot), \tker_{\dm}^{\pr{a}}(y, \cdot) } \le \eta_\dm \norm{x - y}$,
    the assumption that for all $\dm \in \N$, $\eps \in \err$, $x \in \R^\dm$, $a \in A$ it holds that $\pr[\big]{ \int_{\R^\dm} \abs*{ \pr[\big]{ \cR_\activation(\bF_{\dm, \eps}) } (y, 0) }^2 \tkerANN_{\dm, \eps}^{\pr{a}}(x, \dxx y)    }^{\nicefrac{1}{2}} \le \genConst \dm^\genConst \abs{\w_\dm(x)}$,
    \eqref{eq:main_simple_existence_exact_solution},
    \eqref{eq:main_simple_existence_perturbed_solution_2},
    and \cref{item:stability_kernel_w_result} in \cref{prop:stability_kernel_w} (applied for every $\dm \in \N$, $\eps \in \err$ with $A \with A$, $\eta \with \eta_\dm$, $c \with \lambda_\dm$, $K \with K_\dm$, $L_1 \with \lipFANN_{\dm, \eps}$, $L_2 \with L_\dm$, $\w \with \w_\dm$, $f_1 \with \cR_\activation(\bF_{\dm, \eps})$, $f_2 \with f_\dm$, $(\tker_1^{\pr{a}})_{a \in A} \with (\tkerANN_{\dm, \eps}^{\pr{a}})_{a \in A}$, $(\tker_2^{\pr{a}})_{a \in A} \with (\tker_\dm^{\pr{a}})_{a \in A}$, $u_1 \with v_{\dm, \eps}$, $u_2 \with u_\dm$ in the notation of \cref{prop:stability_kernel_w})
    yield for all $\dm \in \N$, $\eps \in \err$ that 
    \begin{align}
        \begin{split} \label{eq:main_simple_stability_error}
            \sup_{(x, a) \in \R^\dm \times A} \!\!\!\!\! \tfrac{\abs*{ (v_{\dm, \eps}(x))(a) - (u_\dm(x))(a)}}{\abs{ \w_\dm (x)}} &\le \tfrac{\lambda_\dm}{1 - \lambda_\dm (\min\{ L_\dm, \lipFANN_{\dm, \eps} \})} \br*{ \sup_{(y,r) \in \R^\dm \times \R^{\abs{A}}} \!\!\!\!\! \tfrac{\abs*{ \pr*{ \cR_\activation( \bF_{\dm, \eps}) } (y, r) - f_\dm(y,r) } }{\abs{\w_\dm (y)}} } \\
            &\hspace{0.5cm}+ \tfrac{K_\dm}{(1 - \eta_\dm L_\dm)(1 - \lambda_\dm \lipf_\dm)} \br*{ \sup_{(y, b) \in \R^\dm \times A}  \!\!\!\!\! \tfrac{\wdist_\dm \pr[]{ \tkerANN_{\dm, \eps}^{\pr{b}} (y, \cdot), \tker_\dm^{\pr{b}}(y, \cdot) }}{\abs{ \w_\dm(y)}} }.
        \end{split}
    \end{align}
    Next, note that
    the assumption that for all $\dm \in \N$, $\eps \in \err$ it holds that $\lambda_\dm \lipFANN_{\dm, \eps} < 1$,
    the assumption that for all $\dm \in \N$, $\eps \in \err$, $x \in \R^\dm$, $a \in A$ it holds that $\pr{ \int_{\R^\dm} \abs{\w_\dm(y)}^2 \tkerANN_{\dm, \eps}^{\pr{a}}(x, \dxx y) }^{\nicefrac{1}{2}} \le \lambda_\dm \abs{ \w_\dm(x) }$
    and $\pr{ \int_{\R^\dm} \abs{ \pr[\big]{ \cR_\activation(\bF_{\dm, \eps}) } (y, 0) }^2 \tkerANN_{\dm, \eps}^{\pr{a}}(x, \dxx y)}^{\nicefrac{1}{2}} \le \genConst \dm^\genConst \abs{\w_\dm(x)}$, 
    the fact that for all $\dm \in \N$, $\eps \in \err$ it holds that $ \sup_{(x,a) \in \R^\dm \times A} \abs{\w_\dm (x)}^{-1} \abs{(v_{\dm, \eps}(x))(a)} < \infty$,
    the assumption that for all $\dm \in \N$, $\eps \in \err$, $x \in \R^\dm$, $r,s \in \R^{\abs{A}}$ it holds that $\abs*{ \pr{ \cR_\activation (\bF_{\dm, \eps}) }(x,r) - \pr{ \cR_\activation (\bF_{\dm, \eps}) }(x,s) } \le \lipFANN_{\dm, \eps} \max_{b \in A} \abs*{ r(b) - s(b) }$,
    \eqref{eq:main_simple_existence_perturbed_solution_2},
    and \cite[Lemma~2.4]{beck2023nonlinctrl} (applied for every $\dm \in \N$, $\eps \in \err$ with
    $c_f \with \genConst \dm^\genConst$,
    $c_\w \with \lambda_\dm$,
    $L \with \lipFANN_{\dm, \eps}$, $(\X, \mathcal{X}) \with (\R^\dm, \Borel(\R^\dm))$,
    $(\Omega, \mathcal{F}, \mathbb{P}) \with (\Omega, \mathcal{F}, \mathbb{P})$,
    $A \with A$,
    $X \with \pr{ \cR_\activation \pr{ \bX_{\dm, \eps}^{0, a} } }_{a \in A}$,
    $f \with \cR_\activation(\bF_{\dm, \eps})$,
    $\w \with \br{\R^\dm \ni x \mapsto \br{A \ni a \mapsto \w_\dm(x) \in \pr{0, \infty} } \in \pr{0,\infty}^{\abs{A}}}$,
    $v \with v_{\dm, \eps}$ in the notation of \cite[Lemma~2.4]{beck2023nonlinctrl}) prove that for all $\dm \in \N$, $\eps \in \err$ it holds that
    \begin{align}
        \begin{split}
            &\sup_{(x,a) \in \R^\dm \times A} \!\!\!\! \tfrac{\abs*{ (v_{\dm, \eps}(x))(a) }}{\abs{\w_\dm(x)}} \le \tfrac{\genConst \dm^\genConst }{1- \lambda_{\dm} \lipFANN_{\dm, \eps}} \qquad \text{and}\\
            &\sup_{(x, a) \in \R^\dm \times A} \!\!\!\! \tfrac{\pr*{ \E \br*{ \max_{b \in A} \abs*{ \pr*{v_{\dm, \eps} \pr*{ \pr{\cR_\activation \pr{ \bX_{\dm, \eps}^{0, a} }}\pr{x} }}(b) }^2 }  }^{\!\!\nicefrac{1}{2}}  }{\abs{\w_\dm(x)}} \le \tfrac{\lambda_\dm \genConst \dm^\genConst }{1 - \lambda_\dm \lipFANN_{\dm, \eps}}.
        \end{split}
    \end{align}
    Combining this,
    the assumption that for all $\dm \in \N$, $\eps \in \err$, $x \in \R^\dm$, $a \in A$ it holds that $\pr{ \int_{\R^\dm} \abs{\w_\dm(y)}^2 \tkerANN_{\dm, \eps}^{\pr{a}}(x, \dxx y) }^{\nicefrac{1}{2}} \le \lambda_\dm \abs{ \w_\dm(x) }$
    and $\pr{ \int_{\R^\dm} \abs{ \pr[\big]{ \cR_\activation(\bF_{\dm, \eps}) } (y, 0) }^2 \tkerANN_{\dm, \eps}^{\pr{a}}(x, \dxx y)}^{\nicefrac{1}{2}} \le \genConst \dm^\genConst \abs{\w_\dm(x)}$,
    the assumption that for all $\dm \in \N$, $\eps \in \err$, $x \in \R^\dm$, $r,s \in \R^{\abs{A}}$ it holds that
    \begin{equation}
        \begin{split}
            \abs*{ \pr{ \cR_\activation (\bF_{\dm, \eps}) }(x,r) - \pr{ \cR_\activation (\bF_{\dm, \eps}) }(x,s) } \le \lipFANN_{\dm, \eps} \max_{b \in A} \abs*{ r(b) - s(b) },
        \end{split}
    \end{equation}
    \eqref{eq:main_simple_existence_perturbed_solution_2},
    and \cite[Proposition~3.8]{beck2023nonlinctrl} (applied for every $\dm \in \N$, $\eps \in \err$ with
    $M \with M$,
    $\Theta \with \Theta$,
    $(\X, \mathcal{X}) \with ( \R^\dm, \Borel(\R^\dm))$,
    $A \with A$,
    $\w \with \w_\dm$, $f \with \cR_\activation(\bF_{\dm, \eps})$,
    $v \with v_{\dm, \eps}$,
    $(\mathcal{F}^\theta)_{\theta \in \Theta} \with (\mathcal{F}^{\theta}_{\dm, \eps})_{\theta \in \Theta}$,
    $(X^{\theta})_{\theta \in \Theta} \with \pr{ \pr{ \cR_\activation\pr{\bX_{\dm, \eps}^{\theta, a} }}_{a \in A} }_{\theta \in \Theta} $,
    $(V^\theta_n)_{(\theta, n) \in \Theta \times \N_0} \with (\mlfp^\theta_{\dm, \eps, n})_{(\theta, n) \in \Theta \times \N_0}$,
    $c \with \gamma_{\dm, \eps}$,
    $c_\w \with \lambda_\dm$,
    $L \with \lipFANN_{\dm, \eps}$,
    $c_f \with \genConst \dm^\genConst$,
    $c_v \with \tfrac{\lambda_\dm \genConst \dm^\genConst}{1- \lambda_\dm \lipFANN_{\dm, \eps}}$ in the notation of \cite[Proposition~3.8]{beck2023nonlinctrl})
    demonstrates for all $\dm \in \N$, $\eps \in \err$, $n \in \N$, $x \in \R^\dm$ that
    \begin{align} \label{eq:main_simple_mlfp_error} 
        \pr*{\frac{ \E \br[\big]{ \max_{a \in A} \abs{ (v_{\dm, \eps}(x))(a) - \pr{\mlfp^0_{\dm, \eps,n}(x)}(a) }^2 } }{ \abs*{\w_\dm(x)}^2 }}^{\!\!\nicefrac{1}{2}} \le \dm^\genConst \gamma_{\dm, \eps} \alpha_{\dm, \eps}^n \le \dm^\genConst \mlfpconst \mlfppower^n.
    \end{align}
    This, \eqref{eq:main_simple_error_decomposition},
    and \eqref{eq:main_simple_stability_error},
    ensure for every $\dm, n \in \N$, $\eps \in \err$ and every probability measure $\mu \colon \Borel\pr{\R^\dm} \to \br{0,1}$ with $\pr{\int_{\R^\dm} \abs{\w_\dm(x)}^2 \mu( \dxx x)}^{\nicefrac{1}{2}} \le \genConst d^\genConst$ that
    \begin{align}
        \begin{split}
            &\pr[\bigg]{ \int_{\R^\dm} \E \br[\Big]{ \max_{a \in A} \abs[\big]{ (u_{\dm}(x))(a) - \pr*{ \mlfp_{\dm, \eps, n}^0(x) }(a) }^2 } \mu (\dxx x)  }^{\!\!\nicefrac{1}{2}} \\
            &\le \pr[\bigg]{ \dm^\genConst \mlfpconst \mlfppower^n + \tfrac{\lambda_\dm}{1 - \lambda_\dm (\min\{ L_\dm, \lipFANN_{\dm, \eps} \})}\br[\bigg]{ \sup_{(y,r) \in \R^\dm \times \R^{\abs{A}}} \!\!\!\!\!\!\! \tfrac{ \abs{ (\cR_\activation(\bF_{\dm, \eps}))(y,r) - f_\dm (y,r) } }{\abs{\w(y)}}   }\\
            &\quad+ \tfrac{K_\dm}{(1 - \eta_\dm L_\dm)(1- \lambda_\dm L_\dm)} \br[\bigg]{ \sup_{(y,b) \in \R^\dm \times A} \!\!\!\!\!\!\! \tfrac{\wdist_\dm\pr*{ \tkerANN_{\dm, \eps}^{\pr{b}}(y, \cdot), \tker_{\dm}^{\pr{b}}(y, \cdot) }}{\abs{\w_\dm(y)}} } } \pr*{\int_{\R^\dm} \abs{\w_\dm(x)}^2 \mu\pr{\dxx x}}^{\!\!\nicefrac{1}{2}}\\
            &\le \pr[\bigg]{ \dm^\genConst \mlfpconst \mlfppower^n + \tfrac{\lambda_\dm}{1 - \lambda_\dm (\min\{ L_\dm, \lipFANN_{\dm, \eps} \})}\br[\bigg]{ \sup_{(y,r) \in \R^\dm \times \R^{\abs{A}}} \!\!\!\!\!\!\! \tfrac{ \abs{ (\cR_\activation(\bF_{\dm, \eps}))(y,r) - f_\dm (y,r) } }{\abs{\w(y)}}   }\\
            &\quad+ \tfrac{K_\dm}{(1 - \eta_\dm L_\dm)(1- \lambda_\dm L_\dm)} \br[\bigg]{ \sup_{(y,b) \in \R^\dm \times A} \!\!\!\!\!\!\! \tfrac{\wdist_\dm\pr*{ \tkerANN_{\dm, \eps}^{\pr{b}}(y, \cdot), \tker_{\dm}^{\pr{b}}(y, \cdot) }}{\abs{\w_\dm(y)}} } } \genConst \dm^\genConst.
        \end{split}
    \end{align}
    Combining this,
    the assumption that $\genConst \in [1, \infty)$,
    the assumptions that for all $\dm \in \N$, $\eps \in \err$, $x \in \R^\dm$, $a \in A$, $r \in \R^{\abs{A}}$ it holds that $\abs*{ \pr*{ \cR_\activation(\bF_{\dm, \eps}) }(x,r) - f_\dm(x,r) } \le \eps \genConst \dm^\genConst \abs{\w_\dm(x)}$ and $\wdist_\dm\pr[\big]{ \tkerANN_{\dm, \eps}^{\pr{a}}(x, \cdot), \tker_\dm^{\pr{a}}(x, \cdot) } \le \eps \genConst \dm^\genConst \abs{ \w_\dm(x) }$,
    and the assumption that for all $\dm \in \N$ it holds that $\max\cu{\lambda_\dm, K_\dm} \le \genConst \dm^\genConst$
    establishes for every $\dm \in \N$, $\delta, \eps \in \err$ and every probability measure $\mu \colon \Borel\pr{\R^\dm} \to \br{0,1}$ with $\pr{\int_{\R^\dm} \abs{\w_\dm(x)}^2 \mu( \dxx x)}^{\nicefrac{1}{2}} \le \genConst d^\genConst$ that
    \begin{align}
        \begin{split}
            &\pr[\bigg]{ \int_{\R^\dm} \E \br[\Big]{ \max_{a \in A} \abs[\big]{ (u_{\dm}(x))(a) - \pr*{ \mlfp_{\dm, \eps, \scrn(\delta)}^0(x) }(a) }^2 } \mu (\dxx x)  }^{\!\!\nicefrac{1}{2}}\\
            &\le \genConst \dm^\genConst \pr[\bigg]{ \dm^\genConst \mlfpconst \mlfppower^{\scrn (\delta)} + \tfrac{\lambda_\dm}{1 - \lambda_\dm (\min\{ L_\dm, \lipFANN_{\dm, \eps} \})}\br[\bigg]{ \sup_{(y,r) \in \R^\dm \times \R^{\abs{A}}} \!\!\!\!\!\!\! \tfrac{ \abs{ (\cR_\activation(\bF_{\dm, \eps}))(y,r) - f_\dm (y,r) } }{\abs{\w(y)}}   }\\
            &\quad+ \tfrac{K_\dm}{(1 - \eta_\dm L_\dm)(1- \lambda_\dm L_\dm)} \br[\bigg]{ \sup_{(y,b) \in \R^\dm \times A} \!\!\!\!\!\!\! \tfrac{\wdist_\dm\pr*{ \tkerANN_{\dm, \eps}^{\pr{b}}(y, \cdot), \tker_{\dm}^{\pr{b}}(y, \cdot) }}{\abs{\w_\dm(y)}} }} \\
            &\le \genConst \dm^\genConst \pr*{ \dm^\genConst \delta +  \tfrac{\eps \genConst^2 \dm^{2\genConst} }{1 - \lambda_\dm (\min\{ L_\dm, \lipFANN_{\dm, \eps} \})} + \tfrac{\eps \genConst^2 \dm^{2\genConst}}{(1 - \eta_\dm L_\dm)(1- \lambda_\dm L_\dm)}}\\
            &\le \genConst^3\dm^{3\genConst} \pr*{\delta +\tfrac{\eps}{1 - \lambda_\dm (\min\{ L_\dm, \lipFANN_{\dm, \eps} \})} + \tfrac{\eps }{(1 - \eta_\dm L_\dm)(1- \lambda_\dm L_\dm)}}.
        \end{split}
    \end{align}
    This
    ensures for every $\dm \in \N$, $\delta \in \err$ and every probability measure $\mu \colon \Borel\pr{\R^\dm} \to \br{0,1}$ with $\pr{\int_{\R^\dm} \abs{\w_\dm(x)}^2 \mu( \dxx x)}^{\nicefrac{1}{2}} \le \genConst d^\genConst$ that 
    \begin{equation} \label{eq:main_simple_intermediate_error}
        \begin{split}        
            &\pr[\bigg]{ \int_{\R^\dm} \E \br[\Big]{ \max_{a \in A} \abs[\big]{ (u_{\dm}(x))(a) - \pr*{ \mlfp_{\dm, \delta, \scrn(\delta)}^0(x) }(a) }^2 } \mu (\dxx x)  }^{\!\!\nicefrac{1}{2}} \\
            &\le \delta \genConst^3 \dm^{3\genConst} \pr*{ 1 +  \tfrac{1}{1 - \lambda_\dm (\min\{ L_\dm, \lipFANN_{\dm, \delta} \})} + \tfrac{1}{(1 - \eta_\dm L_\dm)(1- \lambda_\dm L_\dm)}} \\
            &\le \delta \genConst^3 \dm^{3\genConst} \pr*{1 + \frac{1}{1 - \lambda_\dm L_\dm} + \frac{1}{\pr{1 - \eta_\dm L_\dm} \pr{1 - \lambda_\dm L_\dm}}}.
        \end{split}
    \end{equation}
    Throughout the remainder of the proof let $\pr*{\delta_{\dm, \eps}}_{\pr{\delta, \eps} \in \N \times \err} \subseteq \pr{0,\infty}$ satisfy for all $\dm \in \N$, $\eps \in \err$ that
    \begin{align}
        \begin{split}\label{eq:main_simple_definition_delta_dm_eps}
            \delta_{\dm, \eps} &= \frac{\eps}{\genConst^3 \dm^{3\genConst}} \pr*{ 1 + \frac{1}{(1 - \lambda_\dm L_\dm)} + \frac{1}{(1 - \lambda_{\dm} L_\dm)(1 - \eta_\dm L_\dm)} }^{-1}.
        \end{split}
    \end{align}
    Note that the assumptions that for all $\dm \in \N$ it holds that $\lambda_\dm L_\dm < 1$ and $\eta_\dm L_\dm < 1$
    and the assumption that $\genConst \ge 1$
    ensure that for all $\dm \in \N$, $\eps \in \err$ it holds that $\delta_{\dm, \eps} \in \err$.
    Combining this with \eqref{eq:main_simple_intermediate_error} demonstrates for every $\dm \in \N$, $\eps \in \err$ and every probability measure $\mu \colon \Borel\pr{\R^\dm} \to \br{0,1}$ with $\pr{\int_{\R^\dm} \abs{\w_\dm(x)}^2 \mu( \dxx x)}^{\nicefrac{1}{2}} \le \genConst d^\genConst$ that
    \begin{equation}
        \begin{split}
            &\int_{\R^\dm} \E \br[\Big]{ \max_{a \in A} \abs[\big]{ (u_{\dm}(x))(a) - \pr[\big]{ \mlfp_{\dm, \delta_{\dm, \eps}, \scrn(\delta_{\dm, \eps})}^0(x) }(a) }^2 } \mu (\dxx x)\\
            &\le \pr*{\delta_{\dm, \eps} \genConst^3 \dm^{3\genConst} \pr[\Big]{1 + \tfrac{1}{1 - \lambda_\dm L_\dm} + \tfrac{1}{\pr{1 - \eta_\dm L_\dm} \pr{1 - \lambda_\dm L_\dm}}}}^2 \le \eps^2.
        \end{split}
    \end{equation}
    This and Fubini's theorem establish for every $\dm \in \N$, $\eps \in \err$ and every probability measure $\mu \colon \Borel\pr{\R^\dm} \to \br{0,1}$ with $\pr{\int_{\R^\dm} \abs{\w_\dm(x)}^2 \mu( \dxx x)}^{\nicefrac{1}{2}} \le \genConst d^\genConst$ that 
    \begin{align}
        \begin{split}
            &\E \br*{ \int_{\R^\dm} \max_{a \in A} \abs[\big]{ (u_{\dm}(x))(a) - \pr[\big]{ \mlfp_{\dm, \delta_{\dm, \eps}, \scrn(\delta_{\dm, \eps})}^0(x) }(a) }^2 \mu(\dxx x) }\\
            &= \int_{\R^\dm} \E \br[\Big]{ \max_{a \in A} \abs[\big]{ (u_{\dm}(x))(a) - \pr[\big]{ \mlfp_{\dm, \delta_{\dm, \eps}, \scrn(\delta_{\dm, \eps})}^0(x) }(a) }^2 } \mu (\dxx x) \le \eps^2.
        \end{split}
    \end{align}
    Hence, for every $\dm \in \N$, $\eps \in \err$ and every probability measure $\mu \colon \Borel\pr{\R^\dm} \to \br{0,1}$ with $\pr{\int_{\R^\dm} \abs{\w_\dm(x)}^2 \mu( \dxx x)}^{\nicefrac{1}{2}} \le \genConst d^\genConst$ there exists $\omega_{\dm, \eps, \mu} \in \Omega$, which is assumed to be fixed for the remainder of this proof, such that
    \begin{equation}
        \begin{split} \label{eq:main_simple_error_of_realization_of_mlfp}
            \int_{\R^\dm} \max_{a \in A} \abs*{ \pr{u_\dm\pr{x}}\pr{a} - \pr*{\mlfp_{\dm, \delta_{\dm, \eps}, \scrn\pr{\delta_{\dm, \eps}}} \pr{x, \omega_{\dm, \eps, \mu}} }\pr{a} }^2 \mu\pr{\dxx x} \le \eps^2.
        \end{split}
    \end{equation}
    Next, observe that \cref{lem:mlfp_as_ann_general} (applied for every $\dm \in \N$, $\eps \in \err$, $\omega \in \Omega$ with
    $\dm \with \dm$,
    $\fd \with \fd$,
    $M \with M$, 
    $\Theta \with \Theta$,
    $A \with A$,
    $\action \with \action$,
    $\activation \with \activation$,
    $\fJ \with \fJ$,
    $\bF \with \bF_{\dm, \eps}$,
    $\pr{ \bX_{a}^\theta }_{\pr{\theta, a} \in \Theta \times A} \with \pr{ \bX_{\dm, \eps}^{\theta, a} \pr{\omega} }_{\pr{\theta, a} \in \Theta \times A}$,
    $\pr{\varmlfp_n^\theta}_{\pr{\theta, n} \in \Theta \times \N_0} \with \pr{\mlfp_{\dm, \eps, n}^{\theta} \pr{\cdot, \omega}}_{\pr{\theta, n} \in \Theta \times \N_0}$
    in the notation of \cref{lem:mlfp_as_ann_general}) establishes that for all $\dm \in \N$, $\eps \in \err$, $\omega \in \Omega$ there exists $\Phi_{\dm, \eps, n}^\theta\pr{\omega} \in \N$, $\theta \in \Theta$, $n \in \N_0$ such that 
    \begin{enumerate}[label = (\Roman *)] 
        \item \label{item:main_simple_mlfp_as_ann_realization} it holds for all $\theta \in \Theta$, $n \in \N_0$ that $ \cR_\activation \pr{ \Phi_{\dm, \eps, n}^\theta \pr{\omega} } \in C\pr{\R^\dm, \R^{\abs{A}}}$,
        \item \label{item:main_simple_mlfp_as_ann_realization_evaluated} it holds for all $\theta \in \Theta$, $n \in \N_0$, $x \in \R^\dm$ that $\pr{ \cR_{\activation}\pr{\Phi_{\dm, \eps, n}^\theta \pr{\omega}} } \pr{x} = \mlfp_{\dm, \eps, n}^\theta \pr{x, \omega}$,
        \item \label{item:main_simple_mlfp_as_ann_architecture} it holds for all $\theta \in \Theta$, $n \in \N_0$ that $\cD \pr{\Phi_{\dm, \eps, n}^\theta \pr{\omega}} = \cD \pr{\Phi_{\dm, \eps, n}^0 \pr{\omega}}$,
        \item \label{item:main_simple_mlfp_as_ann_depth} it holds for all $\theta \in \Theta$, $n \in \N_0$ that $\cL\pr{ \Phi_{\dm, \eps, n}^\theta \pr{\omega} } \le n \pr*{\cL \pr{\bF_{\dm, \eps}} + \cL \pr{\bX_{\dm, \eps}^{0, \action}\pr{\omega}} - 1} + 1$,
        \item \label{item:main_simple_mlfp_as_ann_width} it holds for all $\theta \in \Theta$, $n \in \N_0$ that 
        \begin{equation}\label{eq:main_simple_mlfp_as_ann_width}
            \fnorm[\big]{\cD \pr[\big]{ \Phi_{\dm, \eps, n}^\theta \pr{\omega} } } \le 2 \max \cu[\big]{ \fnorm{\cD \pr{ \bF_{\dm, \eps} }}, \fnorm{ \cD\pr{ \bX_{\dm, \eps}^{0, \action} \pr{\omega} } }, \fd d + \fd \abs{A} } \pr{4 \abs{A} M}^{2n},
        \end{equation}
        and
        \item \label{item:main_simple_mlfp_as_ann_size} it holds for all $\theta \in \Theta$, $n \in \N_0$ that 
        \begin{equation}
            \begin{split}
                \paramANN \pr{ \Phi_{\dm, \eps, n}^{\theta} \pr{\omega}} &\le 2 \pr[\big]{ n \pr{ \cL \pr{\bF_{\dm, \eps}} + \cL \pr{ \bX_{\dm, \eps}^{0, \action} \pr{\omega} } - 1} + 1}\\
                &\quad\cdot \pr[\big]{2 \max \cu[\big]{ \fnorm{\cD \pr{ \bF_{\dm, \eps} }}, \fnorm{ \cD\pr{ \bX_{\dm, \eps}^{0, \action} \pr{\omega} } }, \fd d + \fd \abs{A} }}^2 \pr{4\abs{A} M}^{2 n}.                
            \end{split}
        \end{equation}
    \end{enumerate}
    For every $\dm \in \N$, $\eps \in \err$ and every probability measure $\mu \colon \Borel\pr{\R^\dm} \to \br{0,1}$ with $\pr{\int_{\R^\dm} \abs{\w_\dm(x)}^2 \mu( \dxx x)}^{\nicefrac{1}{2}} \le \genConst d^\genConst$ let $\Psi_{\dm, \eps, \mu} \in \bN$ satisfy that
    \begin{equation} \label{eq:main_simple_definition_Psi}
        \Psi_{\dm, \eps, \mu} = \Phi^0_{\dm, \delta_{\dm, \eps}, \scrn\pr{\delta_{\dm, \eps}}} \pr{\omega_{\dm, \eps, \mu}}.
    \end{equation}
    Combining this,
    \cref{eq:main_simple_error_of_realization_of_mlfp},
    and \cref{item:main_simple_mlfp_as_ann_realization,item:main_simple_mlfp_as_ann_realization_evaluated}
    demonstrates for every $\dm \in \N$, $\eps \in \err$, $y \in \R^\dm$, $b \in A$ and every probability measure $\mu \colon \Borel\pr{\R^\dm} \to \br{0,1}$ with $\pr{\int_{\R^\dm} \abs{\w_\dm(x)}^2 \mu( \dxx x)}^{\nicefrac{1}{2}} \le \genConst d^\genConst$
    that $\cR_\activation \pr{ \Psi_{\dm, \eps, \mu}} \in C\pr{\R^\dm, \R^{\abs{A}}}$, $\pr*{\pr[\big]{\cR_\activation \pr{ \Psi_{\dm, \eps, \mu} } }\pr{y}}\pr{b} = \pr[\big]{\mlfp_{\dm, \delta_{\dm, \eps}, \scrn\pr{\delta_{\dm, \eps}}}^0\pr{y, \omega_{\dm, \eps, \mu}}}\pr{b}$, and
    \begin{equation}
        \begin{split} \label{eq:main_simple_error_Psi}
            &\pr*{ \int_{\R^\dm} \max_{a \in A} \abs*{ \pr{u_\dm\pr{x}}\pr{a} - \pr*{\pr[\big]{\cR_\activation \pr{ \Psi_{\dm, \eps, \mu} } }\pr{x}}\pr{a} }^2 \mu\pr{\dxx x} }^{\!\!\nicefrac{1}{2}} \\
            &\le \pr*{ \int_{\R^\dm} \max_{a \in A} \abs*{ \pr{u_\dm \pr{x}}\pr{a} - \pr*{\mlfp_{\dm, \delta_{\dm, \eps}, \scrn\pr{\delta_{\dm, \eps}}} \pr{x, \omega_{\dm, \eps, \mu}} }\pr{a} }^2 \mu\pr{\dxx x} }^{\!\!\nicefrac{1}{2}} \le \eps. 
        \end{split}
    \end{equation}
    This establishes \cref{eq:main_simple_existence_approximating_ann}.
    Furthermore, note that \cref{eq:main_simple_definition_Psi},
    \cref{item:main_simple_mlfp_as_ann_size},
    the fact that for all $\eps \in \err$ it holds that $\scrn \pr{\eps} \in \N$,
    and the assumption that for all $\theta \in \Theta$, $\dm \in \N$, $\eps \in \err$, $a \in A$, $\omega \in \Omega$ it holds that $\cD\pr{ \bX_{\dm, \eps}^{\theta, a}\pr{\omega} } = \cD\pr{\bX_{\dm, \eps}^{0, \action} \pr{\ol{\omega}}}$
    yield that for every $\dm \in \N$, $\eps \in \err$ and every probability measure $\mu \colon \Borel\pr{\R^\dm} \to \br{0,1}$ with $\pr{\int_{\R^\dm} \abs{\w_\dm(x)}^2 \mu( \dxx x)}^{\nicefrac{1}{2}} \le \genConst d^\genConst$ it holds that
    \begin{equation}
        \begin{split} \label{eq:main_simple_effort_base_estimate}
            \paramANN \pr{\Psi_{\dm, \eps, \mu}} &\le 2 \pr[\big]{ \scrn \pr{\delta_{\dm, \eps}} \pr{ \cL \pr{\bF_{\dm, \delta_{\dm, \eps}}} + \cL \pr{ \bX_{\dm, \delta_{\dm, \eps}}^{0, \action} \pr{\ol{\omega}} } - 1} + 1}\\
            &\quad\cdot \pr[\big]{2 \max \cu[\big]{ \fnorm{\cD \pr{ \bF_{\dm, \delta_{\dm, \eps}} }}, \fnorm{ \cD\pr{ \bX_{\dm, \delta_{\dm, \eps}}^{0, \action} \pr{\ol{\omega}} } }, \fd d + \fd \abs{A} }}^2 \pr{4\abs{A} M}^{2 \scrn \pr{\delta_{\dm, \eps}}}\\
            &\le 8 \pr[\big]{ \max \cu[\big]{ \fnorm{\cD \pr{ \bF_{\dm, \delta_{\dm, \eps}} }}, \fnorm{ \cD\pr{ \bX_{\dm, \delta_{\dm, \eps}}^{0, \action} \pr{\ol{\omega}} } }, \fd d + \fd \abs{A} } }^2  \\
            &\quad \cdot \pr[\big]{\cL \pr{\bF_{\dm, \delta_{\dm, \eps}}} + \cL \pr{ \bX_{\dm, \delta_{\dm, \eps}}^{0, \action} \pr{\ol{\omega}} }} \scrn \pr{\delta_{\dm, \eps}} \pr{4\abs{A} M}^{2 \scrn \pr{\delta_{\dm, \eps}}}.
        \end{split}
    \end{equation}
    Next, observe that for all $\dm \in \N$ it holds that $\fd \dm + \fd \abs{A} \le \pr{\fd + \fd\abs{A}} \dm \le \fK \dm$.
    This, 
    the assumption that $\max\cu{\genConst, \fd} \ge 1$,
    the assumption that $A$ is nonempty,
    and the assumption that for all $\dm \in \N$, $\eps \in \err$ it holds that $\max\cu[\big]{ \cL (\bF_{\dm, \eps}), \cL(\bX_{\dm, \eps}^{0, \action}(\ol{\omega})), \fnorm{\cD (\bF_{\dm, \eps}) }, \fnorm{ \cD (\bX_{\dm, \eps}^{0, \action}(\ol{\omega})) } } \le \genConst \dm^{\genConst}\eps^{-\genConst}$
    assure that for all $\dm \in \N$, $\eps \in \err$ it holds that
    \begin{equation}
        \begin{split} \label{eq:main_simple_width_depth_polynomial_growth}
            &\pr[\big]{\cL \pr{\bF_{\dm, \eps}}} + \cL \pr{ \bX_{\dm, \eps}^{0, \action} \pr{\ol{\omega}} } \pr[\big]{ \max \cu[\big]{ \fnorm{\cD \pr{ \bF_{\dm, \eps} }}, \fnorm{ \cD\pr{ \bX_{\dm, \eps}^{0, \action} \pr{\ol{\omega}} } }, \fd d + \fd \abs{A} } }^2 \\
            &\le 2\genConst \dm^\genConst \eps^{-\genConst} \pr{ \fK \dm^\genConst \eps^{-\genConst} }^2 \le 2\fK^3 \dm^{3\genConst} \eps^{-3\genConst}.
        \end{split}
    \end{equation}
    Moreover, note that \cite[Lemma~3.15]{beck2024overcoming} (applied with $m \with 1$,
    $\alpha \with \mlfppower$,
    $\beta \with (4\abs{A}M)^2$,
    $\kappa_1 \with \mlfpconst$,
    $\kappa_2 \with 1$,
    $ (c_n)_{n \in \N} \with \pr[\big]{((4\abs{A}M)^2)^n}_{n \in \N}$,
    $(e_n)_{n \in \N} \with (\mlfpconst \mlfppower^n)_{n \in \N}$,
    $N \with \scrn$ in the notation of \cite[Lemma~3.15]{beck2024overcoming})
    and the fact that $\gamma \in [\nicefrac{3}{2}, \infty)$
    yield
    for all $\eps \in \err$ that
    \begin{equation} \label{eq:main_simple_effort_power_of_M}
        \begin{split}    
            \pr*{ (4 \abs{A} M)^2 }^{\scrn( \eps) } &\le (4 \abs{A} M)^2 \max\Big\{ 1, \mlfpconst^{\frac{2 \ln( 4 \abs{A} M )}{\ln( \nicefrac{1}{\mlfppower})}} \Big\} \eps^{- \frac{2 \ln( 4 \abs{A} M )}{\ln( \nicefrac{1}{\mlfppower})}}\\
            &= (4 \abs{A} M)^2 \mlfpconst^{\frac{2 \ln( 4 \abs{A} M )}{\ln( \nicefrac{1}{\mlfppower})}} \eps^{- \frac{2 \ln( 4 \abs{A} M )}{\ln( \nicefrac{1}{\mlfppower})}}.
        \end{split}
    \end{equation}
    Furthermore,
    observe that the fact that for all $\eps \in \err$ it holds that $\scrn(\eps) = \min\{ n \in \N : \mlfpconst \mlfppower^n \le \eps \}$ implies for all $\eps \in \err$ with $\scrn(\eps) \in \N \cap [2, \infty)$ that $\mlfpconst \mlfppower^{\scrn(\eps) - 1} > \eps$.
    Hence, it holds for all $\eps \in \err$ with $\scrn(\eps) \in \N \cap [2, \infty)$ that $\scrn(\eps) < \frac{\ln( \eps^{-1} ) + \ln( \mlfpconst )}{\ln( \nicefrac{1}{\mlfppower})} + 1$.
    This,
    the fact that $\mlfppower \in \pr{0,1}$,
    and the fact that $\mlfpconst \in [\nicefrac{3}{2}, \infty)$
    establish that for all $\eps \in \err$ it holds that 
    \begin{equation}
        \begin{split} \label{eq:main_simple_iterate_number_poly_bound}
            \scrn\pr{\eps} \le \pr*{\frac{\ln\pr{\nicefrac{1}{\mlfppower}} + \ln \pr{\mlfpconst} + 1}{\ln \pr{ \nicefrac{1}{\mlfppower}}}} \eps^{-1}.
        \end{split}
    \end{equation}
    Combining this,
    \eqref{eq:main_simple_effort_base_estimate},
    \eqref{eq:main_simple_width_depth_polynomial_growth},
    and \eqref{eq:main_simple_effort_power_of_M}
    ensures for every $\dm \in \N$, $\eps \in \err$ and every probability measure $\mu \colon \Borel\pr{\R^\dm} \to \br{0,1}$ with $\pr{\int_{\R^\dm} \abs{\w_\dm(x)}^2 \mu( \dxx x)}^{\nicefrac{1}{2}} \le \genConst d^\genConst$ that
    \begin{equation}
        \begin{split} \label{eq:main_simple_effort_intermediate_estimate}
            \paramANN\pr{\Psi_{\dm, \eps, \mu}} &\le 16 \fK^3\dm^{3\genConst}\delta_{\dm, \eps}^{-3\genConst} \pr*{\frac{\ln\pr{\nicefrac{1}{\mlfppower}} + \ln \pr{\mlfpconst} + 1}{\ln \pr{ \nicefrac{1}{\mlfppower}}}} \delta_{\dm, \eps}^{-1} \pr{4\abs{A}M}^2 \gamma^{\tfrac{2 \ln \pr{4\abs{A}M}}{\ln \pr{\nicefrac{1}{\alpha}}}} \delta_{\dm, \eps}^{- \tfrac{2 \ln \pr{ 4 \abs{A} M}}{\ln \pr{\nicefrac{1}{\alpha}}}} \\
            &= 256 \fK^3 (\abs{A}M)^2 \gamma^{\tfrac{2 \ln \pr{4\abs{A}M}}{\ln \pr{\nicefrac{1}{\alpha}}}} \pr*{\frac{\ln\pr{\nicefrac{1}{\mlfppower}} + \ln \pr{\mlfpconst} + 1}{\ln \pr{ \nicefrac{1}{\mlfppower}}}} \dm^{3\genConst} \delta_{\dm, \eps}^{- \pr*{1 + 3\genConst + \tfrac{2 \ln \pr{ 4\abs{A} M }}{\ln \pr{\nicefrac{1}{\alpha}}}}}.
        \end{split}
    \end{equation}
    Moreover, observe the assumption that $\sup_{\dm \in \N} \lambda_{\dm} L_\dm < 1$,
    the assumption that $\sup_{\dm \in \N} \eta_\dm L_\dm < 1$,
    and \cref{eq:main_simple_definition_delta_dm_eps}
    establish for all $\dm \in \N$, $\eps \in \err$ that
    \begin{equation}
        \begin{split}
            \delta_{\dm, \eps}^{-1} \le \tfrac{3\genConst^3}{\pr{1 - \sup_{u \in \N}\lambda_u L_u} \pr{1 - \sup_{u \in \N} \eta_u L_u}} \dm^{3\genConst} \eps^{-1}.
        \end{split}
    \end{equation}
    Combining this with \eqref{eq:main_simple_effort_intermediate_estimate} demonstrates for every $\dm \in \N$, $\eps \in \err$ and every probability measure $\mu \colon \Borel\pr{\R^\dm} \to \br{0,1}$ with $\pr{\int_{\R^\dm} \abs{\w_\dm(x)}^2 \mu( \dxx x)}^{\nicefrac{1}{2}} \le \genConst d^\genConst$ that
    \begin{equation}
        \begin{split}
            \paramANN\pr{\Psi_{\dm, \eps, \mu}} &\le 256 \fK^3 \pr{\abs{A}M}^2 \gamma^{\tfrac{2 \ln \pr{4\abs{A}M}}{\ln \pr{\nicefrac{1}{\alpha}}}} \pr*{\frac{\ln\pr{\nicefrac{1}{\mlfppower}} + \ln \pr{\mlfpconst} + 1}{\ln \pr{ \nicefrac{1}{\mlfppower}}}} \\
            &\quad\cdot \pr*{\tfrac{3\genConst^3}{\pr{1 - \sup_{u \in \N}\lambda_u L_u} \pr{1 - \sup_{u \in \N} \eta_u L_u}}}^{\pr*{1 + 3\genConst + \tfrac{2 \ln \pr{ 4\abs{A} M }}{\ln \pr{\nicefrac{1}{\alpha}}}}} \dm^{3\genConst + 3\genConst \pr*{1 + 3\genConst + \tfrac{2 \ln \pr{ 4\abs{A} M }}{\ln \pr{\nicefrac{1}{\alpha}}}}} \\
            &\quad \cdot \eps^{-\pr*{1 + 3\genConst + \tfrac{2 \ln \pr{ 4\abs{A} M }}{\ln \pr{\nicefrac{1}{\alpha}}}}}.
        \end{split}
    \end{equation}
    This proves \cref{item:main_simple_ex_ann}.
    The proof of \cref{thm:main_simple} is thus complete.
\end{mproof}

\subsection{ANN approximations for functional fixed-point equations with specific activation functions and one-step transitions}
\label{subsec:activation_transition}

\cfclear
\begin{lemma}
    \label{lem:properties_wdist}
    \cfconsiderloaded{lem:properties_wdist}
    Let $\dm, m \in \N$,
    let $\eta \in [0, \infty)$,
    let $(\Omega, \cF, \mathbb{P})$ be a probability space,
    let $\xi \colon \Omega \to \R^m$ be measurable,
    let $\ft_i \colon \R^\dm \times \R^m \to \R^\dm$, $i \in \cu{1,2}$, be measurable,
    assume that for all $x, y \in \R^\dm$, $p \in \R^m$ it holds that $\norm{\ft_1(x, p) - \ft_1(y, p)} \le \eta \norm{x - y}$,
    let $\kappa_i \colon \R^\dm \times \Borel(\R^\dm) \to \br{0,1}$, $i \in \cu{1,2}$, satisfy for all $i \in \cu{1,2}$, $x \in \R^\dm$, $B \in \Borel(\R^\dm)$ that $\kappa_i\pr{x, B} = \mathbb{P}\br{\ft_i\pr{x, \xi} \in B}$,
    and assume for all $i \in \cu{1,2}$, $x \in \R^\dm$ that $\E \br{\norm{\ft_i\pr{x, \xi}}} < \infty$.
    Then 
    \begin{enumerate}[label = (\roman*)]
        \item \label{item:properties_wdist_error} it holds for all $x \in \R^\dm$ that
        \begin{equation}
            \wdist_\dm \pr{\kappa_1(x, \cdot), \kappa_2(x, \cdot)} \le \int_{\R^m} \norm{\ft_1 (x, p) - \ft_2(x,p)} \pr{\mathbb{P} \circ \xi^{-1}}\pr{\dxx p}
        \end{equation}
        and
        \item \label{item:properties_wdist_lipschitz} it holds for all $x, y \in \R^\dm$ that
        \begin{equation}
            \wdist_\dm \pr{ \kappa_1(x, \cdot), \kappa_1(y, \cdot) } \le \eta \norm{x - y}
        \end{equation}
    \end{enumerate}
    \cfout .
\end{lemma}

\begin{mproof}{\cref{lem:properties_wdist}}
    Let $X_i \colon \R^\dm \times \Omega \to \R^\dm$, $i \in \cu{1,2}$, satisfy for all $i \in \cu{1,2}$, $x \in \R^\dm$ that $X_i^x = \ft_i\pr{x, \xi}$.
    Note that the assumption that for all $x \in \R^\dm$, $i \in \cu{1,2}$ it holds that $\E \br{ \norm{ \ft_i\pr{x, \xi}}} < \infty$ ensures for all $i \in \cu{1,2}$, $x \in \R^\dm$ that $\kappa_i \pr{x, \cdot} \in \meas_\dm$.
    Observe that for all $i \in \cu{1,2}$, $x \in \R^\dm$, it holds that $X^x_i \sim \kappa_i(x, \cdot)$.
    Furthermore, it holds for all $x \in \R^\dm$ that $\mathbb{P} \circ (X_1^x, X_2^x)^{-1} \in \coup_\dm \pr{\kappa_1\pr{x, \cdot}, \kappa_2\pr{x, \cdot}}$.
    This demonstrates for all $x \in \R^\dm$ that
    \begin{equation}
        \begin{split}
            \wdist_\dm \pr{\kappa_1(x, \cdot), \kappa_2(x, \cdot)} &= \inf_{\gamma \in \coup_\dm \pr{ \kappa_1(x, \cdot), \kappa_2(x, \cdot) }} \int_{\R^\dm \times \R^\dm} \norm{\bfx - \bfy} \gamma\pr{\dxx \pr{\bfx, \bfy}} \\
            &\le \int_{\R^\dm \times \R^\dm} \norm{\bfx - \bfy} \pr*{\mathbb{P} \circ (X_1^x, X_2^x)^{-1}} \pr{ \dxx (\bfx, \bfy)} \\
            &= \int_{\Omega} \norm*{ X^x_1(\omega) - X^x_2(\omega)} \mathbb{P}\br{\dxx \omega} \\
            &= \int_{\Omega} \norm*{ \ft_1(x, \xi(\omega)) - \ft_2(x, \xi(\omega)) } \mathbb{P}\br{\dxx \omega}\\
            &= \int_{\R^m} \norm*{\ft_1(x, p) - \ft_2(x, p)} \pr*{\mathbb{P} \circ \xi^{-1}} \pr{\dxx p}.
        \end{split}
    \end{equation}
    This proves \cref{item:properties_wdist_error}.
    Moreover, it holds for all $x, y \in \R^\dm$ that $\mathbb{P} \circ (X_1^x, X_1^y)^{-1} \in \coup_\dm\pr{ \kappa_1(x, \cdot), \kappa_1(y, \cdot)}$.
    Combining this with the assumption that for all $x, y \in \R^\dm$, $p \in \R^m$ it holds that $\norm{\ft_1(x,p) - \ft_1(y, p)} \le \eta \norm{x-y}$
    yields for all $x, y \in \R^\dm$ that
    \begin{equation}
        \begin{split}
            \wdist_\dm \pr{ \kappa_1(x, \cdot), \kappa_1(y, \cdot) } &= \inf_{\gamma \in \coup_\dm \pr{ \kappa_1(x, \cdot), \kappa_1(y, \cdot) }} \int_{\R^\dm \times \R^\dm} \norm{\bfx - \bfy} \gamma\pr{\dxx \pr{\bfx, \bfy}} \\
            &\le \int_{\R^m} \norm*{\ft_1(x, p) - \ft_1(y, p)} \pr*{\mathbb{P} \circ \xi^{-1}} \pr{\dxx p} \\
            &\le \eta \norm{x - y}.
        \end{split}
    \end{equation} 
    This proves \cref{item:properties_wdist_lipschitz}.
    The proof of \cref{lem:properties_wdist} is thus complete.
\end{mproof}

\cfclear
\begin{corollary}
    \label{cor:activation_transition}
    \cfconsiderloaded{cor:activation_transition}
    Let $A$ be a nonempty finite set,
    let $(\Omega, \mathcal{F}, \mathbb{P})$ be a probability space,
    let $\genConst \in [1, \infty)$, 
    let $\nu \in \cu{0,1}$,
    let $\beta \in [0, \infty) \backslash \cu{1}$,
    for every $\dm \in \N$ let $ \lambda_\dm, \eta_\dm, K_\dm, L_\dm \in \pr{0,\infty}$,
    for every $\dm \in \N$, $\eps \in \err$ let $\lipFANN_{\dm, \eps} \in \pr{0,\infty}$,
    assume $\sup_{\dm \in \N} \lambda_\dm L_\dm < 1$, $\sup_{\dm \in \N} \eta_\dm L_\dm < 1$, $\sup_{(\dm, \eps) \in \N \times \err} \lambda_\dm \lipFANN_{\dm, \eps} < 1$,
    for every $\dm \in \N$ let $\w_\dm \colon \R^\dm \to \pr{0, \infty}$ be measurable,
    for every $\dm \in \N$, $x \in \pr{x_1, \dots, x_\dm} \in \R^\dm$ let $\fnorm{x} \in \R$ satisfy $\fnorm{x} = \max_{i \in \cu{1,\dots, \dm}} \abs{x_i}$,
    for every $\dm \in \N$ let $f_\dm \colon \R^\dm \times \R^{\abs{A}} \to \R$ be measurable,
    for every $\dm \in \N$, $a \in A$ let $t^{\pr{a}}_\dm \colon \R^{2\dm} \to \R^\dm$ be measurable,
    for every $\dm \in \N$ let $\xi_\dm \colon \Omega \to \R^\dm$ be measurable,
    let $\action \in A$, 
    let $\activation \colon \R \to \R$ satisfy for all $x \in \R$ that $\activation(x) = \nu \max\cu{x, \beta x} + \pr{1 - \nu} \ln\pr{1 + \exp(x)}$,
    let $(\bF_{\dm, \eps})_{(\dm, \eps) \in \N \times \err} \subseteq \bN$ satisfy for all $\dm \in \N$, $\eps \in \N$ that $\cR_\activation(\bF_{\dm, \eps}) \in C(\R^{\dm + \abs{A}}, \R)$,
    let $(\ft_{\dm, \eps}^{\pr{a}})_{(\dm, \eps, a) \in \N \times \err \times A} \subseteq \bN$ satisfy for all $\dm \in \N$, $\eps \in \err$, $a \in A$ that $\cD\pr{ \ft_{\dm, \eps}^{\pr{a}}} = \cD \pr{ \ft_{\dm, \eps}^{\pr{\action}}}$ and $\cR_\activation\pr{\ft_{\dm, \eps}^{\pr{a}}} \in C(\R^{2\dm}, \R^\dm)$,
    assume for all $\dm \in \N$, $\eps \in \err$ that $\max\cu{\lambda_\dm, K_\dm} \le \genConst \dm^\genConst$ and
    \begin{equation}
        \begin{split}
            \max\cu*{\cL\pr{ \bF_{\dm, \eps} }, \cL\pr{ \ft_{\dm, \eps}^{\pr{\action}} }, \fnorm[\big]{\cD\pr{ \bF_{\dm, \eps} }}, \fnorm[\big]{ \cD\pr{ \ft_{\dm, \eps}^{\pr{\action}} } } } \le \genConst \dm^\genConst \eps^{- \genConst},
        \end{split}
    \end{equation}
    assume for all $\dm \in \N$, $\eps \in \err$, $x,y, p \in \R^\dm$, $r, s \in \R^{\abs{A}}$, $a \in A$ that 
    \begin{align}
        \begin{split} \label{eq:activation_transition_approximation_assumption}
            &\abs[\big]{ \pr[\big]{ \cR_\activation (\bF_{\dm, \eps}) }(x,r) -  f_\dm (x,r) } \le \eps \genConst \dm^{\genConst} \abs{\w_\dm(x)},\\
            &\norm[\big]{ \pr[\big]{\cR_\activation \pr{\ft_{\dm, \eps}^{\pr{a}}}} \pr{x, p} - t_\dm^{\pr{a}}\pr{x, p} } \le \eps \genConst \dm^{\genConst} \abs{\w_\dm(x)},
        \end{split} \\
        \begin{split} \label{eq:activation_transition_lipschitz}
            & \abs*{ f_\dm (x,r) - f_\dm(x,s) } \le L_\dm \max_{b \in A} \abs*{ r(b) - s(b) }, \qquad \abs*{f_\dm(x, r) - f_\dm(y,r)} \le K_\dm \norm{ x - y}, \\
            &\abs*{ \pr[\big]{ \cR_\activation (\bF_{\dm, \eps}) }(x,r) - \pr[\big]{ \cR_\activation (\bF_{\dm, \eps}) }(x,s) } \le \lipFANN_{\dm, \eps} \max_{b \in A} \abs*{ r(b) - s(b) },\quad \text{and} \\
            &\norm[\big]{ t_\dm^{\pr{a}}\pr{ x, p } - t_\dm^{\pr{a}}\pr{ y, p }  }  \le \eta_\dm \norm{x - y},
        \end{split}
    \end{align}
    assume for all $\dm \in \N$, $\eps \in \err$, $a \in A$ that
    \begin{equation}
        \begin{split}
            \sup_{x \in \R^\dm} \abs{\w_\dm\pr{x}}^{-1} \E \br[\big]{ \norm{t_\dm^{\pr{a}} \pr{x, \xi_\dm}} } < \infty \quad \text{ and } \quad \sup_{x \in \R^\dm } \abs{\w_\dm\pr{x}}^{-1} \E \br[\big]{ \norm[\big]{ \pr[\big]{ \cR_\activation \pr{ \ft_{\dm, \eps}^{\pr{a}} } } \pr{x, \xi_\dm} }  } < \infty,
        \end{split}
    \end{equation}
    assume for all $\dm \in \N$, $x \in \R^\dm$, $\eps \in \err$, $a \in A$ that
    \begin{align}
        \begin{split}
            &\pr*{ \E\br*{ \abs*{ \w_\dm \pr*{ \pr[\big]{ \cR_\activation \pr{ \ft_{\dm, \eps}^{\pr{a}}} } \pr{x, \xi_\dm} } }^2 } }^{\!\!\nicefrac{1}{2}} \le \lambda_\dm \abs{\w_\dm\pr{x}}, \\
            &\pr*{ \E \br*{ \abs*{ \pr[\big]{ \cR_\activation \pr{\bF_{\dm, \eps}} } \pr[\big]{ \pr[\big]{ \cR_\activation \pr{\ft_{\dm, \eps}^{\pr{a}}  } } \pr{x, \xi_\dm} , 0  }  }^2 } }^{\!\!\nicefrac{1}{2}} \le \genConst \dm^\genConst \abs{\w_\dm\pr{x}},
        \end{split}\\
        \begin{split}
            & \E \br[\Big]{ \abs[\big]{ \w_\dm \pr[\big]{ t_\dm^{\pr{a}}(x, \xi_\dm) } } } \le \lambda_\dm \abs{\w_\dm\pr{x}}, \quad \text{ and } \quad \sup_{\substack{y \in \R^\dm \\ b \in A}} \abs{\w_\dm\pr{y}}^{-1} \E \br[\Big]{ \abs[\big]{ f_\dm \pr{ t_\dm^{\pr{b}} \pr{y, \xi_\dm} , 0 } } } < \infty.
        \end{split}
    \end{align}
    Then there exists $c \in \R$ such that for every $\dm \in \N$, $\eps \in \err$ and every probability measure $\mu \colon \Borel\pr{\R^\dm} \to \br{0,1}$ with $\pr{\int_{\R^\dm} \abs{\w_\dm(x)}^2 \mu( \dxx x)}^{\nicefrac{1}{2}} \le \genConst d^\genConst$ it holds that
    \begin{enumerate}[label = (\roman*)]
        \item \label{item:activation_transition_ex_uniq_sol} there exists a unique measurable $u \colon \R^\dm \to \R^{\abs{A}}$ which satisfies for all $x \in \R^\dm$, $a \in A$ that 
        $\sup_{(z, b) \in \R^\dm \times A} \abs{\w_\dm(z)}^{-1} \abs{(u(z))(b)} < \infty$,
        $\E \br[\Big]{ \abs[\big]{ f_\dm \pr[\big]{ t_\dm^{\pr{a}} \pr{x, \xi_\dm}, u \pr{ t_\dm^{\pr{a}} \pr{x, \xi_\dm} } } } } < \infty$, and
        \begin{equation}
            \begin{split}
               \pr{u\pr{x}}\pr{a} = \E \br*{ f_\dm \pr[\big]{ t_\dm^{\pr{a}} \pr{x, \xi_\dm}, u \pr{ t_\dm^{\pr{a}} \pr{x, \xi_\dm} } } }
            \end{split}
        \end{equation}
        and
        
        \item \label{item:activation_transition_ex_ann} there exists $\Psi \in \bN$ which satisfies $\cH \pr{\Psi} > 0$, $\paramANN(\Psi) \le c \dm^c \eps^{-c}$, $\cR_\activation(\Psi) \in C(\R^\dm, \R^{\abs{A}})$, and
        \begin{equation} \label{eq:activation_transition_existence_approximating_ann}
            \pr*{ \int_{\R^\dm} \max_{a \in A} \abs*{ (u(x))(a) - \pr[]{\pr[]{\cR_\activation(\Psi)}(x)}(a) }^2 \mu(\dxx x) }^{\!\!\nicefrac{1}{2}} \le \eps
        \end{equation}
    \end{enumerate}
    \cfout.
\end{corollary}

\begin{mproof}{\cref{cor:activation_transition}}
    Throughout this proof let $\Theta = \cup_{n \in \N} \Z^n$,
    for every $\dm \in \N$ let $\xi_\dm^\theta \colon \Omega \to \R^\dm$, $\theta \in \Theta$, be i.i.d. random variables,
    assume for all $\dm \in \N$, $\theta \in \Theta$ that $\xi_\dm$ and $\xi_\dm^\theta$ are identically distributed,
    for every $\dm \in \N$, $\theta \in \Theta$ let $\F_\dm^\theta = \sigma \pr{\xi_\dm^\theta}$,
    let $\ol{\omega} \in \Omega$,
    let $\bI, \fJ \in \bN$ satisfy
    \begin{align}
        \bI &= \pr*{ \pr*{ \begin{pmatrix}
                1 \\ - 1
            \end{pmatrix}, \begin{pmatrix}
                0 \\ 0
            \end{pmatrix} }, \pr[\Big]{ \begin{pmatrix}
                1 & -1
            \end{pmatrix}, 0 }} \in \pr{\R^{2\times 1} \times \R^2} \times \pr{\R^{1 \times 2} \times \R^1} \qquad \text{and}\\
        \fJ &= \begin{cases}
            \bI &: \nu = 0 \\
            (1 + \beta)^{-1} \circledast \bI &: \nu = 1,
        \end{cases}
    \end{align}
    and let $\Xi_\dm^\theta \colon \Omega \to \bN$, $\theta \in \Theta$, $\dm \in \N$, satisfy for every $\theta \in \Theta$, $\dm \in \N$, $\omega \in \Omega$ that
    \begin{equation} \label{eq:activation_transition_def_shock}
        \Xi_\dm^\theta(\omega) = \pr*{ \begin{pmatrix}
            \id_\dm \\ 0
        \end{pmatrix}, \begin{pmatrix}
            0 \\ \xi_\dm^\theta(\omega)
        \end{pmatrix} } \in \R^{2\dm \times \dm} \times \R^{2\dm}.
    \end{equation}
    First, observe that \cite[Lemma~3.5]{ackermann2023deep} and \cite[Lemma~3.8]{ackermann2023deep} ensure for all $x \in \R$ that $\cD(\fJ) = (1,2,1) \in \N^3$ and $\pr{\cR_\activation (\fJ)}(x) = x$.
    Next, note that \eqref{eq:activation_transition_def_shock} shows for all $\dm \in \N$, $\theta \in \Theta$, $x \in \R^\dm$, $\omega \in \Omega$ that
    $ \cD \pr{ \Xi_{\dm}^\theta \pr{\omega}} = \pr{\dm, 2 \dm} = \cD \pr{ \Xi_\dm^0 \pr{ \ol{\omega}}} \in \N^2 $, $\cR_\activation \pr{ \Xi^\theta_\dm \pr{\omega} } \in C\pr{\R^\dm, \R^{2\dm}} $, and $\pr[\big]{ \cR_\activation \pr{ \Xi_\dm^\theta \pr{\omega}} }(x) = \pr{x, \xi^\theta_\dm \pr{\omega}}$.
    This,
    the assumption that for all $\dm \in \N$, $\eps \in \err$, $a \in A$ it holds that $\cD \pr{\ft_{\dm, \eps}^{\pr{a}}} = \cD \pr{\ft_{\dm, \eps}^{\pr{\action}}}$,
    the assumption that for all $\dm \in \N$, $\eps \in \err$ it holds that $\max \cu{ \cL\pr{\ft_{\dm, \eps}^{\pr{a}}}, \fnorm{\cD \pr{\ft_{\dm, \eps}^{\pr{a}} } }} \le \genConst \dm^\genConst \eps^{-\genConst}$,
    and \cref{lem:compositions_of_anns} ensure for all $\theta \in \Theta$, $\dm \in \N$, $\eps \in \err$, $x \in \R^\dm$, $a \in A$, $\omega \in \Omega$ that
    $\cD \pr{ \ft_{\dm, \eps}^{\pr{a}} \bullet \Xi_{\dm}^\theta \pr{\omega} } = \cD \pr{ \ft_{\dm, \eps}^{\pr{\action}} \bullet \Xi_\dm^0\pr{\ol{\omega}} }$,
    $\cL \pr{ \ft_{\dm, \eps}^{\pr{\action}} \bullet \Xi_\dm^0\pr{\ol{\omega}} } = \cL \pr{ \ft_{\dm, \eps}^{\pr{\action}} } + \cL \pr{ \Xi_{\dm}^{0} \pr{\ol{\omega}} } - 1 = \cL \pr{ \ft_{\dm, \eps}^{\pr{\action}} } $,
    $\fnorm{ \cD \pr{ \ft_{\dm, \eps}^{\pr{\action}} \bullet \Xi_\dm^0\pr{\ol{\omega}}} } \le \max \cu{ \fnorm{ \cD\pr{ \ft_{\dm, \eps}^{\action} } }, \fnorm{\cD \pr{ \Xi_\dm^0\pr{ \ol{\omega} } }} } = \fnorm{ \cD\pr{ \ft_{\dm, \eps}^{\pr{\action}} } }$,
    $ \cR_\activation \pr{ \ft_{\dm, \eps}^{\pr{\action}} \bullet \Xi_\dm^0 \pr{\ol{\omega}} } \in C\pr{\R^\dm, \R^\dm} $, and
    $ \pr[\big]{ \cR_\activation \pr{ \ft_{\dm, \eps}^{\pr{a}} \bullet \Xi_{\dm}^\theta \pr{\omega} } } \pr{x} = \pr[\big]{ \cR_\activation \pr{ \ft_{\dm, \eps}^{\pr{a}} } } \pr{ x, \xi_\dm^\theta \pr{\omega} } $.
    Hence, it holds for all $\theta \in \Theta$, $\dm \in \N$, $\err \in \eps$ that $ \max\cu{ \cL\pr{ \ft_{\dm, \eps}^{\pr{\action}} \bullet \Xi_\dm^0\pr{ \ol{\omega} } }, \fnorm{\cD \pr{ \ft_{\dm, \eps}^{\pr{\action}} \bullet \Xi_\dm^0\pr{ \ol{\omega} } }}  } \le \genConst \dm^\genConst \eps^{-\genConst} $.
    Moreover, observe that the assumption that for all $\theta \in \Theta$, $\dm \in \N$ it holds that $\F_\dm^\theta = \sigma \pr{\xi_\dm^\theta}$
    and the assumption that for all $\dm \in \N$, $\eps \in \err$, $a \in A$ it holds that $\cR_\activation \pr{\ft_{\dm, \eps}^{\pr{\action}}} \in C\pr{\R^{2\dm}, \R^\dm}$
    demonstrate for all $\theta \in \Theta$, $ \dm \in \N $, $\eps \in \err$, $ a \in A$ that $\R^\dm \times \Omega \ni \pr{x, \omega} \mapsto \pr[\big]{\cR_\activation \pr{ \ft_{\dm, \eps}^{\pr{a}} } }\pr{x, \xi_\dm^\theta \pr{\omega}} = \pr[\big]{ \cR_\activation \pr{ \ft_{\dm, \eps}^{\pr{a}} \bullet \Xi_\dm^\theta\pr{\omega} } } \pr{x} \in \R^\dm$ is $\pr{ \Borel\pr{\R^\dm} \otimes \F_\dm^\theta} / \Borel\pr{\R^\dm}$-measurable. 
    This implies for all $\theta \in \Theta$, $\dm \in \N$, $\eps \in \err$ that $ \pr[\big]{ \cR_\activation \pr{ \ft_{\dm, \eps}^{\pr{a}} \bullet \Xi_\dm^\theta } }_{a \in A} \colon \R^\dm \times \Omega \to \pr{\R^\dm}^{\abs{A}} $ is $ \pr{ \Borel\pr{\R^\dm} \otimes \F_\dm^\theta } / \Borel \pr{ \pr{ \R^\dm }^{\abs{A}} } $-measurable. 
    Next, note that the assumption that for all $\dm \in \N$ it holds that $\xi_\dm^\theta$, $\theta \in \Theta$, are i.i.d.\ random variables assures that for all $\dm \in \N$, $\eps \in \err$ it holds that $ \pr[\big]{ \cR_\activation \pr{ \ft_{\dm, \eps}^{\pr{a}} \bullet \Xi_\dm^\theta } }_{a \in A} $, $\theta \in \Theta$, are i.i.d.\ random fields.  
    Throughout the remainder of this proof let $\tker^{\pr{a}}_\dm \colon \R^\dm \times \Borel(\R^\dm) \to \br{0,1}$, $\dm \in \N$, $a \in A$, satisfy for all $\dm \in \N$, $x \in \R^\dm$, $a \in A$, $B \in \Borel(\R^\dm)$ that $\tker_\dm^{\pr{a}}\pr{ x, B } = \mathbb{P} \br{ t_\dm^{\pr{a}} \pr{x, \xi_\dm^0} \in B }$,
    let $\tkerANN_{\dm, \eps}^{\pr{a}} \colon \R^\dm \times \Borel(\R^\dm) \to \br{0,1}$, $\dm \in \N$, $\eps \in \err$, $a \in A$, satisfy for all $\dm \in \N$, $\eps \in \err$, $x \in \R^\dm$, $a \in A$, $B \in \Borel\pr{\R^\dm}$ that $\tkerANN_{\dm, \eps}^{\pr{a}} \pr{x, B} = \mathbb{P} \br{ \pr[\big]{ \cR_\activation\pr{ \ft_{\dm, \eps}^{\pr{a} }}} \pr{x, \xi_\dm^0}  \in B}$.
    Note that for all $\dm \in \N$, $\eps \in \err$, $a \in A$ it holds that $\tker_{\dm}^{\pr{a}}$ and $\tkerANN_{\dm, \eps}^{\pr{a}}$ are stochastic kernels.
    This and \cref{lem:properties_wdist} (applied for every $\dm \in \N$, $\eps \in \err$, $a \in A$ with $\dm \with \dm$,
    $m \with \dm$,
    $\eta \with \eta_\dm$,
    $\pr{\Omega, \cF, \mathbb{P}} \with \pr{\Omega, \cF, \mathbb{P}}$,
    $\xi \with \xi_\dm^0$,
    $\ft_1 \with t_\dm^{\pr{a}}$,
    $\ft_2 \with \cR_\activation \pr{\ft_{\dm, \eps}^{\pr{a}}}$
    in the notation of \cref{lem:properties_wdist})
    establishes that for all $\dm \in \N$, $x, y \in \R^\dm$, $\eps \in \err$, $a \in A$ it holds that
    $\wdist_\dm \pr{ \tker_{\dm}^{\pr{a}} \pr{x, \cdot}, \tker_{\dm}^{\pr{a}} \pr{y, \cdot} } \le \eta_\dm \norm{x - y}$ and 
    \begin{equation}
        \begin{split}
            \wdist_\dm \pr{ \tker_{\dm}^{\pr{a}} \pr{x, \cdot}, \tkerANN_{\dm, \eps}^{\pr{a}} \pr{x, \cdot} } \le \int_{\R^\dm} \norm{ t_\dm^{\pr{a}} \pr{x, p} - \pr[\big]{ \cR_\activation \pr{ \ft_{\dm, \eps}^{\pr{a}} } } \pr{x, p} } \pr[\big]{ \mathbb{P} \circ \pr{\xi_\dm^0}^{-1} } \pr{\dxx p}.
        \end{split}
    \end{equation}
    This and the assumption that for all $\dm \in \N$, $\eps \in \err$, $x, p \in \R^\dm$, $a \in A$ it holds that $\norm{ t_\dm^{\pr{a}} \pr{x, p} - \pr[\big]{ \cR_\activation \pr{\ft_{\dm, \eps}^{\pr{a}} } \pr{x, p} } } \le \eps \genConst \dm^\genConst \abs{\w_\dm \pr{x}}$ implies for all $\dm \in \N$, $x \in \R^\dm$, $\eps \in \err$, $a \in A$ that $\wdist_\dm \pr{ \tker_\dm^{\pr{a}} \pr{x, \cdot},  \tkerANN_{\dm, \eps}^{\pr{a}}\pr{x, \cdot} } \le \eps \genConst \dm^\genConst \abs{\w_\dm\pr{x}}$.
    Combining this 
    with   
    \cref{thm:main_simple} (applied with
    \begin{align}
        &\Theta \with \Theta, \quad A \with A, \quad \pr{\Omega, \cF, \mathbb{P}} \with \pr{\Omega, \cF, \mathbb{P}}, \quad \genConst \with \genConst, \quad \pr{ \lambda_\dm, \eta_\dm, K_\dm, L_\dm }_{\dm \in \N} \with \pr{ \lambda_\dm, \eta_\dm, K_\dm, L_\dm }_{\dm \in \N}  \nonumber\\
        &\pr{\lipFANN_{\dm, \eps}}_{\pr{\dm, \eps} \in \N \times \err} \with \pr{\lipFANN_{\dm, \eps}}_{\pr{\dm, \eps} \in \N \times \err}, \quad \pr{\w_\dm}_{\dm \in \N} \with \pr{\w_\dm}_{\dm \in \N}, \quad \pr{f_\dm}_{\dm \in \N} \with \pr{f_\dm}_{\dm \in \N}, \nonumber \\
        &\pr{\tker_\dm^{\pr{a}}}_{\pr{\dm, a} \in \N \times A} \with \pr{\tker_\dm^{\pr{a}}}_{\pr{\dm, a} \in \N \times A}, \quad \pr{\tkerANN_{\dm, \eps}^{\pr{a}}}_{\pr{\dm, a} \in \N \times \err \times A} \with \pr{\tkerANN_{\dm, \eps}^{\pr{a}}}_{\pr{\dm, \eps, a} \in \N \times \err \times A},  \nonumber \\
        & \action \with \action, \quad \ol{\omega} \with \ol{\omega}, \quad \activation \with \activation, \quad \fJ \with \fJ, \quad \pr{\bF_{\dm, \eps}}_{\pr{\dm, \eps} \in \N, \err} \with \pr{\bF_{\dm, \eps}}_{\pr{\dm, \eps} \in \N, \err}, \nonumber\\
        &\pr{ \bX_{\dm, \eps}^{ \theta, a } }_{ \pr{ \theta, \dm, \eps, a } \in \Theta \times \N \times \err \times A } \with \pr{ \ft_{\dm, \eps}^{\pr{a}} \bullet \Xi_\dm^\theta }_{ \pr{ \theta, \dm, \eps, a } \in \Theta \times \N \times \err \times A }, \nonumber \\
        &\pr{ \mathcal{F}_{\dm, \eps}^\theta }_{\pr{\theta, \dm, \eps} \in \Theta \times \N \times \err} \with \pr{ \F_{\dm}^\theta }_{\pr{ \theta, \dm, \eps } \in \Theta \times \N \times \err} 
    \end{align}
    in the notation of \cref{thm:main_simple})
    yields that there exists $c \in \R$ such that for every $\dm \in \N$, $\eps \in \err$ and every probability measure $\mu \colon \Borel\pr{\R^\dm} \to \br{0,1}$ with $\pr{\int_{\R^\dm} \abs{\w_\dm(x)}^2 \mu( \dxx x)}^{\nicefrac{1}{2}} \le \genConst d^\genConst$ it holds that
    \begin{enumerate}[label=(\Roman*)]
        \item \label{item:activation_transition_item_I_main_theorem} there exists a unique measurable $u_\dm \colon \R^\dm \to \R^{\abs{A}}$ which satisfies for all $x \in \R^\dm$, $a \in A$ that 
        $\sup_{(z, b) \in \R^\dm \times A} \abs{\w_\dm(z)}^{-1} \abs{(u_\dm(z))(b)} < \infty$,
        $\int_{\R^\dm} \abs{f_\dm(y, u_\dm(y))} \tker_\dm^{\pr{a}}(x, \dxx y) < \infty$, and
        \begin{equation}
            \pr*{ u_\dm (x) } \pr{a} = \int_{\R^\dm} f_\dm \pr{ y, u_\dm(y) } \tker^{\pr{a}}_\dm(x, \dxx y)
        \end{equation}
        and
        \item \label{item:activation_transition_item_II_main_theorem} there exists $\Psi_{\dm, \eps, \mu} \in \bN$ which satisfies $\cH\pr{\Psi_{\dm, \eps, \mu}} > 0$, $\paramANN(\Psi_{\dm, \eps, \mu}) \le c \dm^c \eps^{-c}$, $\cR_\activation(\Psi_{\dm, \eps, \mu}) \in C(\R^\dm, \R^{\abs{A}})$, and
        \begin{equation} 
            \pr*{ \int_{\R^\dm} \max_{a \in A} \abs*{ (u_\dm(x))(a) - \pr*{\pr[\big]{\cR_\activation(\Psi_{\dm, \eps, \mu})}(x)}(a) }^2 \mu(\dxx x) }^{\!\!\nicefrac{1}{2}} \le \eps.
        \end{equation}
    \end{enumerate}
    Note that \cref{item:activation_transition_item_I_main_theorem} shows for all $\dm \in \N$, $x \in \R^\dm$, $a \in A$ that $\E \br[\Big]{ \abs[\big]{ f_\dm \pr[\big]{ t_\dm^{\pr{a}} \pr{x, \xi_\dm^0}, u_\dm \pr{ t_\dm^{\pr{a}} \pr{x, \xi_\dm^0} } } } } < \infty$ and 
    \begin{equation}
        \pr{u_\dm\pr{x}}\pr{a} = \E \br*{ f_\dm \pr[\big]{ t_\dm^{\pr{a}} \pr{x, \xi_\dm^0}, u_\dm \pr{ t_\dm^{\pr{a}} \pr{x, \xi_\dm^0} } } }.
    \end{equation}
    This and the assumption that for all $\dm \in \N$ it holds that $\xi_\dm$ and $\xi_\dm^0$ are identically distributed prove \cref{item:activation_transition_ex_uniq_sol}.
    Moreover, observe that \cref{item:activation_transition_item_II_main_theorem} establishes \cref{item:activation_transition_ex_ann}.
    The proof of \cref{cor:activation_transition} is thus complete.
\end{mproof}

\subsection{ANN approximations for Bellman equations with specific activation functions and one-step transitions}
\label{subsec:bellman}

\cfclear
\begin{corollary}
    \label{cor:bellman}
    \cfconsiderloaded{cor:bellman}
    Let $A$ be a nonempty finite set,
    let $(\Omega, \mathcal{F}, \mathbb{P})$ be a probability space,
    let $\beta \in [0, \infty) \backslash \cu{1}$,
    let $ \lambda, \eta, K, L \in \pr{0,\infty}$,
    assume $L \max\cu{ \lambda, \eta } < 1$,
    let $\genConst \in [ \max\cu{1, \lambda, K} , \infty)$,
    for every $\dm \in \N$, $x \in \pr{x_1, \dots, x_\dm} \in \R^\dm$ let $\fnorm{x} \in \R$ satisfy $\fnorm{x} = \max_{i \in \cu{1,\dots, \dm}} \abs{x_i}$, 
    for every $\dm \in \N$ let $\w_\dm \colon \R^\dm \to \pr{0, \infty}$ and $g_\dm \colon \R^\dm \to \R^{\abs{A}}$ be measurable,
    for every $\dm \in \N$, $a \in A$ let $t^{\pr{a}}_\dm \colon \R^{2\dm} \to \R^\dm$ be measurable,
    for every $\dm \in \N$ let $\xi_\dm \colon \Omega \to \R^\dm$ be measurable,
    let $\action \in A$, 
    let $\activation \colon \R \to \R$ satisfy for all $x \in \R$ that $\activation(x) = \max\cu{x, \beta x} $,
    let $(\bG_{\dm, \eps})_{(\dm, \eps) \in \N \times \err} \subseteq \bN$ satisfy for all $\dm \in \N$, $\eps \in \N$ that $\cR_\activation(\bG_{\dm, \eps}) \in C(\R^{\dm}, \R^{\abs{A}})$,
    let $(\ft_{\dm, \eps}^{\pr{a}})_{(\dm, \eps, a) \in \N \times \err \times A} \subseteq \bN$ satisfy for all $\dm \in \N$, $\eps \in \err$, $a \in A$ that $\cD\pr{ \ft_{\dm, \eps}^{\pr{a}}} = \cD \pr{ \ft_{\dm, \eps}^{\pr{\action}}}$ and $\cR_\activation\pr{\ft_{\dm, \eps}^{\pr{a}}} \in C(\R^{2\dm}, \R^\dm)$,
    assume for all $\dm \in \N$, $\eps \in \err$ that 
    \begin{equation}
        \begin{split}
            \max\cu*{\cL\pr{ \bG_{\dm, \eps} }, \cL\pr{ \ft_{\dm, \eps}^{\pr{\action}} }, \fnorm[\big]{\cD\pr{ \bG_{\dm, \eps} }}, \fnorm[\big]{ \cD\pr{ \ft_{\dm, \eps}^{\pr{\action}} } } } \le \genConst \dm^\genConst \eps^{- \genConst},
        \end{split}
    \end{equation}
    assume for all $\dm \in \N$, $\eps \in \err$, $x,y, p \in \R^\dm$, $a \in A$ that 
    \begin{align}
        \begin{split} \label{eq:bellman_approximation_assumption}
            &\abs[\big]{ \pr[\big]{\pr[\big]{ \cR_\activation (\bG_{\dm, \eps}) }(x)}\pr{a} -  \pr{g_\dm (x)}\pr{a} } \le \eps \genConst \dm^{\genConst} \abs{\w_\dm(x)},\\
            &\norm[\big]{ \pr[\big]{\cR_\activation \pr{\ft_{\dm, \eps}^{\pr{a}}}} \pr{x, p} - t_\dm^{\pr{a}}\pr{x, p} } \le \eps \genConst \dm^{\genConst} \abs{\w_\dm(x)},
        \end{split} \\
        \begin{split} \label{eq:bellman_lipschitz}
            & \abs*{ \pr{g_\dm \pr{x}} \pr{a} - \pr{g_\dm\pr{y}}\pr{a}} \le K \norm{ x - y}, \quad \text{and} \\
            &\norm[\big]{ t_\dm^{\pr{a}}\pr{ x, p } - t_\dm^{\pr{a}}\pr{ y, p }  }  \le \eta \norm{x - y},
        \end{split}
    \end{align}
    assume for all $\dm \in \N$, $\eps \in \err$, $a \in A$ that
    \begin{equation}
        \begin{split}
            \sup_{x \in \R^\dm} \abs{\w_\dm\pr{x}}^{-1} \E \br[\big]{ \norm{t_\dm^{\pr{a}} \pr{x, \xi_\dm}} } < \infty \quad \text{ and } \quad \sup_{x \in \R^\dm } \abs{\w_\dm\pr{x}}^{-1} \E \br[\big]{ \norm[\big]{ \pr[\big]{ \cR_\activation \pr{ \ft_{\dm, \eps}^{\pr{a}} } } \pr{x, \xi_\dm} }  } < \infty,
        \end{split}
    \end{equation}
    assume for all $\dm \in \N$, $x \in \R^\dm$, $\eps \in \err$, $a \in A$ that
    \begin{align}
        \begin{split}
            &\pr*{ \E\br*{ \abs*{ \w_\dm \pr*{ \pr[\big]{ \cR_\activation \pr{ \ft_{\dm, \eps}^{\pr{a}}} } \pr{x, \xi_\dm} } }^2 } }^{\!\!\nicefrac{1}{2}} \le \lambda \abs{\w_\dm\pr{x}}, \\
            &\pr*{ \E \br*{ \max_{b \in A} \abs*{ \pr[\Big]{\pr[\big]{ \cR_\activation \pr{\bG_{\dm, \eps}} } \pr[\big]{ \pr[\big]{ \cR_\activation \pr{\ft_{\dm, \eps}^{\pr{a}}  } } \pr{x, \xi_\dm} }} \pr{b}  }^2 } }^{\!\!\nicefrac{1}{2}} \le \genConst \dm^\genConst \abs{\w_\dm\pr{x}}, 
        \end{split}\\
        \begin{split}
            & \E \br[\Big]{ \abs[\big]{ \w_\dm \pr[\big]{ t_\dm^{\pr{a}}(x, \xi_\dm) } } } \le \lambda \abs{\w_\dm\pr{x}}, \quad \text{ and } \quad \sup_{\substack{y \in \R^\dm \\ b \in A}} \abs{\w_\dm\pr{y}}^{-1} \abs{ \pr{g_\dm\pr{y}} \pr{b}} < \infty.
        \end{split}
    \end{align}
    Then there exists $c \in \R$ such that for every $\dm \in \N$, $\eps \in \err$ and every probability measure $\mu \colon \Borel\pr{\R^\dm} \to \br{0,1}$ with $\pr{\int_{\R^\dm} \abs{\w_\dm(x)}^2 \mu( \dxx x)}^{\nicefrac{1}{2}} \le \genConst d^\genConst$ it holds that
    \begin{enumerate}[label = (\roman*)]
        \item \label{item:bellman_ex_uniq_sol} there exists a unique measurable $\qfunc \colon \R^\dm \to \R^{\abs{A}}$ which satisfies for all $x \in \R^\dm$, $a \in A$ that 
        $\sup_{(z, b) \in \R^\dm \times A} \abs{\w_\dm(z)}^{-1} \abs{(\qfunc(z))(b)} < \infty$,
        $\E \br[\Big]{ \max_{b \in A} \abs[\big]{ \pr{ \qfunc \pr{ t_{\dm}^{\pr{a}}\pr{x, \xi_\dm} } } \pr{b} }  } < \infty$, and
        \begin{equation}
            \begin{split}
                \pr{\qfunc\pr{x}}\pr{a} = \pr{g_\dm(x)}\pr{a} + L \E \br*{ \max_{b \in A}  \pr{ \qfunc \pr{ t_{\dm}^{\pr{a}}\pr{x, \xi_\dm} } } \pr{b} }
            \end{split}
        \end{equation}
        and
        
        \item \label{item:bellman_ex_ann} there exists $\bQ \in \bN$ which satisfies $\cH \pr{\bQ} > 0$, $\paramANN (\bQ) \le c \dm^c \eps^{-c}$, $\cR_\activation (\bQ) \in C( \R^\dm, \R^{\abs{A}})$, and
        \begin{equation} \label{eq:bellman_approximating_ann}
            \pr*{ \int_{\R^\dm} \max_{a \in A} \abs*{ (\qfunc(x))(a) - \pr*{ \pr*{ \cR_\activation(\bQ) }(x) }(a) }^2 \mu(\dxx x) }^{\!\!\nicefrac{1}{2}} \le \eps
        \end{equation}
    \end{enumerate}
    \cfout.
\end{corollary}

\begin{mproof}{\cref{cor:bellman}}
    Throughout this proof let $\genConstVar = \max\cu{L, 2\abs{A} + 1} \genConst$, let $\bI, \fJ \in \bN$ satisfy
    \begin{align}
        \bI &= \pr*{ \pr*{ \begin{pmatrix}
                1 \\ - 1
            \end{pmatrix}, \begin{pmatrix}
                0 \\ 0
            \end{pmatrix} }, \pr[\Big]{ \begin{pmatrix}
                1 & -1
            \end{pmatrix}, 0 }} \in \pr{\R^{2\times 1} \times \R^2} \times \pr{\R^{1 \times 2} \times \R^1} \qquad \text{and}\\
        \fJ &= (1 + \beta)^{-1} \circledast \bI
    \end{align}
    let $\I_k \in \bN$, $k \in \N$, satisfy for all $k \in \N$ that $\I_k = \paraANN{k}(\fJ, \fJ, \dots, \fJ)$,
    for every $\dm \in \N$ let $f_\dm \colon \R^\dm \times \R^{|A|} \to \R$ satisfy for all $x \in \R^\dm$, $r\in \R^{|A|}$ that $f_\dm (x, r) = L \max_{a \in A} \{ ( g_\dm(x) )(a) + r(a) \}$.
    Observe that \cite[Lemma~3.5]{ackermann2023deep} ensures for all $x \in \R$ that $\cD\pr{\fJ} = \pr{1,2,1} \in \N^3$ and $\pr{ \cR_\activation \pr{\fJ} }\pr{x} = x$.
    This and \cref{lem:parallelizations_of_anns_same_depth} imply for all $k \in \N$, $ x \in \R^k $ that $\cD \pr{\I_k} = \pr{k, 2k, k} \in \N^3$ and $\pr{\cR_\activation \pr{\I_k} } \pr{x} = x$.
    Next, note that \cref{lem:max_of_finite_set_leaky} ensures that there exists $\bM \in \bN$ such that for all $x = (x_1, x_2, \dots, x_{|A|}) \in \R^{|A|}$ it holds that $\cR_\activation(\bM) \in C(\R^{|A|}, \R)$, $ \big( \cR_\activation (\bM)  \big)(x) = \max_{i \in [1, |A|] \cap \N} x_i$, $\cL(\bM) = \ceil{\log_2(|A|)} + 1$, and $\fnorm{ \cD (\bM) } \le 2|A|$.
    Throughout the remainder of this proof let $(\bF_{\dm, \eps})_{(\dm, \eps) \in \N \times \err} \subseteq \bN$ satisfy for all $\dm \in \N$, $\eps \in \err$ that
    \begin{equation} \label{eq:bellman_def_nonlin_net}
        \bF_{\dm, \eps} = L \circledast \pr[\big]{\bM \bullet \fS_{|A|, 2} \bullet \paraLANN{2}{\pr{\I_{\abs{A}}, \I_{\abs{A}}}} \pr[\big]{ \bG_{\dm, \eps}, \I_{|A|}} } .
    \end{equation}
    Next, observe that
    the assumption that $A$ is nonempty,
    the assumption that for all $\dm \in \N$, $\eps \in \err$ it holds that $ \max\cu{ \cL \pr{ \bG_{\dm, \eps} }, \fnorm{\cD \pr{ \bG_{\dm, \eps} }}} \le \genConst \dm^\genConst \eps^{-\genConst} $,
    the assumption that for all $r = (r(a))_{a \in A} \in \R^{\abs{A}}$ it holds that $\pr[\big]{ \cR_\activation (\bM) }(r) = \max_{a \in A} r(a)$,
    the fact that for all $k \in \N$ it holds that $\cD (\I_k) = (k , 2 k, k)$,
    the fact that for all $n \in \N$ it holds that $\ceil{\log_2(n)} + 1 \le n$,
    \cref{lem:compositions_of_anns},
    \cref{lem:extensions_of_anns_involving_identities},
    \cref{cor:parallelizations_of_anns_with_different_layer_structure_realization_involving_identities}, and
    \cref{lem:scalar_multiplications_of_anns} 
    assure for all $\dm \in \N$, $\eps \in \err$, $x \in \R^\dm$, $r = (r(a))_{a \in A} \in \R^{\abs{A}}$ that
    \begin{align}
        \begin{split} \label{eq:bellman_nonlin_architecture}
            \cI (\bF_{\dm, \eps}) &= \cI \pr[\big]{ \paraLANN{2}{\pr{\I_{\abs{A}}, \I_{\abs{A}}}} \pr[\big]{ \bG_{\dm, \eps}, \I_{|A|}}  } = \dm + |A|,\\
            \cO (\bF_{\dm, \eps}) &= \cO (\bM) = 1, \qquad \cR_\activation (\bF_{\dm, \eps}) \in C(\R^{\dm + |A|}, \R),\\
            \cL(\bF_{\dm, \eps}) &= \cL(\bM) + \cL\pr{ \fS_{\abs{A}, 2} } + \cL\pr{ \paraLANN{2}{\pr{\I_{\abs{A}}, \I_{\abs{A}}}} \pr[\big]{ \bG_{\dm, \eps}, \I_{|A|}} } - 2\\ 
            &\le \ceil{ \log_2(\abs{A}) } + 1 + \max\cu{ \cL\pr{ \bG_{\dm, \eps} }, \cL \pr{ \I_{\abs{A}}} } - 1 \le \abs{A} + \cL(\bG_{\dm, \eps}) + 1 \\
            &\le \pr{\abs{A} + 2} \genConst \dm^\genConst \eps^{-\genConst} \le \pr{2\abs{A} + 1} \genConst \dm^\genConst \eps^{-\genConst} \le \genConstVar \dm^\genConstVar \eps^{-\genConstVar},\\
            \fnorm{\cD(\bF_{\dm, \eps})} &\le \max \cu[\big]{ \fnorm{\cD(\bM)}, \fnorm{\cD (\fS_{\abs{A}, 2})}, \fnorm{ \cD \pr{ \paraLANN{2}{\pr{\I_{\abs{A}}, \I_{\abs{A}}}} \pr{ \bG_{\dm, \eps}, \I_{|A|}} } } } \\
            &\le \max \cu[\big]{ \fnorm{\cD(\bM)},  2\abs{A}, 2\abs{A} + \fnorm{\cD(\bG_{\dm, \eps})} }  \\
            &= 2\abs{A} + \fnorm{\cD(\bG_{\dm, \eps})} \le \pr{2\abs{A} + 1} \genConst \dm^\genConst \eps^{-\genConst} \le \genConstVar \dm^\genConstVar \eps^{- \genConstVar}, \qquad \qquad \text{and} 
        \end{split} \\
        \begin{split} \label{eq:bellman_nonlin_realization} 
            \pr*{\cR_\activation (\bF_{\dm, \eps})}(x,r) &= L \pr*{ \pr[\big]{\cR_\activation (\bM) \circ \cR_\activation (\fS_{|A|, 2}) \circ \cR_\activation \pr[\big]{ \paraLANN{2}{\pr{\I_{\abs{A}}, \I_{\abs{A}}}} \pr{ \bG_{\dm, \eps}, \I_{|A|}} }} (x, r) } \\
            &= L \pr[\bigg]{ \pr[\big]{\cR_\activation (\bM) \circ \cR_\activation(\fS_{\dm, \eps}) } \begin{pmatrix}
                \pr[\big]{\cR_\activation(\bG_{\dm, \eps})} (x) \\r
            \end{pmatrix} } \\
            &= L \pr[\big]{\pr[\big]{\cR_\activation(\bM)}  \pr[\big]{(\cR_\activation (\bG_{\dm, \eps}))(x) + r} }  \\
            &= L \max_{a \in A}  \cu[\big]{  \pr[\big]{\pr[\big]{ \cR_\activation(\bG_{\dm, \eps})}(x)}(a) + r(a) }.
        \end{split}
    \end{align}
    Next, observe that for all $\dm \in \N$, $x \in \R^\dm$, $r,s \in \R^{\abs{A}}$, $\eps \in \err$ it holds that $\abs{ f_\dm\pr{x, r} - f_\dm \pr{x, s} } \le L \max_{a \in A} \abs{r\pr{a} - s\pr{a}}$ and $\abs{ \pr{ \cR_\activation \pr{\bF_{\dm, \eps}} } \pr{x, r} - \pr{ \cR_\activation \pr{\bF_{\dm, \eps}} } \pr{x, s} } \le L \max_{a \in A} \abs{ r\pr{a} - s\pr{a} }$.
    Moreover, note that the assumption that for all $\dm \in \N$, $x,y\in \R^\dm$, $a \in A$ it holds that $\abs{ \pr{g_\dm\pr{x}}\pr{a} - \pr{g_\dm\pr{y}}\pr{a} } \le K \norm{x - y}$
    implies for all $\dm \in \N$, $x,y \in \R^\dm$, $r\in \R^{\abs{A}}$ that
    \begin{equation}
        \begin{split} \label{eq:bellman_nonlin_lipschitz_1st_arg}
            \abs{ f_\dm \pr{x,r} - f_\dm \pr{y, r} } \le L \max_{a \in A} \abs{ \pr{g_\dm \pr{x}}\pr{a} - \pr{g_\dm \pr{y}}\pr{a} } \le K L \norm{x - y}.
        \end{split}
    \end{equation}
    In addition, observe that the assumption that for all $\dm \in \N$, $\eps \in \err$, $x \in \R^\dm$, $a \in A$ it holds that
    $\abs[\big]{ \pr[\big]{\pr[\big]{ \cR_\activation (\bG_{\dm, \eps}) }(x)}\pr{a} -  \pr{g_\dm (x)}\pr{a} } \le \eps \genConst \dm^{\genConst} \abs{\w_\dm(x)}$
    establishes for all $\dm \in \N$, $\eps \in \err$, $x \in \R^\dm$, $r \in \R^\dm$ that
    \begin{equation}
        \begin{split}
            \abs[\big]{ \pr{ \cR_\activation \pr{\bF_{\dm, \eps}} } \pr{x, r} - f_\dm \pr{x, r} } &\le L \max_{a \in A} \abs[\big]{ \pr[\big]{ \pr{\cR_\activation \pr{\bG_{\dm, \eps}}} \pr{x} } \pr{a} - \pr{ g_\dm\pr{x} }\pr{a}  }\\
            &\le \eps \genConst L \dm^\genConst \abs{\w_\dm \pr{x}} \le \eps \genConstVar \dm^\genConstVar \abs{\w_\dm\pr{x}}.
        \end{split}
    \end{equation} 
    This
    and \cref{cor:activation_transition} (applied with
    \begin{equation}
        \begin{split}
            &A \with A, \quad \pr{\Omega, \mathcal{F}, \mathbb{P}} \with \pr{\Omega, \mathcal{F}, \mathbb{P}}, \quad \genConst \with \genConstVar, \quad \nu \with 1, \quad \beta \with \beta, \\
            &\pr{\lambda_\dm, \eta_\dm, K_\dm, L_\dm}_{\dm \in \N} \with (\lambda, \eta, KL, L)_{\dm \in \N}, \quad \pr{\lipFANN_{\dm, \eps}}_{\pr{\dm, \eps} \in \N \times \err} \with \pr{ L}_{\pr{\dm, \eps} \in \N \times \err}, \\
            &\pr{\w_\dm}_{\dm \in \N} \with \pr{\w_\dm}_{\dm \in \N}, \quad \pr{f_\dm}_{\dm \in \N} \with \pr{f_\dm}_{\dm \in \N}, \\
            &\pr{ t_\dm^{\pr{a}} }_{\pr{\dm, a } \in \N \times A} \with \pr{ t_\dm^{\pr{a}} }_{\pr{\dm, a } \in \N \times A}, \quad \pr{\xi_\dm}_{\dm \in \N} \with \pr{\xi_\dm}_{\dm \in \N}, \quad \action \with \action, \quad \activation \with \activation, \\
            &\pr{\bF_{\dm, \eps}}_{\pr{\dm, \eps} \in \N \times \err} \with \pr{\bF_{\dm, \eps}}_{\pr{\dm, \eps} \in \N \times \err}, \quad (\ft_{\dm, \eps}^{\pr{a}})_{(\dm, \eps, a) \in \N \times \err \times A} \with (\ft_{\dm, \eps}^{\pr{a}})_{(\dm, \eps, a) \in \N \times \err \times A}
        \end{split}
    \end{equation}
    in the notation of \cref{cor:activation_transition})
    demonstrate that there exists $c \in \R$ such that for every $\dm \in \N$, $\eps \in \err$ and every probability measure $\mu \colon \Borel\pr{\R^\dm} \to \br{0,1}$ with $\pr{\int_{\R^\dm} \abs{\w_\dm(x)}^2 \mu( \dxx x)}^{\nicefrac{1}{2}} \le \genConst d^\genConst$ it holds that
    \begin{enumerate}[label = (\Roman*)]
        \item \label{item:bellman_item_I_cor_activation_transition} there exists a unique measurable $u_\dm \colon \R^\dm \to \R^{\abs{A}}$ which satisfies for all $x \in \R^\dm$, $a \in A$ that 
        $\sup_{(z, b) \in \R^\dm \times A} \abs{\w_\dm(z)}^{-1} \abs{(u_\dm(z))(b)} < \infty$,
        \begin{equation}
            \E \br[\Big]{  \max_{b \in A} \abs[\big]{ \pr[\big]{ g_\dm\pr{ t_\dm^{\pr{a}} \pr{x, \xi_\dm} } }\pr{b} + \pr[\big]{u_\dm \pr{ t_\dm^{\pr{a}} \pr{x, \xi_\dm} }} \pr{b} } } < \infty,
        \end{equation}
        and
        \begin{equation}
            \pr{u_\dm\pr{x}}\pr{a} = L \E \br[\Big]{ \max_{b \in A} \pr[\Big]{ \pr[\big]{ g_\dm\pr{ t_\dm^{\pr{a}} \pr{x, \xi_\dm} } }\pr{b} + \pr[\big]{u_\dm \pr{ t_\dm^{\pr{a}} \pr{x, \xi_\dm} }} \pr{b} } } 
        \end{equation}
        and
        
        \item \label{item:bellman_item_II_cor_activation_transition} there exists $\Psi_{\dm, \eps, \mu} \in \bN$ which satisfies $ \cH\pr{\Psi_{\dm, \eps, \mu}} > 0$, $\paramANN(\Psi_{\dm, \eps, \mu}) \le c \dm^c \eps^{-c}$, $\cR_\activation(\Psi_{\dm, \eps, \mu}) \in C(\R^\dm, \R^{\abs{A}})$, and
        \begin{equation} \label{eq:bellman_approximating_ann_cor_activation_transition}
            \pr*{ \int_{\R^\dm} \max_{a \in A} \abs*{ (u_\dm(x))(a) - \pr*{\pr[\big]{\cR_\activation(\Psi_{\dm, \eps, \mu})}(x)}(a) }^2 \mu(\dxx x) }^{\!\!\nicefrac{1}{2}} \le \eps.
        \end{equation}
    \end{enumerate}
    Throughout the remainder of this proof let $\qfunc_\dm \colon \R^\dm \to \R^{\abs{A}}$, $\dm \in \N$, satisfy for all $\dm \in \N$, $x \in \R^\dm$, $a \in A$ that $\pr{\qfunc_\dm \pr{x}}\pr{a} = \pr{g_\dm\pr{x}}\pr{a} + \pr{u_\dm\pr{x}}\pr{x}$,
    let $\pr{\delta_{\dm, \eps}}_{\pr{\dm, \eps} \in \N \times \err} \subseteq \R$ satisfy for all $\dm \in \N$, $\eps \in \err$ that $\delta_{\dm, \eps} = \eps \pr{1 + \genConst^2 \dm^{2\genConst}}^{-1}$,
    let $\fK = \max\cu{2\abs{A}, \genConst, c}$,
    and for every $\dm \in \N$, $\eps \in \err$ and every probability measure $\mu \colon \Borel\pr{\R^\dm} \to \br{0,1}$ with $\pr{\int_{\R^\dm} \abs{\w_\dm(x)}^2 \mu( \dxx x)}^{\nicefrac{1}{2}} \le \genConst d^\genConst$ let $\bQ_{\dm, \eps, \mu} \in \bN$
    satisfy that $\bQ_{\dm, \eps, \mu} = \bG_{\dm, \delta_{\dm, \eps}} \,\bSum_{\, \I_{|A|}}\, \Psi_{\dm, \delta_{\dm, \eps}, \mu}$.
    Observe that for all $\dm \in \N$, $\eps \in \err$ it holds that $\delta_{\dm, \eps} \in (0,1]$.
    Moreover, it follows from 
    the assumption that for all $\dm \in \N$ it holds that $g_\dm$ is measurable,
    the assumption that for all $\dm \in \N$ it holds that $\sup_{(x, a) \in \R^\dm \times A} \abs{\w_\dm(x)}^{-1} \abs{(g_\dm(x))(a)} < \infty$,
    and \cref{item:bellman_item_I_cor_activation_transition} that for all $\dm \in \N$, $x \in \R^\dm$, $a \in A$ it holds that 
    $\qfunc_\dm$ is measurable,
    $\sup_{\pr{y,b} \in \R^\dm \times A} \abs{\w_\dm\pr{y}}^{-1} \abs{\pr{\qfunc_\dm\pr{y}}\pr{b}} < \infty$,
    $\E \br[\Big]{ \max_{b \in A} \abs[\big]{ \pr{ \qfunc_\dm \pr{ t_{\dm}^{\pr{a}}\pr{x, \xi_\dm} } } \pr{b} }  } < \infty$, and
    \begin{equation}
        \begin{split}
            \pr{\qfunc_\dm\pr{x}}\pr{a} &= \pr{g_\dm\pr{x}}\pr{a} + \pr{u_\dm\pr{x}}\pr{x} \\
            &= \pr{g_\dm\pr{x}}\pr{a} + L \E \br[\Big]{ \max_{b \in A} \cu[\big]{ \pr[\big]{ g_\dm\pr{ t_\dm^{\pr{a}} \pr{x, \xi_\dm} } }\pr{b} + \pr[\big]{u_\dm \pr{ t_\dm^{\pr{a}} \pr{x, \xi_\dm} }} \pr{b} } } \\
            &= \pr{g_\dm\pr{x}}\pr{a} + L \E \br[\Big]{ \max_{b \in A} \cu[\big]{ \pr[\big]{ \qfunc_\dm\pr{ t_\dm^{\pr{a}} \pr{x, \xi_\dm} } }\pr{b} } }.
        \end{split}
    \end{equation}
    This proves \cref{item:bellman_ex_uniq_sol}.
    Furthermore, note that
    \cref{lem:linear_combination_anns_different_length}
    establishes for every $\dm \in \N$, $\eps \in \err$, $x \in \R^\dm$, $a \in A$ and every probability measure $\mu \colon \Borel\pr{\R^\dm} \to \br{0,1}$ with $\pr{\int_{\R^\dm} \abs{\w_\dm(x)}^2 \mu( \dxx x)}^{\nicefrac{1}{2}} \le \genConst d^\genConst$ that
    $ \cI(\bQ_{\dm, \eps, \mu}) = \cI ( \bG_{\dm, \delta_{\dm, \eps}}) = \dm$, $\cO(\bQ_{\dm, \eps, \mu}) = \cO(\bG_{\dm, \delta_{\dm, \eps}}) = \abs{A}$, $\cR_\activation (\bQ_{\dm, \eps, \mu}) \in C(\R^\dm, \R^{\abs{A}})$, and $\pr*{ \cR_\activation (\bQ_{\dm, \eps, \mu})} (x) = \pr*{ \cR_\activation(\bG_{\dm, \delta_{\dm, \eps}}) }(x) + \pr*{ \cR_\activation(\Psi_{\dm, \delta_{\dm, \eps}, \mu}) }(x) \in \R^{\abs{A}}$.
    This,
    the assumption that for all $\dm \in \N$, $\eps \in \err$, $x \in \R^\dm$, $a \in A$ it holds that $\abs{ (g_\dm(x))(a) - \pr*{\pr*{ \cR_\activation(\bG_{\dm, \eps}) }(x)}(a) } \le \eps \genConst \dm^\genConst \abs{\w_\dm(x)}$,
    the triangle inequality,
    and \cref{item:bellman_item_II_cor_activation_transition}
    yield for every $\dm \in \N$, $\eps \in \err$ and every probability measure $\mu \colon \Borel\pr{\R^\dm} \to \br{0,1}$ with $\pr{\int_{\R^\dm} \abs{\w_\dm(x)}^2 \mu( \dxx x)}^{\nicefrac{1}{2}} \le \genConst d^\genConst$ that
    \begin{align}
        &\pr*{\int_{\R^\dm} \max_{a \in A} \abs*{ (\qfunc_\dm(x))(a) -  \pr[\big]{ \pr*{ \cR_\activation (\bQ_{\dm, \eps, \mu}) }(x) }(a) }^2   \mu (\dxx x)}^{\!\!\nicefrac{1}{2}} \nonumber \\
        &= \pr[\bigg]{ \int_{\R^\dm} \!\!  \max_{a \in A} \abs[\big]{ (g_\dm(x))(a) -  \pr*{\pr*{\cR_\activation(\bG_{\dm, \delta_{\dm, \eps}})}(x)}(a) + (u_\dm(x))(a) - \pr*{\pr*{\cR_\activation(\Psi_{\dm, \delta_{\dm, \eps}, \mu})}(x)}(a) }^2  \mu(\dxx x)\! }^{\!\!\nicefrac{1}{2}} \nonumber\\
        &\le \pr*{\int_{\R^\dm} \max_{a \in A} \abs*{ (g_\dm(x))(a) -  \pr*{\pr*{\cR_\activation(\bG_{\dm, \delta_{\dm, \eps}})}(x)}(a) }^2 \mu(\dxx x)}^{\!\!\nicefrac{1}{2}} \nonumber \\
        &\quad+ \pr*{\int_{\R^\dm} \max_{a \in A} \abs*{ (u_\dm(x))(a) - \pr*{\pr*{\cR_\activation(\Psi_{\dm, \delta_{\dm, \eps}, \mu})}(x)}(a) }^2 \mu(\dxx x)}^{\!\!\nicefrac{1}{2}} \nonumber \\
        &\le \delta_{\dm, \eps} \genConst\dm^\genConst \pr*{ \int_{\R^\dm} \abs{\w_\dm(x)}^2 \mu (\dxx x)}^{\!\! \nicefrac{1}{2}} + \delta_{\dm, \eps} \le \delta_{\dm, \eps} (1 + \genConst^2 \dm^{2\genConst}  ) \le \eps.
    \end{align}
    This proves \eqref{eq:bellman_approximating_ann}. 
    Moreover, observe that the assumption that $\fK = \max\cu{2\abs{A}, \genConst, c}$,
    the assumption that for all $ \dm \in \N $, $\eps \in \err$ it holds that $\max\cu{ \cL\pr{\bG_{\dm, \eps}}, \fnorm{\cD \pr{ \bG_{\dm, \eps} }} } \le \genConst \dm^\genConst \eps^{-\genConst}$,
    \cref{item:bellman_item_II_cor_activation_transition},
    and \cref{lem:ANN_size_estimate}
    establish for every $\dm \in \N$, $\eps \in \err$ and every probability measure $\mu \colon \Borel\pr{\R^\dm} \to \br{0,1}$ with $\pr{\int_{\R^\dm} \abs{\w_\dm(x)}^2 \mu( \dxx x)}^{\nicefrac{1}{2}} \le \genConst d^\genConst$ that
    \begin{equation}
        \begin{split}
            \cL\pr{ \bQ_{\dm, \eps, \mu} } &\le \max \cu{ \cL \pr{ \bG_{\dm, \delta_{\dm, \eps}} }, \cL\pr{ \Psi_{\dm, \delta_{\dm, \eps}, \mu} } } \le \max\cu{ \genConst \dm^\genConst \delta_{\dm, \eps}^{-\genConst}, c\dm^c \delta_{\dm, \eps}^{-c} } \le \fK \dm^\fK \delta_{\dm, \eps}^{-\fK}.
        \end{split}
    \end{equation}
    and 
    \begin{equation}
        \begin{split}
            \fnorm{\cD (\bQ_{\dm, \eps, \mu})} &\le \fnorm{\cD( \longerANN{ \cL( \bQ_{\dm, \eps, \mu}), \I_{\abs{A}}} ( \bG_{\dm, \delta_{\dm, \eps}} )  )} + \fnorm{ \cD( \longerANN{\cL(\bQ_{\dm, \eps, \mu}), \I_{\abs{A}}} (\Psi_{\dm, \delta_{\dm, \eps}, \mu}) ) } \\
            &\le \max\cu[\big]{ \fnorm{\cD(\I_{\abs{A}}) }, \fnorm{\cD (\bG_{\dm, \delta_{\dm, \eps}} )} } + \max\cu[\big]{ \fnorm{\cD(\I_{\abs{A}}) }, \fnorm{\cD (\Psi_{\dm, \delta_{\dm, \eps}, \mu} )} } \\
            &\le \max \cu{ 2\abs{A}, \genConst \dm^\genConst \delta_{\dm, \eps}^{- \genConst}} + \max \cu{ 2 \abs{A}, c\dm^c \delta_{\dm, \eps}^{- c }} \le \fK \dm^\fK \delta_{\dm, \eps}^{- \fK}.
        \end{split}
    \end{equation}
    This,
    the assumption that for all $\dm \in \N$, $\eps \in \err$ it holds that $\delta_{\dm, \eps} = \eps \pr{1 + \genConst^2 \dm^{2 \genConst}}^{-1}$,
    and \cref{lem:ANN_size_estimate}
    demonstrate for every $\dm \in \N$, $\eps \in \err$ and every probability measure $\mu \colon \Borel\pr{\R^\dm} \to \br{0,1}$ with $\pr{\int_{\R^\dm} \abs{\w_\dm(x)}^2 \mu( \dxx x)}^{\nicefrac{1}{2}} \le \genConst d^\genConst$ that
    \begin{equation}
        \begin{split}
            \paramANN \pr{ \bQ_{\dm, \eps, \mu} } &\le 2 \cL \pr{\bQ_{\dm, \eps, \mu}} \fnorm{\cD \pr{ \bQ_{\dm, \eps, \mu} }}^2 \\
            &\le 2 \fK \dm^\fK \delta_{\dm, \eps}^{- \fK} \pr{ 2 \fK \dm^\fK \delta_{\dm, \eps}^{- \fK} }^2 \\
            &= 8 \fK^3 \dm^{3\fK} \pr{ 1+ \genConst^2 \dm^{2\genConst}}^{3 \fK} \eps^{- 3\fK} \\
            &\le 8 \fK^3 \dm^{3\fK} 2^{3\fK} \genConst^{6\fK} \dm^{6\genConst \fK} \eps^{- 3\fK}\\
            &= 8^{\fK + 1} \fK^3 \genConst^{6\fK} \dm^{3\fK \pr{ 1 + 2\genConst }} \eps^{-3\fK}.
        \end{split}
    \end{equation}
    This proves \cref{item:bellman_ex_ann}.
    The proof of \cref{cor:bellman} is thus complete.
\end{mproof}

\subsubsection*{Acknowledgements}
This work has been partially funded by the Deutsche Forschungsgemeinschaft (DFG, German \mbox{Research} Foundation) in the frame of the priority programme SPP 2298 ``Theoretical \mbox{Foundations} of Deep Learning'' – Project Number 464101154. Moreover, this work has been partially funded by the European Union (ERC, MONTECARLO, 101045811). The views and the opinions expressed in this work are however those of the authors only and do not necessarily reflect those of the European Union or the European Research Council (ERC). Neither the \mbox{European} Union nor the granting authority can be held responsible for them. We also gratefully acknowledge the \mbox{Cluster} of Excellence EXC 2044-390685587, Mathematics Münster: Dynamics-Geometry-Structure funded by the Deutsche Forschungsgemeinschaft (DFG, German Research \mbox{Foundation}).

{
\bibliographystyle{acm}
\bibliography{bibfile}
}

\end{document}